\documentclass[a4paper,11pt]{article}
\usepackage{url}
\usepackage{amssymb}
\usepackage{amsmath}
\usepackage{amsthm}
\usepackage{geometry}
\usepackage{stmaryrd}
\usepackage{siunitx}
\usepackage{commath}
\usepackage{algorithm}
\usepackage{algorithmic}
\usepackage{subfig}
\usepackage{enumerate}
\usepackage{amsmath}
\usepackage{multirow}
\usepackage{makecell}
\usepackage{lipsum}
\usepackage{amsfonts}
\usepackage{graphicx}
\usepackage{epstopdf}
\usepackage{algorithm}
\usepackage{algorithmic}
\usepackage{amssymb}
\usepackage{bm}
\usepackage{color}
\usepackage{caption}
\usepackage{subfig}

\usepackage{amsmath}
\usepackage{multirow}
\usepackage{makecell}

\usepackage{amsfonts}
\usepackage{graphicx}
\usepackage{epstopdf}

\usepackage{amssymb}
\usepackage{bm}
\usepackage{color}

\usepackage{caption}
\usepackage{subfig}

\usepackage[hidelinks]{hyperref}
\hypersetup{
colorlinks   = true, 
urlcolor     = red, 
linkcolor    = red, 
citecolor   = red 
}

\usepackage[misc]{ifsym}
\usepackage{color}
\usepackage{cleveref}

\theoremstyle{definition}
\newtheorem{definition}{Definition}[section]

\newtheorem{remark}{Remark}[section]

\makeatletter
\def\ps@pprintTitle{%
\let\@oddhead\@empty
\let\@evenhead\@empty
\def\@oddfoot{}%
\let\@evenfoot\@oddfoot}

\newcommand{\tnorm}{\@ifstar\@tnorms\@tnorm}
\newcommand{\@tnorm}[2][]{%
\mathopen{#1|\mkern-1.5mu#1|\mkern-1.5mu#1|}
#2
\mathclose{#1|\mkern-1.5mu#1|\mkern-1.5mu#1|}
}

\makeatother

\newtheoremstyle{mystyle}
{3pt} 
{1pt} 
{} 
{} 
{\bfseries} 
{.} 
{.5em} 
{} 

\theoremstyle{mystyle} 
\newtheorem{theorem}{Theorem}[section]
\newtheorem{lemma}{Lemma}[section]

\allowdisplaybreaks[4]
\DeclareMathOperator{\diag}{diag}

\newtheorem{assumption}[theorem]{Assumption}

\theoremstyle{remark}

\numberwithin{equation}{section}

\usepackage{float}
\usepackage{comment}

\makeatletter

\@addtoreset{equation}{section}

\makeatother

\begin{document}
%
\title{Stochastic gradient descent based variational inference for infinite-dimensional inverse problems}

\author{Jiaming Sui
\footnote{Email: sjming1997327@stu.xjtu.edu.cn, School of Mathematics and Statistics, Xi'an Jiaotong University, Xi'an, Shaanxi 710049, China},
Junxiong Jia 
\footnote{Email: jjx323@xjtu.edu.cn, School of Mathematics and Statistics, Xi'an Jiaotong University, Xi'an, Shaanxi 710049, China}, 
Jinglai Li
\footnote{Email: j.li.10@bham.ac.uk, School of Mathematics, University of Birmingham, Birmingham, B15 2TT, United Kingdom}}




\maketitle
\begin{abstract}

This paper introduces two variational inference approaches for infinite-dimensional inverse problems, developed through gradient descent with a constant learning rate. The proposed methods enable efficient approximate sampling from the target posterior distribution using a constant-rate stochastic gradient descent (cSGD) iteration. Specifically, we introduce a randomization strategy that incorporates stochastic gradient noise, allowing the cSGD iteration to be viewed as a discrete-time process. This transformation establishes key relationships between the covariance operators of the approximate and true posterior distributions, thereby validating cSGD as a variational inference method. We also investigate the regularization properties of the cSGD iteration and provide a theoretical analysis of the discretization error between the approximated posterior mean and the true background function. Building on this framework, we develop a preconditioned version of cSGD to further improve sampling efficiency. Finally, we apply the proposed methods to two practical inverse problems: one governed by a simple smooth equation and the other by the steady-state Darcy flow equation. Numerical results confirm our theoretical findings and compare the sampling performance of the two approaches for solving linear and non-linear inverse problems.

~\\ 
\textbf{keyword}: inverse problems, infinite-dimensional variational inference, Bayesian analysis for functions, partial differential equations, stochastic gradient descent

\end{abstract}


\section{Introduction}

Due to its multidisciplinary applications in areas like seismic exploration \cite{weglein2003inverse} and medical imaging \cite{zhou2020bayesian}, the study of inverse problems involving partial differential equations (PDEs) has seen significant advancements in recent decades \cite{arridge2019solving}. 
When addressing inverse problems for PDEs, uncertainties are common, such as measurement uncertainty and epistemic uncertainty.
The Bayesian inverse approach offers a flexible framework for solving these problems by converting them into statistical inference tasks, allowing for the analysis of uncertainties in the solutions.

Since inverse problems for PDEs are often defined in infinite-dimensional spaces, traditional finite-dimensional Bayesian methods \cite{kaipio2006statistical}, which are well understood, cannot be directly applied to these infinite-dimensional cases \cite{Stuart2010ActaNumerica}. 
To overcome this challenge, there are two different approaches wildly employed:
\begin{itemize}
	\item \emph{Discretize-then-Bayesianize}: The PDEs are first discretized to approximate the original problem in a finite-dimensional space. The resulting simplified problem is then solved using finite-dimensional Bayesian methods \cite{kaipio2006statistical}.
	\item \emph{Bayesianize-then-discretize}: In this approach, Bayesian formulas and algorithms are first constructed in infinite-dimensional spaces. Some infinite-dimensional algorithms are developed or built, and then finite-dimensional approximations are applied \cite{Stuart2010ActaNumerica}.
\end{itemize}
Both approaches have their advantages and disadvantages.
Using the first approach, we can apply various Bayesian inference methods established in the statistical literature to solve inverse problems.
However, since the original problems are defined in infinite-dimensional spaces, certain issues arise, such as discretization errors \cite{Stuart2010ActaNumerica} and non-uniform convergence issues \cite{Stuart2010ActaNumerica}, posing significant challenges to this approach.
On the other hand, the Bayesianize-then-discretize approach offers several advantages: 
first an understanding of the structures within function spaces is crucial for designing optimal numerical schemes for PDEs, particularly regarding the gradient information discussed in \cite{hinze2008optimization};
second, rigorously formulating infinite-dimensional theory can help avoid misleading intuitive notions that may arise from finite-dimensional inverse methods, such as the total variation prior \cite{Lassas2004IP,yao2016tv}.
Taking these advantages, Bayesianize-then-discretize approach has attracted considerable attention in recent years \cite{bui2013computational,Cotter2009IP,Stuart2010ActaNumerica}. 

A central task of the Bayesian approach is to compute the posterior distributions.
For many practical problems, the posterior distribution is not feasible to compute directly, necessitating the use of an approximated posterior measure. 
A popular strategy to compute the posterior distribution is Markov chain Monte Carlo (MCMC) \cite{cotter2013}. 
From the \emph{Bayesianize-then-discretize} perspective, several infinite-dimensional MCMC  methods have been proposed, where a notable example is 
the preconditioned Crank-Nicolson (pCN) \cite{Pillai2014SPDE}. 
Furthermore, sampling methods that utilize gradient, sample history and geometrical structure have been developed, such as the infinite-dimensional Metropolis-adjusted Langevin algorithm \cite{Thanh2016IPI}, the adaptive pCN logarithm \cite{hu2017adaptive, zhou2017hybrid}, and the geometric pCN algorithm \cite{Beskos2017JCP}.
However, the computational cost for MCMC methods becomes prohibitive for many applications, making it challenging to apply them to large-scale problems \cite{Fichtner2011Book}.

Variational inference (VI) \cite{blei2017variational, zhang2018advances}, which aims to estimate an approximation of the actual posterior,   
offers an efficient alternative to MCMC.
Some studies on VI methods for inverse problems related to PDEs have adopted a \emph{discretize-then-Bayesianize} approach in finite-dimensional spaces. 
For example, mean-field assumption based VI approach was utilized to solve finite-dimensional inverse problems involving hyperparameters in prior and noise distributions \cite{Guha2015JCP, jin2010hierarchical}. 
However, there has been less research on VI in the context of infinite-dimensional settings for solving inverse problems related to PDEs from the \emph{Bayesianize-then-discretize} perspective.
Specifically, when the approximated measures are constrained to be Gaussian, a Robbins-Monro algorithm was developed from a calculus-of-variations perspective in the works of \cite{pinski2015algorithms, Pinski2015SIAMMA}. 
To accommodate non-Gaussian approximated measures, a comprehensive infinite-dimensional mean-field variational inference  theoretical framework was established in \cite{jia2021variational}, which was further extended to hierarchical inverse problems in \cite{Sui2024MOC}.
Additionally, building on this theoretical framework, a generative deep neural network model, referred to as VINet, was constructed in\cite{Jia2022VINet}, 
and an infinite-dimensional Stein variational gradient descent approach was proposed in \cite{jia2021stein}.

In this paper, we consider a different line of VI methods, which are motivated by the  Stochastic gradient descent (SGD) optimization algorithm.
SGD facilitates efficient optimization in the presence of large datasets, as stochastic gradients can often be computed inexpensively through random sampling.
Remarkably, it was shown in \cite{Stephan2016PMLR, Stephan2017JMLR} that SGD with a constant learning rate can be posed as a hybrid method that combines Monte Carlo and variational algorithms.
Instead of computing an analytical form of the approximate posterior,
this approach samples from an asymptotically approximated posterior measure, offering minimal implementation effort while typically achieving faster mixing times.
Recently, an infinite-dimensional SGD method for linear inverse problems was proposed in \cite{lu2022moc}, which provides 
an estimation error analysis taking into account the discretization levels and learning rate. This method seeks a point estimate of the unknown parameter, rather than
computing the posterior distribution. 
To our knowledge, there have been few studies on the SGD based VI methods in infinite-dimensional spaces. 
Therefore, we seek to develop an infinite-dimensional variational inference (iVI) method, extended from the constant SGD (cSGD) method proposed in \cite{Stephan2016PMLR, Stephan2017JMLR}. 
In this paper, we present our method within the linear inverse problem setting (but do not require its posterior to be Gaussian), as some of our theoretical results are developed in this context. 
Nevertheless, as demonstrated by numerical results, our proposed method can also be applied to nonlinear problems.
We summarize the main technical contributions of the work as the following: 
\begin{itemize}
	\item 
	By introducing stochastic gradient noise, we develop a randomization strategy for the gradient that differs from the approach in \cite{lu2022moc}. 
	Consequently, the cSGD iteration is approximated by an infinite-dimensional discrete-time process, whose stationary probability measures are chosen to estimate the posterior measure.
	We obtain a theoretical result characterizing the covariance operator of the estimated posterior measure, which provides the theoretical foundation for applying cSGD-iVI in infinite-dimensional spaces.
	\item We determine the optimal learning rate of cSGD by minimizing the Kullback-Leibler divergence between the estimated and real posterior measures. 
	Taking consideration into discretization, we also prove the regularization properties \cite{Peter2003IP1} based on the cSGD formulation. 
	Subsequently, we obtain a discretization error bound between samples and the background truth function, which is controlled by the learning rate and the discretization levels.
	\item To further improve the sampling efficiency,  we develop a  preconditioned version of cSGD (termed as pcSGD-iVI).
	Similarly, theoretical results on the relationship between the estimated and real posteriors are also provided,
    and the optimal learning rate is  derived accordingly.
	Moreover, we discuss various choices of the stochastic gradient noise.
	\item { Moreover, we consider the loss function by applying the projection operator.
	And we calculate the corresponding measure of the stochastic gradient to determine prior operator $\mathcal{Q}$.
	Then we bound scale \( S \) through the eigenvalues of \( \mathcal{Q} \), which guarantees the stable and well-defined solution of discrete-time Lyapunov equation.}
	
\end{itemize}

The outline of this paper is listed as follows. 
In Section $\ref{sec2}$, we develop the cSGD-iVI and pcSGD-iVI methods in order to sample from the estimated posterior measures of linear inverse problems.
{ In Subsection $\ref{subsec2.1}$, we introduce the inverse problem, and SGD method on Hilbert spaces.
In Subsection $\ref{subsec2.2}$, we construct the general framework of cSGD-iVI approach.
In Subsection $\ref{subsec2.3}$, we establish pcSGD-iVI approach based on the framework proposed in $\ref{subsec2.2}$, and discuss the strategies of determining the stochastic gradient noise.}
In Section \ref{sec3}, we apply the developed cSGD-iVI and pcSGD-iVI approaches to two inverse problems governed by a simple smooth equation and the steady-state Darcy flow equation, with both noise and prior measures assumed to be Gaussian.
Section \ref{sec4} concludes with a summary of our findings, highlights certain limitations, and proposes potential directions for future research.

\section{Infinite dimensional constant SGD VI method }\label{sec2}

In this section, we present a variational inference framework based on constant-step stochastic gradient descent (cSGD) for generating samples from an estimated posterior distribution.
By appropriately controlling the strength of the stochastic gradient noise, the stationary distribution of the SGD iterates can be characterized through a discrete-time Lyapunov equation. 
To approximate the posterior distribution, we minimize the Kullback–Leibler (KL) divergence between the estimated and true posterior measures, allowing samples to be generated via SGD updates with an optimal learning rate. 
Our method extends the theoretical framework introduced in \cite{Stephan2016PMLR, Stephan2017JMLR} to infinite-dimensional settings. 
Furthermore, we derive the mean function and covariance operator of the stationary distribution by treating the stochastic gradient as a random variable, which enables us to identify the covariance structure of the gradient noise and determine its scaling factor. 
Based on this analysis, we propose an iterative algorithm. 
In addition, we introduce the pcSGD-based infinite-dimensional variational inference (iVI) method, which builds upon the cSGD framework. 
The theoretical proofs for the results in this section are provided in Section C of the supplementary material.

\subsection{Problem setup: Inverse problems in a Hilbert space}\label{subsec2.1}
In this section, the following two goals need to be achieved: introduce the infinite-dimensional Bayesian inverse problem; and provide the loss function under the current Bayesian formulation.
We consider a generic inverse problem defined on a Hilbert space $\mathcal{H}_u$: 
\begin{align}\label{eq:invproblem}
	\bm{d} = H(u) + \bm{\epsilon},
\end{align}
where $\bm{d}\in \mathbb{R}^{N_d}$ is the measurement data, $u\in \mathcal{H}_u$ is the parameter of interest, 
$H: \mathcal{H}_u \rightarrow \mathbb{R}^{N_d}$ is the forward operator.
The measurement noise $\bm{\epsilon}$ is a Gaussian random vector defined on $\mathbb{R}^{N_d}$ with zero mean and variance $\bm{\Gamma}_{\text{noise}}$, indicating
$\bm{\epsilon} \sim \mathcal{N}(0, \bm{\Gamma}_{\text{noise}})$.

We adopt a Bayesian framework for solving this inverse problem. 
First we assume that the prior distribution of $u$ is 
$\mathcal{N}(0, \mathcal{C}_0)$, where $\mathcal{C}_0 : \mathcal{H}_u \rightarrow \mathcal{H}_u$ is a self-adjoint, positive-definite and trace-class operator.
Let $\lbrace c_i, e_i \rbrace^\infty_{i = 1}$ be the eigensystem of operator $\mathcal{C}_0$ such that $\mathcal{C}_0e_i = c_ie_i$ for $i = 1, 2, ...$.
And without loss of generality, we assume that the eigenmodes $\lbrace e_i \rbrace^\infty_{i = 1}$ are orthonormal and the eigenvalues $\lbrace c_i\ \rbrace^\infty_{i = 1}$ are in descending order.
Finally note that $\lbrace e_i \rbrace^\infty_{i = 1}$ form
a complete basis of $\mathcal{H}_u$, which means
any $u \in \mathcal{H}_u$ can be linearly represented by $\{e_i\}_{i=1}^\infty$. 

The likelihood function associated with the inverse problem is expressed as:
\begin{align*}
	\pi(\bm{d}|u) \propto \exp(-\Phi(u,\bm{d})),
\end{align*}
where the potential $ \Phi(u, \bm{d}) =\frac{1}{2}\| H(u)-\bm{d}\|^2_{\bm{\Gamma}_{\text{noise}}} $.
Based on \cite{dashti2013bayesian}, the following theorem establishes the existence and explicit form of the posterior measure in the infinite-dimensional Bayesian inverse problem setting.
\begin{theorem}\label{BayesTheorem}
	Let $ \mathcal{H}_u $ be a separable Hilbert space and $ N_d $ a positive integer. Suppose $ \mu_0 $ is a Gaussian measure on $ \mathcal{H}_u $ and for some constants $ M_1, M_2 > 0 $, $ \Phi : \mathcal{H}_u \rightarrow \mathbb{R} $ satisfies:
	\begin{align*}
		\Phi(u) \geq -M_1, \quad |\Phi(u, \bm{d}_1) - \Phi(u, \bm{d}_2)| \leq M_2 \|\bm{d}_1 - \bm{d}_2\|_2.
	\end{align*}
	Under these conditions, the posterior measure $ \mu \ll \mu_0 $ on $ \mathcal{H}_u $ is well-defined with Radon-Nikodym derivative:
	\begin{align}\label{eq:post}
		\frac{d\mu}{d\mu_0}(u) = \frac{1}{Z_{\mu}} \exp(-\Phi(u)),
	\end{align}
	where the normalization constant $ Z_{\mu} $ is given by:
	\begin{align*}
		Z_{\mu} = \int_{\mathcal{H}_u} \exp(-\Phi(u)) \mu_0(du).
	\end{align*}
\end{theorem}
Throughout this paper, we assume that the likelihood function satisfies the conditions of Theorem \ref{BayesTheorem},
and so the posterior distribution exists.

Next, we introduce the cost function, defined by:
\begin{align}\label{eq:loss}
\mathcal{L}(u) = \Phi(u) + {R}(u),
\end{align}
where 
\begin{equation}
{R}(u) = \frac12\|\mathcal{C}_0^{-1/2}u\|^2_{\mathcal{H}_u}, 
\end{equation}
 is the regularization term derived from the Gaussian prior measure.
 Here, $\| \cdot \|_{\mathcal{H}_u}$ denoting the usual Hilbert norm defined on $\mathcal{H}_u$.
 For notational convenience, we denote $\| \cdot \|_{\mathcal{H}_u}$ by $\| \cdot \|$.
 The minimizer \( \bar{u} \) of the cost function \eqref{eq:loss} corresponds to the maximum a posteriori (MAP) estimate of the posterior distribution on \( \mathcal{H}_u \) \cite{bui2013computational}.

Consider a sample \( u\in \mathcal{H}_u \) drawn from the prior measure, which can be expanded as \( u = \sum_{i=1}^\infty u_ie_i \).
Let us consider the natural isomorphism \( \gamma \) between \(\mathcal{H}_u\) and the Hilbert space \( l^2 \) of all sequences \( (u_i)_{i=1}^\infty \) of real numbers such that
\[
\sum_{i=1}^\infty | u_i |^2 < \infty,
\]
defined by
\[
H \rightarrow l^2, \quad u \mapsto \gamma(u) = (u_i)_{i=1}^\infty.
\]
Then we indentify \( \mathcal{H}_u \) with \( l^2 \), and we consider the following product measure
\[
\tilde{\mu}_0 := \bigotimes_{i=1}^\infty \mathcal{N}(0, c_i),
\]
where \( \{c_i \}_{i=1}^\infty \) is the eigenvalues of covariance \( \mathcal{C}_0 \).
This means, for a sample \(u\) from the prior measure, we have \( u_i \sim \mathcal{N}(0, c_i) \).
Furthermore, let us assume that operator \( H \) is linear operator, then we rewrite \( H \) on \( l^2 \) as a matrix given by
\begin{align}\label{eq:H}
H =
\left(
\begin{array}{cccc}
h_{1,1} & h_{1,2} & \cdots & h_{1,N_d} \\
h_{2,1} & h_{2,2} & \cdots & h_{2,N_d} \\
\vdots & \vdots & \ddots & \vdots \\
h_{M,1} & h_{M,2} & \cdots & h_{M,N_d} \\
0 & 0 & \cdots & 0 \\
\vdots & \vdots & ~ & \vdots
\end{array}
\right),
\end{align}
Recalling the eigensystem \( \{ c_i, e_i \}_{i=1}^\infty \), prior covariance operator \( \mathcal{C}_0 \) on \( l^2 \) as a martix is rewritten by
\begin{align}\label{eq:C_0}
\mathcal{C}_0 = \text{diag}(c_1, c_2, \cdots, c_M, \cdots).
\end{align}
For notational convenience, we will use the same symbol for the operator in spaces \( \mathcal{H}_u \) and \( l^2 \).
For example, we use the symbol \( H \) to denote the forward operator \( H \) on \( \mathcal{H}_u \) as shown in Eq.~\eqref{eq:invproblem}, and matrix on \( l^2 \) as shown in Eq.~\eqref{eq:H}.

\subsection{Variational inference based on SGD method}\label{subsec2.2}
In this section, we present the theoretical framework of the constant SGD-based iVI (cSGD-iVI) approach, which enables sampling from the estimated posterior distribution by executing the cSGD iteration. 
This method extends the approach proposed in \cite{Stephan2016PMLR, Stephan2017JMLR} for finite-dimensional problems. 
In the cSGD iteration, the estimated posterior is determined by two key components: the learning rate and the stochastic gradient, as defined in equation $(\ref{eq:SGD})$. 
The choice of these components critically influences how well the estimated posterior approximates the true posterior distribution of 
$u$.
We first introduce an alternative randomization strategy for generating stochastic gradients (SG), which incorporates stochastic gradient noise, following the formulation in \cite{Stephan2016PMLR, Stephan2017JMLR} for finite-dimensional settings. 
Once the stochastic gradient is defined, we then turn to the learning rate. In the context of Bayesian inference, the optimal learning rate is obtained by minimizing the Kullback–Leibler (KL) divergence between the estimated and true posterior distributions.
Our goals in this section are listed as follows:
\begin{itemize}
	\item propose the random strategy to generate SG from the full gradient, and clarify the SGD iteration form based on SG;
	\item calculate the mean function and covariance operator of the estimated posterior;
	\item determine the optimal learning rate using a Bayesian inference approach; 
	\item present the algorithm for the cSGD-iVI approach. 
\end{itemize}

\subsubsection{The SGD iteration in Hilbert Spaces}\label{subsec:2.2.1}
We begin by introducing the stochastic gradient descent (SGD) iteration.
When the goal is to compute the MAP estimator, the function $\mathcal{L}(u)$ can be minimized using the standard gradient descent scheme:
\begin{align}\label{e:gd}
	u_{k+1} = u_k - \eta_k{ \mathcal{G}}(u_k),
\end{align}
where $\mathcal{G}(u) = \nabla \mathcal{L}(u)$ and $\eta_k$ is the step size (a.k.a.,  the learning rate) at step $k$.
A stochastic version of this update takes the form:
\begin{align}\label{eq:SGD}
	u_{k+1} = u_k - \eta_k\widetilde{ \mathcal{G}}(u_k),
\end{align}
where $\widetilde{ \mathcal{G}}(u_k)$ is the stochastic gradient, which serves as a randomized approximation of the full gradient $\mathcal{G}(u_k)$. 

In \cite{Stephan2016PMLR, Stephan2017JMLR}, the stochastic gradient is constructed using mini-batches of data points. In our setting, we will define the stochastic gradient differently; for now, we assume it takes the following mathematical form:
\begin{definition}\label{assump1}
	Let $\mathcal{Q}$ be an operator defined on $\mathcal{H}_u$ as, 
	\begin{align*}
		\mathcal{Q} = \sum^M_{i = 1}q_ie_i\otimes e_i + \sum^\infty_{i = M+1}e_i\otimes e_i,
	\end{align*}
	and let $\mathcal{C}_\text{GN}$ be,
	\begin{align*}
		\mathcal{C}_\text{GN} = \mathcal{C}_0\mathcal{Q} = \sum^M_{i = 1}c_iq_ie_i\otimes e_i + \sum^\infty_{i = M+1}c_ie_i\otimes e_i.
	\end{align*}
	Then we define the \emph{stochastic gradient} $\widetilde{\mathcal{G}}(u)$ as
	\begin{align}\label{eq:SG}
		\widetilde{\mathcal{G}}_S(u) = \mathcal{G}(u) - \frac{1}{\sqrt{S}}\Delta \mathcal{G}(u), 
	\end{align}
	where $\Delta \mathcal{G}(u) \sim \mu_{\text{GN}}=\mathcal{N}(0, \mathcal{C}_\text{GN})$ is referred to as \emph{gradient noise} (GN),
	and $S$ is a positive constant controlling the strength of the gradient noise. 
\end{definition}
We remark that the parameter $S$ is an analogy of the batch size of the stochastic gradient in  \cite{Stephan2016PMLR, Stephan2017JMLR}. 
It is straightforward to verify that the stochastic gradient $\widetilde{\mathcal{G}}_S(u)$ is an unbiased estimator of the full gradient $\mathcal{G}(u)$. 
Substituting this stochastic gradient into the iteration scheme Eq.~($\ref{eq:SG}$), we rewrite the SGD iteration form ($\ref{eq:SGD}$) with a constant learning rate $\eta$ as, 
\begin{equation} \label{eq:SGDiter}
{ u_{(k+1)} = u_{k} - \eta\widetilde{ \mathcal{G}}_S(u_{k}).}
\end{equation}

{ \begin{remark}
	In Eq.~($\ref{eq:SGDiter}$), we introduce the noisy full-gradient iteration rather than a mini-batch scheme.
	Mini-batching is for finite-sum functions \(\sum_{i=1}^{N}\ell_i(u)\), but the cost function \(\mathcal{L}(u) = \Phi(u) + R(u)\) does not admit such a large-sample finite-sum structure.
	On the other hand, we need to solve one forward and one adjoint PDE to calculate gradient \(\mathcal{G}(u)\) for each iteration step.
	Thus the computational cost is unaffected by parameter \(S\).
\end{remark}}

In what follows address the following two issues. 
First, we want to show that iteration~\eqref{eq:SGDiter} provides an approximate sampling scheme for the posterior distribution.
Second we want to using a variational formulation to determine the optimal value of $S$.

\subsubsection{Assumptions of SGD} \label{subsec:2.2.a}


A critical step in the Bayesian method for solving the inverse problems is to compute the posterior distribution. 
In \cite{Stephan2016PMLR, Stephan2017JMLR}, a method based on SGD iterations was proposed for approximate posterior sampling in finite-dimensional settings. 
In this section, we extend this approach to the infinite-dimensional case and present a set of assumptions, adapted from those in the aforementioned works.

\begin{assumption}\label{a:SG}
The stochastic gradient $\widetilde{ \mathcal{G}}_S$  in  \eqref{eq:SGDiter} is defined according to Definition \ref{assump1}.
\end{assumption}

\begin{assumption} \label{a:quad}
We assume that the stationary distribution of the iterates is supported in a region $\Xi$, within which the Onsager–Machlup functional can be well approximated by a quadratic expansion (up to a constant term):
\begin{equation}\label{eq:gloss}
 \mathcal{L}(u) =  \frac{1}{2}\|\mathcal{A}^{1/2}(u - \bar{u})\|^2 + R(u-\bar{u}),
\end{equation}
where $ \bar{u} $ is the MAP estimator and $ \mathcal{A} $ is the Hessian operator of the data fidelity term $\Phi$ at $\bar{u}$. 
\end{assumption}

\begin{assumption}\label{a:H}
There exist a positive integer $M$ such that
\[
H(u) = H \bigg(\sum^{\infty}_{i = 1}u_ie_i \bigg) \approx H \bigg(u^M :=\sum^{M}_{i = 1}u_ie_i \bigg),
\]
and we assume that the forward operator \( H \) is linear throughout this section.
\end{assumption}

\begin{remark} Some remarks on the above assumptions are provided as follows:
\begin{itemize}
\item Assumption \ref{a:SG} corresponds to Assumptions 1–2 in \cite{Stephan2016PMLR}, adapted to the infinite-dimensional setting. 
\item Assumption \ref{a:quad} is the infinite-dimensional analogue of Assumption 4 in \cite{Stephan2016PMLR}, which allows the posterior to be locally approximated by a Gaussian distribution.
\item Assumption 3 in \cite{Stephan2016PMLR}, which involves a continuous-time approximation, is not adopted here; instead, we work directly with the discrete-time formulation.
\item The truncation level \( M \) in { Assumption 2.4} can be determined heuristically \cite{feng2018adaptive} by
\[
M = \min \bigg\{ M \in \mathbb{N} \bigg|~ \frac{c_M}{c_1} < C_M \bigg\},
\]
where \( C_M \) is a prescribed threshold.
This type of truncation is standard in infinite-dimensional inverse problems \cite{Arendt2020IMA, Cui2014IP, CUI2016JCP, Cui2021IP}, and reflects the fact that finite-resolution observations only inform a finite number of eigenmodes.
\end{itemize}
\end{remark}

Under Assumptions \ref{a:SG}, \ref{a:quad}, \ref{a:H}, we decompose the parameter \( u \in \mathcal{H}_u \) into two orthogonal components:
\[
u=u^M+u^{\perp},
\]
where \( u^M = \sum^{M}_{i = 1}u_ie_i \) represents the projection onto the first \( M \) eigenmodes of the prior covariance operator, and $u^{\perp} = \sum_{i=M+1}^\infty u_i e_i$ lies in the orthogonal complement.
It follows that the cost function becomes
\begin{equation}
\mathcal{L}(u)= \Phi(u^M) + R(u^M) + R(u^{\perp}). 
\end{equation}
This decomposition implies that the posterior measure can be factored into a product of two independent components: one associated with the data-informed subspace and the other with the prior-only complement. 
Then we can calculate the posterior measure of \( u^M, u^\perp \) separately since they are independent.
Before we provide the detailed approach, one natural question arises: whether the product posterior measure of \( u^M\) and \( u^\perp \) equal to the posterior measure \(\mu \) of \( u \)?

Using the isomorphism \( \gamma \) in Subsection \ref{subsec2.1}, we conduct the analysis in \( l^2 \).
To address the question, we examine both the mean function and the covariance operator.
For the covariance operator, let us consider subspace \( u^M \in \mathcal{H}^M_u \), the inverse problem is then given by Eq.\eqref{eq:invproblem}.
The prior measure on \( \mathcal{H}^M_u \) of \( u^M \) is denoted by
\[
\mu_0^M = \mathcal{N}(0, \mathcal{C}_0^M),
\]
where in \( l^2 \),
\begin{align*}
\mathcal{C}_0^M = \diag(c_1, \cdots, c_M).
\end{align*}
According to Proposition 3.1 in \cite{Knapik2011TAS}, the posterior covariance operator for \( u^M \) is given by
\[
\mathcal{C}^M = \mathcal{C}_0^M - \mathcal{C}_0^MH^{M, *}(\bm{\Gamma}_{\text{noise}} + H^M\mathcal{C}_0^MH^{M, *})^{-1}H^M\mathcal{C}_0^M,
\]
where
\begin{align*}
H^M =
\left(
\begin{array}{cccc}
h_{1,1} & h_{1,2} & \cdots & h_{1,N_d} \\
h_{2,1} & h_{2,2} & \cdots & h_{2,N_d} \\
\vdots & \vdots & \ddots & \vdots \\
h_{M,1} & h_{M,2} & \cdots & h_{M,N_d}
\end{array}
\right).
\end{align*}
On the orthogonal complement \( u^\perp \in \mathcal{H}^\perp_u := \text{span}(e_{M+1}, e_{M+2}, \cdots ) \), the prior measure is
\[
\mu^\perp_0 = \mathcal{N}(0, \mathcal{C}^\perp_0),
\]
where in \( l^2 \),
\begin{align*}
\mathcal{C}_0^\perp = \diag(c_{M+1}, c_{M+2}, \cdots,).
\end{align*}
And the posterior covariance remains unchanged:
\[
\mathcal{C}^\perp = \mathcal{C}_0^\perp.
\]
Then the product posterior covariance of \( u^M, u^\perp \) on \( l^2 \) is given by
\begin{align*}
\mathcal{C}^\prime =
\left(
\begin{array}{cc}
\mathcal{C}^M & 0 \\
0 & \mathcal{C}_0^\perp
\end{array}
\right)=
\left(
\begin{array}{cc}
\mathcal{C}_0^M - \mathcal{C}_0^MH^{M, *}(\bm{\Gamma}_{\text{noise}} + H^M\mathcal{C}_0^MH^{M, *})^{-1}H^M\mathcal{C}_0^M & 0 \\
0 & \mathcal{C}_0^\perp
\end{array}
\right)
\end{align*}
Alternatively, the posterior covariance from the full-space Bayesian update is:
\begin{align}\label{eq:C}
    \mathcal{C} = \mathcal{C}_0 - \mathcal{C}_0H^{*}(\bm{\Gamma}_{\text{noise}} + H\mathcal{C}_0H^{*})^{-1}H\mathcal{C}_0
\end{align}
with operators \( \mathcal{C}_0, H \) being defined in Eqs.~\eqref{eq:C_0} and \eqref{eq:H}.
Substituting \( \mathcal{C}_0^M, H^M \) into Eq.~\eqref{eq:C} yields
\[
\mathcal{C} = 
\left(
\begin{array}{cc}
\mathcal{C}_0^M & 0\\
0 & \mathcal{C}_0^\perp
\end{array}
\right) +
\left(
\begin{array}{cc}
-\mathcal{C}_0^MH^{M, *}(\bm{\Gamma}_{\text{noise}} + H^M\mathcal{C}_0^MH^{M, *})^{-1}H^M\mathcal{C}_0^M & 0\\
0 & 0
\end{array}
\right)
= \mathcal{C}^\prime.
\]

We now consider the mean function.
Minimizing the loss function \( \mathcal{L}(u)= \Phi(u^M) + R(u^M) + R(u^{\perp}) \) yields the posterior mean \( \bar{u} \).
Separately, on subspace \( \mathcal{H}_u^M \), the posterior mean \( \bar{u}^M \) corresponding to \( u^M \) is derived from the cost function:
\[
\mathcal{L}^M(u^M)= \Phi(u^M) + R(u^M).
\]
And on space \( \mathcal{H}_u^\perp \), the posterior mean \( \bar{u}^\perp \) corresponding to \( u^\perp \) is derived from the cost function:
\[
\mathcal{L}^\perp(u^M)= R(u^\perp).
\]
Then the product posterior mean of \( u^M, u^\perp \) on \( l^2 \) is given by
\[
\bar{u}^\prime =
\left(
\begin{array}{c}
\bar{u}^M \\
\bar{u}^\perp
\end{array}
\right) = \bar{u}.
\]

Combining the results for the mean and covariance operators, we conclude that the full posterior measure \( \mu \) coincides with the product measure of the posteriors of \( u^M \) and \( u^\perp \) on \( l^2 \).
Moreover, since the posterior distribution of \( u^\perp \) is identical to its prior, the posterior distribution of \( u \) can be effectively approximated by independently approximating the posteriors of \( u^M \) and \( u^\perp \).
This decomposition serves as the foundation for the sampling framework introduced in Subsection \ref{subsec:2.2.2}.


\subsubsection{SGD as a discrete-time process} \label{subsec:2.2.2}
In this section, we demonstrate that the SGD iteration \eqref{eq:SGDiter} is a sampling scheme for an approximate posterior distribution.
Based on Subsection \ref{subsec:2.2.a}, we can then approximate the posterior measure of \( u \) by these two components: \( u^M \) and \( u^\perp \), individually.
Following the variational approximation framework \cite{blei2017variational}, we aim to approximate the true posterior measure with a Gaussian distribution, hereafter referred to as the \textit{estimated posterior} and denoted by $\nu$. 
Next our goal is to calculate the mean function and covariance operator with the help of SGD iteration form.

Combining Assumptions \ref{a:quad} and \ref{a:H}, in region $\Xi$, the gradient of $\mathcal{L}(u^M)$ (first \( M \) modes) is given by
\begin{align*}
	\mathcal{G}(u^M) = \mathcal{A}(u^M - \bar{u}^M) +\mathcal{C}^{-1}_0 (u^M-\bar{u}^M).
\end{align*}
Then, the SGD iteration Eq.~\eqref{eq:SGDiter} restricted to the first \( M \) modes becomes
\begin{equation}\label{eq:OUprocess}
	u_{k+1}^M = (I-\eta \mathcal{A}-\eta \mathcal{C}^{-1}_0)u_k^M + \eta \mathcal{A}\bar{u}^M+\eta\mathcal{C}^{-1}_0\bar{u}^M + \frac{\eta}{\sqrt{S}}\xi_{M, k}.
\end{equation}
\( \xi_{M, k} \) is the first \( M \) modes of gradient noise \( \xi_k =\Delta \mathcal{G}(u_k) \) following $\mathcal{N}(0, \mathcal{C}_\text{GN})$.
This noise satisfies
\begin{align*}
	\xi_k = \sum^M_{i = 1} \sqrt{c_iq_i}\zeta_ie_i +  \sum^\infty_{i = M+1} \sqrt{c_i}\zeta_ie_i, \text{with} \quad \zeta_i \sim \mathcal{N}(0, 1).
\end{align*}
such that \( \xi_k^M = \sum^M_{i = 1} \sqrt{c_iq_i}\zeta_ie_i \).

By leveraging the eigensystem \( \{c_i, e_i\}^{\infty}_{i=1} \) of \( \mathcal{C}_0 \), we are able to calculate the coordinate values \( \{u_i\}^\infty_{i=1} \) separately of corresponding  eigen-functions \( \{e_i\}^\infty_{i=1} \).
Under Assumption \ref{a:H}, Hessian operator $\mathcal{A}$ is represented by
\begin{align*}
	\mathcal{A} = \sum^M_{i = 1}a_ie_i,
\end{align*}
where each \( a_i>0 \) , and $ \mathcal{A} $ is self-adjoint, positive definite and trace-class.
Then, the component-wise update rule for $u_k$ becomes 
\begin{align}\label{eq:SGDcoord}
	u_{k+1,i} = 
	\begin{cases}
		\Bigl(1 - \eta (a_i+\frac{1}{c_i})\Bigr) u_{k,i} + \eta (a_i+\frac{1}{c_i}) \bar{u}_i + \frac{\eta}{S} \sqrt{c_iq_i}\zeta_i, &  i = 1, \cdots, M, \\
		\sqrt{c_i}\zeta_i, & i = M+1, M+2, \cdots.
	\end{cases}
\end{align}

For notation convience, we denote \( \tilde{a}_i = a_i+\frac{1}{c_i} \), for \( i = 1, \cdots, M \).
The iteration step is divided into two parts according to the coordinate index,
as the information is concentrated on the first \( M \) eigenvalues.
As for index \( i \leq M \), the update is actually SGD iteration form Eq.~\eqref{eq:OUprocess}.
When index \( i > M \), due to the lack of data information, the posterior measure of parameter \( u \) is considered as the prior measure.

Following the update of each coordinate Eq.~\eqref{eq:SGDcoord}, the iteration form of $u$ is then divided into two modes: the activate modes for $i \leq M$, and inactivate modes for $i > M$.
Then we discuss these two parts individually.

\subsubsection*{Activate modes(\( i \leq M \))}
For these modes, the update involves both data and prior terms. 
The gradient of $\mathcal{L}(u)$ incorporates the update of $u$, then we have
\begin{align*}
	u_{k+1,i} = \Bigl(1 - \eta \tilde{a}_i\Bigr) u_{k,i} + \eta \tilde{a}_i \bar{u}_i  + \frac{\eta}{S} \sqrt{c_iq_i}\zeta_i.
\end{align*}
Taking expectation on both sides, we have
\begin{align*}
	\mathbb{E}[u_{k+1,i}] &= \Bigl(1 - \eta \tilde{a}_i\Bigr) \mathbb{E}[u_{k,i}] + \eta \tilde{a}_i \bar{u}_i  + \frac{\eta}{S} \sqrt{c_iq_i}\mathbb{E}[\zeta_i] \\
	&= \Bigl(1 - \eta \tilde{a}_i\Bigr) \mathbb{E}[u_{k,i}] + \eta \tilde{a}_i \bar{u}_i.
\end{align*}
where $\zeta_i \sim \mathcal{N}(0, 1)$.
Setting $\mathbb{E}[u_i] := \mathbb{E}[u_{k+1,i}] = \mathbb{E}[u_{k,i}]$, we have
\begin{align*}
	\mathbb{E}[u_i] &= \Bigl(1 - \eta \tilde{a}_i\Bigr) \mathbb{E}[u_i] + \eta \tilde{a}_i \bar{u}_i,
\end{align*}
deriving
\begin{align*}
	\mathbb{E}[u_i] &= \bar{u}_i.
\end{align*}

By denoting $s_i = \mathbb{E}[u_{k, i} - \bar{u}_i]^2$ as the stationary variance in the $i-$th mode, the corresponding discrete-time Lyapunov equation is
\begin{align*}
	s_i = \bigg(1 - \eta \tilde{a}_i \bigg)^2s_i + \bigg(\frac{\eta}{S} \bigg)^2 c_iq_i.
\end{align*}
For the stable and well-defined solution, the coefficient satisfies
\[
\bigg(1 - \eta \tilde{a}_i \bigg)^2 < 1,
\]
deriving
\[
\eta < \frac{2}{\tilde{a}_i}.
\]
As a result, in order to guarantee the conditions of Lyapunov equation for all \( i \leq M \), learning rate satisfies
\begin{align}\label{eq:etacondition}
	\eta < \min \bigg\{ \frac{2}{\tilde{a}_1}, \cdots, \frac{2}{\tilde{a}_M} \bigg\} := \frac{2}{\tilde{a}_{\text{max}}},
\end{align}
where \( \tilde{a}_{\text{max}} \) denotes the biggest value of \( \{ a_i + 1/c_i \}^M_{i=1} \).

Subtracting the drift part, we obtain
\begin{align*}
	\bigg( 1 - (1-\eta \tilde{a}_i)^2 \bigg)s_i &= \frac{\eta^2}{S^2}c_iq_i \\
	\bigg( 2\eta \tilde{a}_i - \eta^2\tilde{a}^2_i \bigg)s_i &= \frac{\eta^2}{S^2}c_iq_i.
\end{align*}
Then $i-$th mode's stationary variance is calculated by
\begin{align*}
	s_i &= \frac{\eta c_iq_i}{S^2(2\tilde{a}_i - \eta \tilde{a}^2_i)}.
\end{align*}

\subsubsection*{Inactivate modes($i > M$)}
These modes are uninformed by the data, and the SGD update reduces to direct sampling from the prior:
\begin{align*}
	u_{k+1,i} = \sqrt{c_i}\zeta_i,
\end{align*}
indicating that the iteration proceed turns to sample from the prior measure.
And $\mathbb{E}[u_i] = 0$.
Then $i-$th mode's stationary variance is derived by
\begin{align*}
	s_i = c_1.
\end{align*}

As a result, the stationary distribution { variance} of iteration form Eq.~\eqref{eq:SGDcoord} is then summarized by
\begin{align*}
	s_i  = 
	\begin{cases}
		\frac{\eta c_iq_i}{S^2(2\tilde{a}_i - \eta \tilde{a}^2_i)} &  i = 1, \cdots, M, \\
		c_i, & i = M+1, M+2, \cdots.
	\end{cases}
\end{align*}
Then the estimated posterior measure \( \nu \) is obtained by
\[
\nu = \mathcal{N}(\bar{u}_{\nu}, \mathcal{C}_{\nu}),
\]
where
\begin{align}\label{eq:estpost}
	\bar{u}_{\nu} = \sum^M_{i=1}\bar{u}_i \quad \mathcal{C}_{\nu} = \frac{\eta c_iq_i}{S^2(2\tilde{a}_i - \eta \tilde{a}^2_i)} e_i \otimes e_i + \sum^\infty_{i = M+1} c_i e_i \otimes e_i.
\end{align}

The mean function and covariance operator of estimated posterior measure are parameterized by the learning rate, indicating \( \eta \) determining how ``close'' \( \mu \) and \( \nu \) is. 
Next, in order to approximate the posterior measure, we employ Bayesian inference method and find the optimal learning rate by solving the inference problem.

\subsubsection{SGD as Bayesian inference problem}\label{subsec:2.2.3}
Let $\mathcal{H}_u$ be a separable Hilbert space, and let $\mathcal{M}(\mathcal{H}_u)$ denote the set of Borel probability measures on $\mathcal{H}_u$.
The posterior measure $\mu$ with respect to a prior $\mu_0$ defined on $\mathcal{H}$ is defined via Bayes’ formula ($\ref{eq:post}$). 
For any $\nu \in \mathcal{M}(\mathcal{H}_u)$ that is absolutely continuous with respect to $\mu$, 
the Kullback-Leibler (KL) divergence is defined as 
\begin{align*}
	D_{\text{KL}}(\nu || \mu) &= \int_{\mathcal{H}} \log \bigg(\frac{d\nu}{d\mu}(u) \bigg)\nu(du)
	=\mathbb{E}^{\nu} \bigg[\log \bigg(\frac{d\nu}{d\mu}(u) \bigg) \bigg].
\end{align*}
Here, the notation $\mathbb{E^{\nu}}$ means taking expectation with respect to the probability measure $\nu$.
If $\nu$ is not absolutely continuous with respect to $\mu$, the KL divergence is defined as $+\infty$. 
With this definition, this paper examines the following minimization problem
\begin{align}\label{eq:inference}
	\mathop{\arg\min}_{\eta \in \Omega}D_{\text{KL}}(\nu || \mu),
\end{align}
where  \( \Omega \) is a set of learning rate satisfying the stable and well-defined condition Eq.~\eqref{eq:etacondition}.
The inference method aims to find the closest probability measure $\nu$ to the posterior measure $\mu$ with respect to the KL divergence.

Before employing the Bayesian inference method, we divide our goal into two parts:
\begin{itemize}
	\item Based on Eq.~\eqref{eq:estpost}, both the mean function and the covariance operator of the estimated posterior are parameterized by the learning rate.
	Then the optimal \( \eta \) can be calculated by solving the minimization inference problem between \( \nu \) and \( \mu \), which is measured by Kullback-Leibler divergence.
	\item To ensure the stationary distribution of the discrete-time process Eq.~\eqref{eq:OUprocess}, learning rate \( \eta \in \Omega \) is required.
	We achieve this goal by carefully choosing scale parameter \( S \) through the explicit expression of \( \eta \).
\end{itemize}
To address the first goal, we present the following theorem for determining the optimal learning rate \( \eta \).

\begin{theorem}\label{the:cSGDeta}
	Inference problem
	\[
	\mathop{\arg\min}_{\eta \in \Omega}D_{\text{KL}}(\nu || \mu)
	\]
	processes a solution \( \eta^\dagger \) with the following form:
	\[
	\eta^\dagger =  \frac{2S^2\sum^M_{i = 1}\frac{\tilde{a}_i}{q_i}}{M^\prime + S^2\sum^M_{i = 1}\frac{\tilde{a}^2_i}{q_i}},
	\] 
	where \( M^\prime = M + \sum_{i=1}^Ma_ic_i \), \( \tilde{a}_i = a_i+\frac{1}{c_i} \), \( \nu, \mu \) are the estimated and real posterior measure of Eq.~\eqref{eq:estpost} and Eq.~\eqref{eq:post}, respectively; \( S \) is the scale parameter of stochastic gradient noise; \( \{ a_i, q_i \}^M_{i=1} \) are the eigenvalues of operators \( \mathcal{A} \) and \( \mathcal{Q} \).
\end{theorem}

{ 
\begin{remark}
	Under the assumptions of Theorem~\ref{the:cSGDeta} and taking the notation, we can obtain the optimal $S^\dagger$, which is calculated by
	\[
	S^\dagger
	=\sqrt{
	\frac{\eta M^\prime}
	{\displaystyle\sum_{i=1}^M \frac{\,2\tilde{a}_i-\eta \tilde{a}_i^{\,2}\,}{q_i}}
	}.
	\]
\end{remark}
}

To address the second goal, we require that the optimal learning rate \( \eta^\dagger \in \Omega \), i.e., 
\[
\eta^\dagger =  \frac{2S^2\sum^M_{i = 1}\frac{\tilde{a}_i}{q_i}}{M^\prime + S^2\sum^M_{i = 1}\frac{\tilde{a}^2_i}{q_i}}
< \frac{2}{\tilde{a}_\text{max}},
\]
deriving
\[
\frac{\tilde{a}_{\text{max}}S^2\sum^M_{i = 1}\frac{\tilde{a}_i}{q_i}}{M^\prime + S^2\sum^M_{i = 1}\frac{\tilde{a}^2_i}{q_i}} < 1.
\]
Then we have
\begin{align}\label{eq:etaineq}
	S^2\sum^M_{i = 1}\frac{\tilde{a}_{\text{max}}\tilde{a}_i}{q_i} < M^\prime + S^2\sum^M_{i = 1}\frac{\tilde{a}^2_i}{q_i},
\end{align}
then
\begin{align}\label{eq:ScSGD}
	S < \sqrt{\frac{M^\prime}{\sum^M_{i = 1}\frac{\tilde{a}_{\text{max}}\tilde{a}_i - \tilde{a}^2_i}{q_i}}}.
\end{align}

By appropriately bounding $S$, we guarantee that the Lyapunov equation for each mode yields a stable and well-defined stationary distribution.
Combining the discussion in Subsections \ref{subsec:2.2.1}, \ref{subsec:2.2.2} and \ref{subsec:2.2.3}, we provide the cSGD-iVI method, see Algorithm \ref{alg A}. 

\subsubsection{Discretization error of SGD}\label{subsec:2.2.4}
Based on the discussion in Subsection \ref{subsec:2.2.1}, the update form of the SGD iteration is divided into two parts: the active modes and inactive modes according to \( M \).
Thus the infinite-dimensional inverse problem Eq.~\eqref{eq:invproblem} can be regarded as two components: finite-dimensional problem with data information in \( \mathcal{H}^M_u \) and infinite-dimensional problem with prior information in \( \mathcal{H}^{\perp}_u \).
For solving finite-dimensional problem, the error bound derived by the discretization arisen.
In this section, by truncating the eigenvalues \( \{ c_i \}_{i=1}^\infty \) of \(\mathcal{C}_0 \), the discretization error is also introduced, which should be considered.
Recalling the main reason of optimizing learning rate $\eta^{\dagger}$ is to find the approximated posterior measure, we intend to illustrate the discretization error between the approximated posterior mean function and background truth.

With the optimal learning rate $\eta^{\dagger}$, SGD iteration form proposed in { $\eqref{eq:SGD}$} is then written by
\begin{align}\label{eq:optimalSGD}
	u_{k+1} = u_{k} - \eta^{\dagger}( \mathcal{G}(u_{k}) - \frac{1}{S}\Delta \mathcal{G}(u_{k})).
\end{align}
As we mentioned in Subsection $\ref{subsec:2.2.3}$, sampling in estimated posterior measure is proceeded by running SGD iteration with constant learning rate.
Considering the first $M$ modes of \( u\), which can be regarded as employing projection operator \( P^M \), equation $(\ref{eq:optimalSGD})$is discretized as
\begin{align}\label{eq:optimalSGDfd}
	u_{k+1}^M = u_k^M - \eta^{\dagger}( {\mathcal{G}}(u_k^M) - \frac{1}{S}\Delta {\mathcal{G}}(u_k^M)).
\end{align}
In the following, we discuss the discretization error between approximated posterior mean function and background truth in finite-dimensional spaces.

We express the loss function $\mathcal{L}(u_k^M)$ at \(k-\)th step by
\begin{align*}
	\mathcal{L}(u_k^M) \approx \frac{1}{2}\bigg(\| {\mathcal{A}}^{1/2}(u_k^M - \bar{{u}}_M) \|^2 +\| {\mathcal{C}}^{-1/2}_0(u_k^M - \bar{u}^M) \|^2 \bigg)
    :=\frac{1}{2}\| {A}(u_k^M - \bar{u}^M) \|^2,
\end{align*}
then and the full gradient ${\mathcal{G}}(u_k^M)$ is
\begin{align*}
	{\mathcal{G}}(u_k^M) = {A}^{1/2}({A}^{1/2}(u_k^M - \bar{u}^M)).
\end{align*}
For notational convenience, we denote $\widetilde{{d}} := {A}^{1/2}\bar{u}^M \in \mathcal{H}^M_u$.
Taking the expectation with respect to transition probability measure ${\nu}_k$ of $(\ref{eq:optimalSGDfd})$, approximated posterior mean function is generated by
\begin{align}\label{eq:lwbiterfd}
	\begin{split}
		\mathbb{E}^{{\nu}_k} \bigg[u_{k+1}^M \bigg] &= \mathbb{E}^{\nu_k^M} \bigg[  u_k^M - \eta^{\dagger}( {\mathcal{G}}(u_k^M) - \frac{1}{S}\Delta {\mathcal{G}}(u_k^M)) \bigg] \\
		&= \mathbb{E}^{{\nu}_k^M} [u_k^M ] - \eta^{\dagger}\mathbb{E}^{{\nu}_k^M} [{\mathcal{G}}(u_k^M) ] \\
		&= \mathbb{E}^{{\nu}_k^M} [u_k^M ] - \eta^{\dagger} {A}^{1/2}({A}^{1/2}\mathbb{E}^{{\nu}_k^M}[{u^M}] - \widetilde{{d}}).
	\end{split}
\end{align}
Recalling $\Delta {\mathcal{G}}(u_k^M)$ is random Gaussian variable with mean zero, we obtain the conclusion.

\begin{assumption}\label{assump:solution}
	Let $u_0$ be the initial value of the iteration form $(\ref{eq:SGDiter})$.
	There is an index function $\varphi$, and $v \in \mathcal{H}_u$, with $\|v \|_{\mathcal{H}^N_u} \leqslant 1$, such that
	\begin{align*}
		u^{\dagger} - u_0  = \varphi(A^{-1})v,
	\end{align*}
	where $u^{\dagger}$ is the background truth function, \( A := \mathcal{A} + \mathcal{C}^{-1}_0 \).
	And function $\varphi^2$ is assumed to a concave operator.
\end{assumption}

\begin{lemma}\label{lemma:pre}
	Let $\rho$ be a index function satisfying $\exp (-t/\beta) \rho(t) \leq C_{\rho}\rho(\beta)$, where $C_{\rho}$ is a constant. 
	Setting $\alpha_k := (\sum^k_{i=1}\eta_i)^{-1}$, define:
	\begin{align*}
		g_{\alpha_k}(t) := \sum^k_{j=1}\eta_j\prod^k_{i=j+1}(1-\eta_it), \quad 0< t < 1/\eta_1.
	\end{align*}
	Then there hold, for $r_{\alpha_k} = 1 - tg_{\alpha_k}(t)$, that
	\begin{itemize}
		\item $	r_{\alpha_k} = \prod^k_{i=1}(1 - \eta_it)$;
		\item $| g_{\alpha_k}(t) | \leq 1/\alpha_k$;
		\item $| r_{\alpha_k} | \leq 1$;
		\item For each index function $\rho$, $| r_{\alpha_k} |\rho(t) \leq C_{\rho}\rho(\alpha_k)$.
	\end{itemize}
\end{lemma}

This lemma provides the regularization property \cite{Peter2003IP1} of the cSGD formulation.
Then the following theorem illustrates the error bound estimation between approximated posterior mean and background truth functions.

\begin{theorem}\label{the:errbound}
    Let $D_{ M} := \| H(I - P^M) \|_{\mathcal{H}_u \rightarrow \mathcal{H}^M_u}$ denote the approximation degree of the projection $P^M$.
    Under Assumption $\ref{assump:solution}$, we assign to the index function $\varphi$ the companion function $\theta(t) := \sqrt{t}\varphi(t)$.
	Then there exist constant $0 < C < \infty$, such that for learning rate $\eta^{\dagger}$ from $(\ref{eq:lwbiterfd})$,
	\begin{align*}
		\bigg \| \mathbb{E}^{\nu} [u_k^M ] - {u}^{M, \dagger} \bigg \|^2 \leq
		C \bigg(\varphi(\alpha_k) + \varphi(D^2_M) + \frac{1}{\sqrt{\eta^{\dagger}}}\theta(D^2_M) \bigg),
	\end{align*}
	where $\alpha_k = (k\eta^{\dagger})^{-1}$, and $ u^{M, \dagger} = P^Mu^{\dagger}$ is the discretization of background truth $u^{\dagger}$.
	
\end{theorem}

Theorem $\ref{the:errbound}$ illustrates that the discretization error bound between the estimated posterior mean function and background truth is controlled by optimal learning rate $\eta^{\dagger}$ and truncated level $M$.

\begin{algorithm}
	\caption{The cSGD-iVI algorithm}
	\label{alg A}
	\begin{algorithmic}[1]
		\STATE{Initialize $u_0 = 0$, specify number of samples $K$, averaging converging steps $J$, scale $S_0$, threshold \( tol \), and discretization level $N$;}
		\STATE{Determine \( M = \min \{ M \in \mathbb{N} |~ \frac{c_M}{c_1} < C_M \} \);}
		\STATE{Calculate the eigenvalues \( \{ a_i, c_i, q_{i} \}^M_{i = 1} \) of operators \( \mathcal{A}, {\mathcal{C}_0}, \mathcal{Q} \);}
		\REPEAT
		\STATE{ Calculate learning rate
				{ \( \eta_k^\dagger = \frac{2S_{k-1}^2\sum^M_{i = 1}\frac{\tilde{a}_i}{q_i}}{M^\prime + S_{k-1}^2\sum^M_{i = 1}\frac{\tilde{a}^2_i}{q_i}} \) and \(S_k
				=\sqrt{
				\frac{\eta_k^\dagger M^\prime}
				{\sum_{i=1}^M \frac{\,2\tilde{a}_i-\eta_k^\dagger \tilde{a}_i^{\,2}\,}{q_i}}
				}\),} \\
				where \( M^\prime = M + \sum_{i=1}^Ma_ic_i \) and \( \tilde{a} = a_i + 1/c_i \) ;}
		\STATE{Run cSGD iteration form by \\
			$u_{(k+1),j} =u_{k,j} - \eta^{\dagger}_k\widetilde{\mathcal{G}}_S(u_k)$;}
		\UNTIL{\( \| u_{(k+1), j} - u_{k, j} \| / \| u_{k, j} \| \leq tol \);} 
		\STATE{Return samples $\lbrace u_{k,J}, u_{(k+1),J}, \cdots, u_{(k+K),J} \rbrace$.}
	\end{algorithmic}
\end{algorithm}

\subsection{Preconditioned SGD as variational inference methods}\label{subsec2.3}
In this section, we aim to establish a sampling method based on preconditioned constant SGD (pcSGD), which is a variant of the constant SGD method. 
Since preconditioning improves convergence and enhances computational accuracy, the choice of the preconditioning operator is crucial. 
Our goal is to determine the optimal learning rate and the corresponding preconditioning operator by minimizing the variational inference objective, following the discussion in Subsection $\ref{subsec2.2}$.

Let us consider a bounded linear symmetric preconditioning operator \( T :\mathcal{H}_u \rightarrow \mathcal{H}_u\) satisfying 
\[
T(u) = T(\sum^\infty_{i=1}u_ie_i) \approx T(\sum^M_{i=1}u_ie_i + \sum^\infty_{i=M+1}0e_i).
\]

Then the iteration takes the form:
\begin{align}\label{eq:SGDiterprecon}
	\begin{split}
		u_{(k+1)} &= u_{k} - \eta T\widetilde{ \mathcal{G}}_S(u_{k\eta}) \\
		&= u_{k} - \eta T\bigg( \mathcal{G}(u_{k}) - \frac{1}{S}\Delta \mathcal{G}(u_{k}) \bigg)\\
		&= u_{k} - \eta T\mathcal{G}(u_{k}) + \frac{\eta}{S}T\Delta \mathcal{G}(u_{k}),
	\end{split}
\end{align}
where and \( \Delta \mathcal{G}(u_{k}) \sim \mathcal{N}(0, \mathcal{C}_{\text{noise}})\) denotes the stochastic gradient noise.
By employing the framework proposed in Subsection \ref{subsec2.2}, pcSGD-iVI is then generated.

Recalling notation \( \tilde{a}_i = a_i + 1/c_i \), the update of each coordinate becomes
\begin{align}\label{eq:SGDcoordprecon}
	u_{k+1,i} = 
	\begin{cases}
		\Bigl(1 - \eta t_i\tilde{a}_i\Bigr) u_{k,i} + \eta t_i\tilde{a}_i \bar{u}_i + \frac{\eta}{S} t_i\sqrt{c_iq_i}\zeta_i, &  i = 1, \cdots, M, \\
		\sqrt{c_i}\zeta_i, & i = M+1, M+2, \cdots,
	\end{cases}
\end{align}
where \( \zeta_i \sim \mathcal{N}(0, 1) \), and \( \{t_i\}_{i=1}^\infty \) are the eigenvalues of preconditioner \( T \).
As in the cSGD-iVI case, we analyze the two parts separately:

\subsubsection*{Activate modes(\( i \leq M \))}
The update is
\[
u_{k+1,i} = \Bigl(1 - \eta t_i\tilde{a}_i\Bigr) u_{k,i} + \eta t_i\tilde{a}_i \bar{u}_i + \frac{\eta}{S} t_i\sqrt{c_iq_i}\zeta_i.
\]
Taking the expectation on both side, we obtain
\[
\mathbb{E}[u_i] = \bar{u}_i.
\]
Let $s_i = \mathbb{E}[u_{k, i} - \bar{u}_i]^2$ as the stationary variance in the $i-$th mode, the corresponding discrete-time Lyapunov equation is
\begin{align*}
	s_i = (1 - \eta t_i\tilde{a}_i)^2s_i + \left(\frac{\eta}{S} \right)^2 t^2_ic_iq_i,
\end{align*}
deriving
\[
s_i = \frac{\eta t_ic_iq_i}{S^2(2\tilde{a}_i - \eta t_i\tilde{a}^2_i)}.
\]

\subsubsection*{Inactivate modes(\( i > M \))}
These modes are unaffected by the data and are sampled directly from the prior:
\[
s_i = c_i.
\]

As a result, the stationary distribution variance of iteration form Eq.~\eqref{eq:SGDcoordprecon} is then summarized by
\begin{align*}
	s_i  = 
	\begin{cases}
		\frac{\eta t_ic_iq_i}{S^2(2\tilde{a}_i - \eta t_i\tilde{a}^2_i)}, &  i = 1, \cdots, M, \\
		c_i, & i = M+1, M+2, \cdots.
	\end{cases}
\end{align*}
Then the estimated posterior measure \( \nu \) is obtained by
\[
\nu = \mathcal{N}(\bar{u}_{\nu}, \mathcal{C}_{\nu}),
\]
where
\begin{align}\label{eq:estpostprecon}
	\bar{u}_{\nu} = \sum^M_{i=1}\bar{u}_i \quad C_{\nu} = \sum^M_{i = 1}  \frac{\eta t_ic_iq_i}{S^2(2\tilde{a}_i - \eta t_i\tilde{a}^2_i)}e_i \otimes e_i + \sum^\infty_{i = M+1} c_i e_i \otimes e_i.
\end{align}

For the stable and well-defined solution, the coefficient satisfies
\[
(1 - \eta t_i\tilde{a}_i)^2 < 1,
\]
deriving
\[
\eta < \frac{2}{t_i\tilde{a}_i}.
\]
As a result, in order to guarantee the conditions of Lyapunov equation for all \( i \leq M \), learning rate satisfies
\begin{align}
	\eta < \min \bigg\{ \frac{2}{t_1\tilde{a}_1}, \cdots, \frac{2}{t_M\tilde{a}_M} \bigg\} := \frac{2}{\tau_\text{max}},
\end{align}
where \( \{\tau_i = t_i\tilde{a}_i\}^M_{i=1} \), and \( \tau_\text{max} = \max \{\tau_i\}^M_{i=1}\).

For the pcSGD-iVI case, we provide the similar theorem as in the cSGD-iVI case as follows.

\begin{theorem}\label{the:pcSGDeta}
	Inference problem
	\[
	\mathop{\arg \min}_{\eta \in \Omega_T} D_{\text{KL}}(\nu || \mu)
	\]	
	processes a solution \( \eta^\dagger \) with the following form:
	\[
	\eta^\dagger = \frac{2S^2\sum^M_{i = 1}\frac{\tilde{a}_i}{t_iq_i}}{M^\prime + S^2\sum^M_{i = 1}\frac{\tilde{a}^2_i}{q_i}}.
	\] 
	where \( M^\prime = M + \sum_{i=1}^Ma_ic_i \), and \( \tilde{a}_i = a_i+\frac{1}{c_i} \), \( \nu, \mu \) are the estimated and real posterior measure of Eq.~\eqref{eq:estpost} and Eq.~\eqref{eq:post}, respectively; \( S \) is the scale parameter of stochastic gradient noise; \( \{ t_i, a_i, c_i, q_i \}^M_{i=1} \) are the eigenvalues of operators \( T, \mathcal{A}, \mathcal{C}_0 \) and \( \mathcal{Q} \).
\end{theorem}


{ 
\begin{remark}
	Under the assumptions of Theorem~\ref{the:pcSGDeta} and taking the notation, we can obtain the optimal $S^\dagger$, which is calculated by
	\[
	S^\dagger 
	= \sqrt{\frac{\eta M^\prime }{\displaystyle\sum^M_{i = 1}\frac{(2\tilde{a}_i - \eta t_i\tilde{a}^2_i)}{t_iq_i}}}.
	\]
\end{remark}
}

Following the discussion of \( S \), learning rate satisfies \( \eta \in \Omega_T \).
That is,
\[
\eta^\dagger = \frac{2S^2\sum^M_{i = 1}\frac{t_i\tilde{a}_i}{q_i}}{M^\prime + S^2\sum^M_{i = 1}\frac{\tilde{a}^2_i}{q_i}} < \frac{2}{\tau_{\text{max}}}.
\]
deriving
\[
\frac{\tau_{\text{max}}S^2\sum^M_{i = 1}\frac{\tilde{a}_i}{t_iq_i}}{M^\prime + S^2\sum^M_{i = 1}\frac{\tilde{a}^2_i}{q_i}} < 1.
\]
Then 
\begin{align}\label{eq:SpcSGD}
	S < \sqrt{\frac{M^\prime}{\sum^M_{i = 1}\bigg(\frac{\tau_{\text{max}}\tilde{a}_i}{t_iq_i} - \frac{\tilde{a}^2_i}{q_i}\bigg)}}.
\end{align}

By bounding \( S\) appropriately, we ensure that the Lyapunov equation with the preconditioner \( T \) yields a stable and well-defined stationary distribution. 
Combining the discussion in Subsections \ref{subsec2.2} and \ref{subsec2.3}, we provide the pcSGD-iVI method, see Algorithm \ref{alg B}.

\begin{algorithm}
	\caption{The pcSGD-iVI algorithm}
	\label{alg B}
	\begin{algorithmic}[1]
		\STATE{Initialize $u_0 = 0$, specify number of samples $K$, averaging converging steps $J$, scale $S$, threshold \( \epsilon \), preconditioning operator \( T \) and tolerance \( tol \);}
		\STATE{Determine \( M = \min \{ M \in \mathbb{N} |~ \frac{c_M}{c_1} < C_M \} \);}
		\STATE{Calculate the eigenvalues \( \{ a_i, c_i, q_{i}, {t_i} \}^M_{i = 1} \) of operators \( \mathcal{A}, \mathcal{C}_0, \mathcal{Q}, {T} \);}
		\REPEAT
		\STATE{ Calculate learning rate
				{ \( \eta_k^\dagger = \frac{2S_{k-1}^2\sum^M_{i = 1}\frac{\tilde{a}_i}{t_iq_i}}{M^\prime + S_{k-1}^2\sum^M_{i = 1}\frac{\tilde{a}^2_i}{q_i}} \) and \(S_k 
				= \sqrt{\frac{\eta_k^\dagger M^\prime}{\sum^M_{i = 1}\frac{(2\tilde{a}_i - \eta_k^\dagger t_i\tilde{a}^2_i)}{t_iq_i}}}\),} \\
				where \( M^\prime = M + \sum_{i=1}^Ma_ic_i \) and \( \tilde{a}_i = a_i+\frac{1}{c_i} \);}
		\STATE{Run pcSGD iteration form by \\
			$u_{(k+1),j} =u_{k,j} - \eta^{\dagger}_kT\widetilde{\mathcal{G}}_S(u)$;}
		\UNTIL{\( \| u_{(k+1), j} - u_{k, j} \| / \| u_{k, j} \| \leq tol \);} 
		\STATE{Return samples $\lbrace u_{k,J}, u_{(k+1),J}, \cdots, u_{(k+K),J} \rbrace$.}
	\end{algorithmic}
\end{algorithm}

\subsection{Discussion about choosing \( \mathcal{Q} \)}\label{subsec2.4}
Based on Theorems $\ref{the:cSGDeta}$ and $\ref{the:pcSGDeta}$, 
optimal learning rate $\eta^{\dagger}$ relies on the eigenvalues \( \{ q_i\}^M_{i=1} \) especially.
In other words, the extent to which the estimated posterior measure approximates the true posterior measure $\mu$ is influenced by the choice of $\mathcal{Q}$. 
To determine operator $\mathcal{Q}$, we present the following discussion under a more specific setting.

We consider the loss function
\begin{align*}
    \mathcal{L}(u) = \frac{1}{2}\|\mathcal{A}^{1/2}(u - \bar{u})\|^2 + \frac{1}{2}\| \mathcal{C}^{-1/2}_0(u - \bar{u})\|^2 := \frac{1}{2}\| A(u - \bar{u}) \|^2.
\end{align*}
where \(A : \mathcal{H}_u \rightarrow \mathcal{H}_u \) is Hessian operator.
The gradient is
\[
\nabla \mathcal{L}(u) = A(u - \bar{u}).
\]
In what follows, we work in the activation space \(\mathcal{H}_u^M\).
Recall the isomorphism \(\gamma:\mathcal{H}_u \to \ell^2\) introduced in
Subsection~\ref{subsec2.1}. 
By restricting \(\gamma\) to
\(\mathcal{H}_u^M\) and retaining its first \(M\) coordinates, we obtain an isometric isomorphism between \(\mathcal{H}_u^M\) and \(\mathbb{R}^M\).
Without risk of confusion, we use the same notation \(u\) both for function \(u \in \mathcal{H}_u^M\) and for its coordinate vector
\(u = (u_1,\ldots,u_M)^T \in \mathbb{R}^M\).
Next, we consider a random projection operator \(P \in \mathbb{R}^{p\times M}\),
which is a random matrix with independent entries satisfying
\(P_{ij} \sim \mathcal{N}(0,1/p)\). 
Applying this random project to the loss function, we derive
\[
\tilde{\mathcal{L}}(u) = \frac{1}{2}\|PA^{1/2}(u - \bar{u})\|^2,
\]
which is random due to $P$.
With the same procedure, we can compute the gradient of $\tilde{\mathcal{L}}(u)$: 
\begin{align*}
	\nabla \tilde{\mathcal{L}}(u) &= (PA^{1/2})^*PA^{1/2}(u - \bar{u}) \\
	&= A^{1/2*}P^*PA^{1/2}(u - \bar{u}).
\end{align*}
Our goal is to compute the expectation and covariance of $\nabla \tilde{\mathcal{L}}({u})$, which we present in the following theorem.

\begin{theorem}\label{the:stogradient}
	Let forward operator \( H \) be a bounded linear operator, and \( P \) be the random projection matrix.
	Then the stochastic gradient \( \nabla \tilde{\mathcal{L}}(u) \) is a Gaussian random variable, with mean function and covariance operator being defined by
	\begin{align}\label{eq:stogradientdis}
		\mathbb{E}[\nabla \tilde{\mathcal{L}}(u)] &= A(u - \bar{u}); \\
		\operatorname{Cov}(\nabla \tilde{\mathcal{L}}(u)) &= \frac{1}{p} A^{1/2*}\Lambda A^{1/2},
	\end{align}
	where \( \Lambda = \|A^{1/2}( u - \bar{u})\|^2 I_{M}+(A^{1/2}( u - \bar{u}))(A^{1/2}( u - \bar{u}))^* \).
\end{theorem}

\begin{figure}
	\centering
	{
		\includegraphics[ keepaspectratio=true, width=0.6\textwidth,  clip=true, trim=80pt 2pt 80pt 30pt]{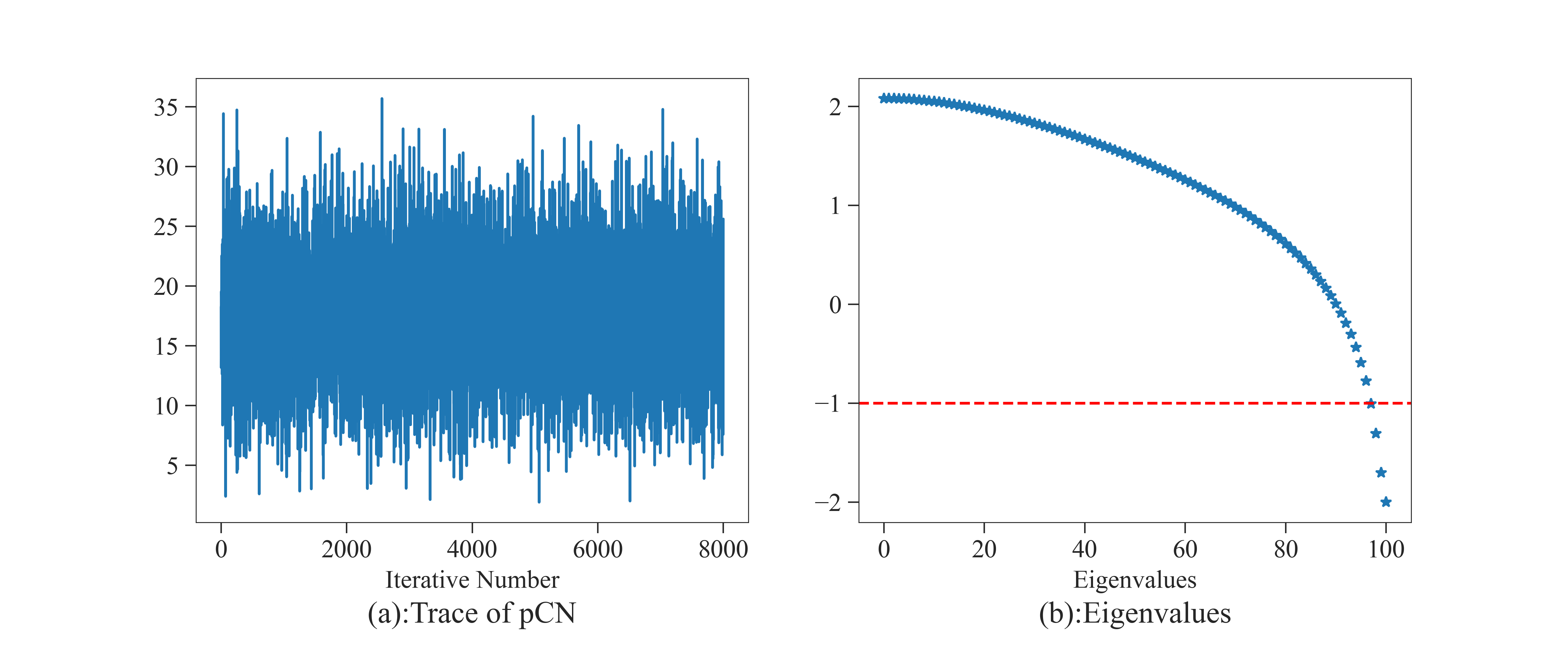}}

	\caption{\emph{\small 
			(a): The trace plot of pCN method under the same mesh size of $100$; 
			(b): Logarithm of the eigenvalues $\lbrace c_k \rbrace^{n}_{k=1}$ of prior measure. 
			The horizontal red dashed line shows the corresponding eigenvalue \( c_M \) satisfying \( c_M / c_1 < 10^{-3} \).}}		
	
	\label{fig:trace}
	
\end{figure}

From Theorem \ref{the:stogradient}, the stochastic gradient used in the SGD update is
\[
\nabla \tilde{\mathcal{L}}(u) = \tilde{\mathcal{G}}_S(u) = A(u - \bar{u}) - \frac{1}{S}\xi,
\]
where $ \xi$ is a mean-zero Gaussian error term in $ \mathcal{H}_u $ with covariance operator $\mathcal{C}_\text{noise} = \mathcal{C}_0\mathcal{Q}$. 
Since random variable \( \xi \) is centered, the expected value of \( \tilde{\mathcal{G}}_S(u)\)  is \(A(u - \bar{u})\), which equals to \( \mathbb{E}[\nabla \tilde{\mathcal{L}}(u)] = A(u - \bar{u}) \).

We now characterize the covariance operator of \(\operatorname{Cov}(\nabla \tilde{\mathcal{L}}(u))\), which is:
\[
\operatorname{Cov}(\nabla \tilde{\mathcal{L}}(u)) = \frac{1}{S^2}\mathcal{C}_0\mathcal{Q}.
\]
Then eigenvalues \( \{ \phi_i \}^M_{i = 1} \) of \(\operatorname{Cov}(\nabla \tilde{\mathcal{L}}(u))\) satisfy
\[
\operatorname{Cov}(\nabla \tilde{\mathcal{L}}(u))e_i = \phi_ie_i = \frac{1}{S^2}c_iq_ie_i,
\]
indicating
\[
\phi_i = \frac{1}{S^2}c_iq_i.
\]

That is, by calculate the eigenvalues of covariance \( \operatorname{Cov}(\nabla \tilde{\mathcal{L}}(u)) \), the eigenvalues \( \{q_i\}^M_{i = 1} \) of operator \( \mathcal{Q} \) is then calculated by
\[
\mathcal{Q} = \sum^M_{i = 1}q_i e_i \otimes e_i 
=  \sum^M_{i = 1}\frac{S^2\phi_i}{c_i}e_i \otimes e_i.
\]
Substituting \( \{q_i\}^M_{i = 1} \) into the bounds for \( S \) in Eqs.~\eqref{eq:ScSGD} and \eqref{eq:SpcSGD}, we obtain
\begin{itemize}
	\item for the cSGD-iVI case, we derive
	\[
	S < \sqrt{\frac{M^\prime}{\sum^M_{i = 1}\frac{\tilde{a}_{\text{max}}\tilde{a}_i - \tilde{a}^2_i}{q_i}}} =
	\sqrt{\frac{M^\prime}{\sum^M_{i = 1}\frac{\tilde{a}_{\text{max}}\tilde{a}_i - \tilde{a}^2_i}{\frac{S^2\phi_i}{c_i}}}} = 	\sqrt{S^2\frac{M^\prime}{\sum^M_{i = 1}\frac{\tilde{a}_{\text{max}}c_i\tilde{a}_i - c_i\tilde{a}^2_i}{\phi_i}}}.
	\]
	Then we obtain the necessary condition
	\[
	1 < \frac{M^\prime}{\sum^M_{i = 1}\frac{\tilde{a}_{\text{max}}c_i\tilde{a}_i - c_i\tilde{a}^2_i}{\phi_i}},
	\]
    where \( M^\prime = M + \sum_{i=1}^Ma_ic_i \), and \( \tilde{a}_i = a_i+\frac{1}{c_i} \).
    
	\item for the pcSGD-iVI case, we derive
	\[
	S < \sqrt{\frac{M^\prime}{\sum^M_{i = 1}\bigg(\frac{\tau_{\text{max}}\tilde{a}_i}{t_i\frac{S^2\phi_i}{c_i}} - \frac{\tilde{a}^2_i}{\frac{S^2\phi_i}{c_i}}\bigg)}} 
	=  \sqrt{\frac{M^\prime}{\frac{1}{S^2}\sum^M_{i = 1}\bigg(\frac{\tau_{\text{max}}c_i\tilde{a}_i}{t_i\phi_i} - \frac{c_i\tilde{a}^2_i}{\phi_i}\bigg)}}.
	\]
	Then we obtain the necessary condition
	\[
	1 < \frac{M^\prime}{\sum^M_{i = 1}\bigg(\frac{\tau_{\text{max}}c_i\tilde{a}_i}{t_i\phi_i} - \frac{c_i\tilde{a}^2_i}{\phi_i}\bigg)},
	\]
    where \( M^\prime = M + \sum_{i=1}^Ma_ic_i \), and \( \tilde{a}_i = a_i+\frac{1}{c_i} \).
\end{itemize}


For both cSGD-iVI and pcSGD-iVI method, the upper bound condition of hyperparameter \(S\) is then transformed reformulated as inequalities satisfied by the eigenvalues \(\{\phi_i\}_{i=1}^M\), and further reduced to inequalities on the hyperparameter \(p\).
Consequently, the distribution $\mathcal{N}(0, 1/p)$ of the entries $P_{ij}$ of the projection operator directly determines the eigenvalues $\{q_i\}_{i=1}^M$ of the operator $\mathcal{Q}$, thereby providing a basis for choosing $\mathcal{Q}$.


\begin{figure}
	\centering
	\includegraphics[ keepaspectratio=true, width=0.6\textwidth,  clip=true, trim=80pt 2pt 80pt 30pt]{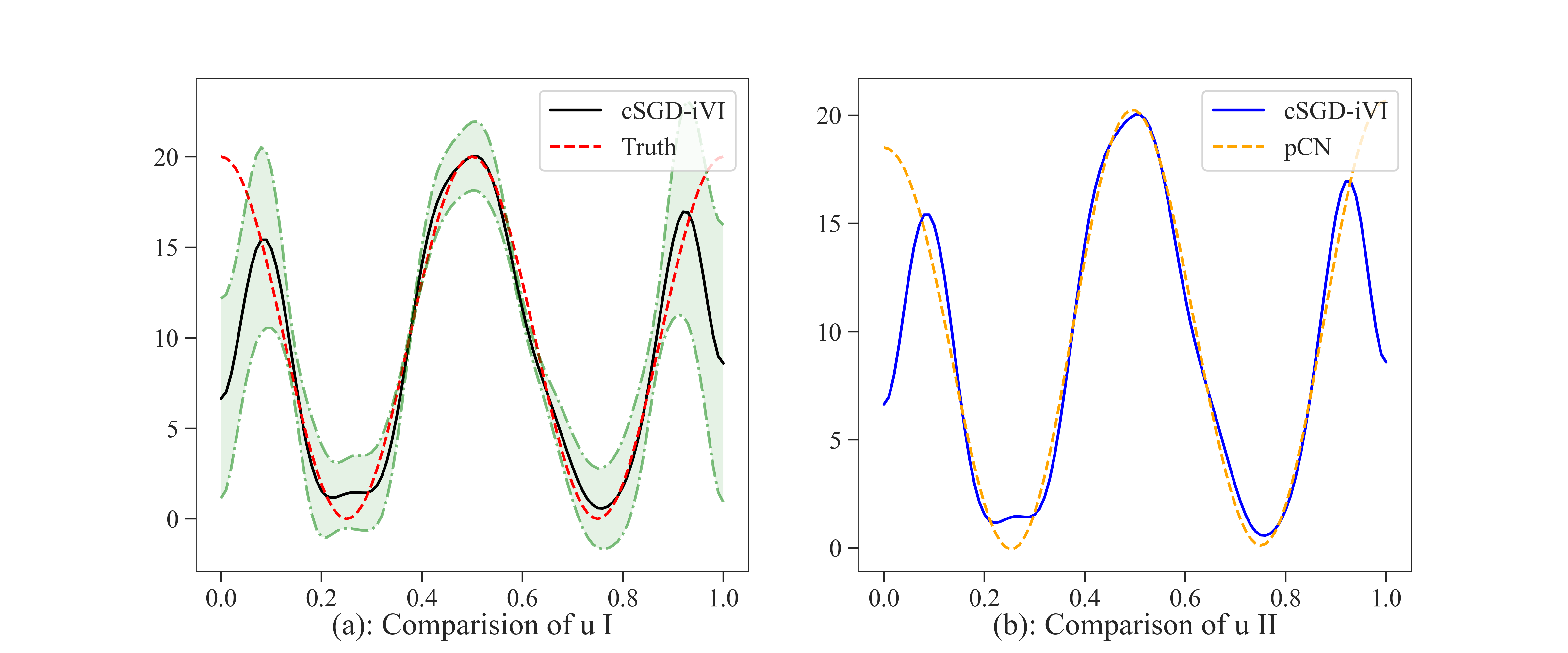}

	\caption{\emph{\small 
			(a): The comparison of the \textcolor{black}{estimated posterior mean function} obtained by cSGD-iVI and the background truth of $u$, respectively.
			The green shade area represents the $95 \%$ credibility region of estimated posterior mean function; 
			(b): The comparison of the \textcolor{black}{estimated posterior mean function} obtained by cSGD-iVI and pCN sampling of $u$, respectively.
	}}
	
	\label{fig:Comparison_cSGD}
	
\end{figure}

\section{Numerical experiments}\label{sec3}
We employ the cSGD-iVI (proposed in Subsection $\ref{subsec2.2}$) and pcSGD-iVI (proposed in Subsection $\ref{subsec2.3}$) methods solve two inverse problems, including one linear and one nonlinear case. 
In Subsection { \ref{subsec3.1}}, we consider a simple smooth model with Gaussian noise; meanwhile, we compare the computational results with those obtained using classical sampling methods, such as the pCN method. 
Due to space limitations, the nonlinear inverse problem governed by the Darcy flow equation is provided in Section B of the supplementary material.

\subsection{The simple elliptic equation with Gaussian noise}\label{subsec3.1}
\subsubsection{Basic settings}\label{subsec3.1.1}
{ The simple elliptic equation we solved in this section is followed from \cite{Sui2024MOC}.}
Consider an inverse source problem of the elliptic equation
\begin{align}\label{prob1}
	\begin{split}
		-\alpha \Delta w + w &= u \quad \text{in}\ \Omega, \\ 
		w &= 0 \quad \text{on}\ \partial \Omega,
	\end{split}
\end{align}
where $\Omega = (0, 1) \subset \mathbb{R}$, $\alpha > 0$ is a positive constant. 
{ Let $A := \text{I} - \alpha\Delta$ and denote the solution operator by $F := A^{-1}$; then $w = F u$.
Define the bounded linear observation operator $O$ by
\[
Ow := \big(w(x_1), \ldots, w(x_{N_d})\big)^{\top}.
\]
The forward operator from the parameter to the data is the composition
\[
H := O \circ F,
\]
satisfying}
\begin{align}
	Hu = (w(x_1), w(x_2), \cdots, w(x_{N_d}))^T,
\end{align}
where $u \in \mathcal{H}_u := L^2(\Omega)$, $w$ denotes the solution of $(\ref{prob1})$, $x_i \in \Omega$ for $i = 1, \cdots, N_d$, {  and \(H : L^2(\Omega) \to \mathbb{R}^{N_d}\) is the bounded linear forward operator}.
With these notations, the problem is written abstractly as:
\begin{align}
	\bm{d} = Hu + \bm{\epsilon},
\end{align}
where $\bm{\epsilon} \sim \mathcal{N}(0, \bm{\Gamma}_{\text{noise}})$ is the random Gaussian noise.
In our implementations, the measurement points $\lbrace x_i \rbrace^{N_d}_{i=1}$ are taken at the coordinates $\lbrace i/20\rbrace^{20}_{i = 1}$. 
To avoid the inverse crime \cite{kaipio2006statistical}, we discretize the elliptic equation by the finite element method on a regular mesh (the grid points are uniformly distributed on the domain $\Omega$) with the number of grid points being equal to $10^4$. 
In our experiments, the prior measure of $u$ is a Gaussian probability measure $\mu_0$ with mean zero and covariance $\mathcal{C}_0$.

\begin{figure}
	\centering
	\includegraphics[ keepaspectratio=true, width=0.6\textwidth,  clip=true, trim=80pt 2pt 80pt 30pt]{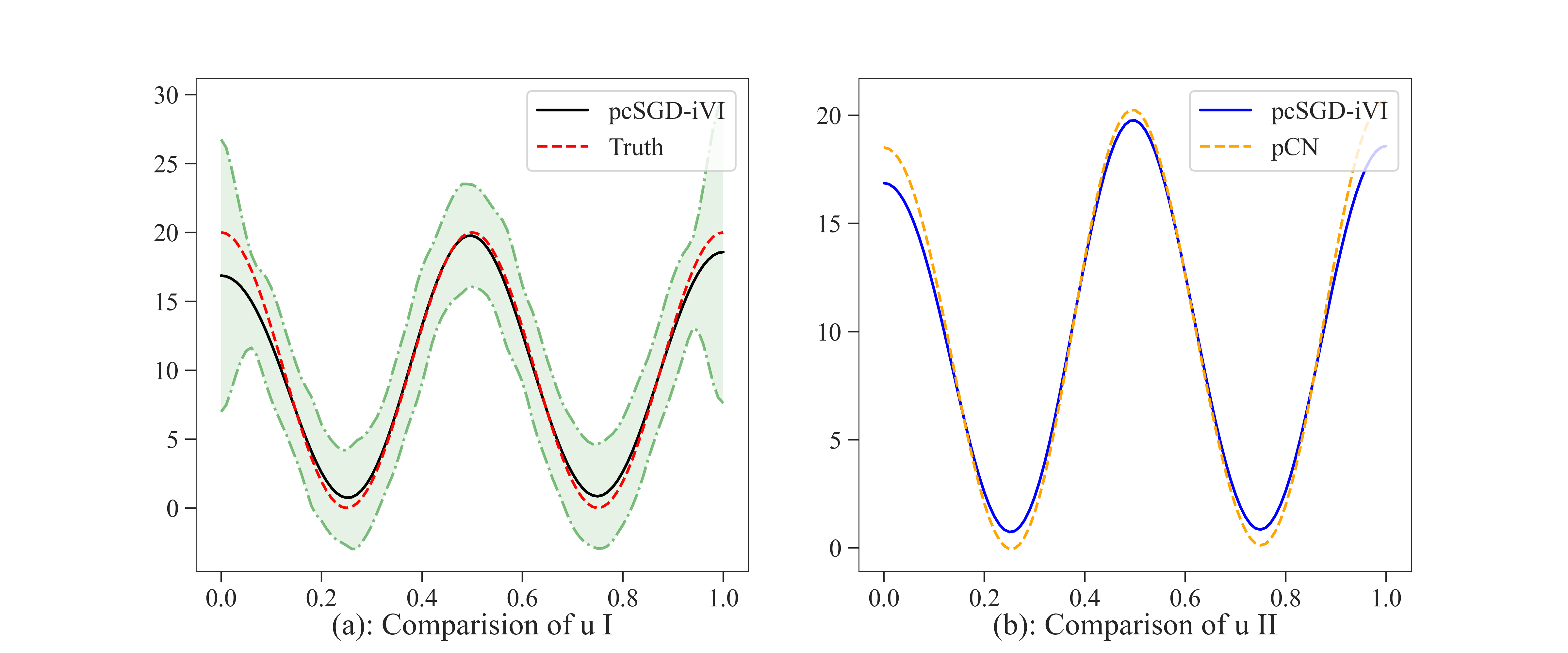}

	\caption{\emph{\small 
			(a): The comparison of the \textcolor{black}{estimated posterior mean function} obtained by pcSGD-iVI and the background truth of $u$, respectively.
			The green shade area represents the $95 \%$ credibility region of estimated posterior mean function; 
			(b): The comparison of the \textcolor{black}{estimated posterior mean function} obtained by pcSGD-iVI and pCN of $u$, respectively.
	}}
	
	\label{fig:Comparison_pcSGD}
	
\end{figure}

For clarity, we list the specific choices for some parameters introduced in this subsection as follows: 
\begin{itemize}
	\item \textcolor{black}{Assume that 5$\%$ random Gaussian noise $\epsilon \sim \mathcal{N}(0, \bm{\Gamma}_{\text{noise}})$ is added, where $\bm{\Gamma}_{\text{noise}} = \tau^{-1}\textbf{I}$, and $\tau^{-1} = (0.05\max(\lvert Hu\rvert))^2$.}
	\item Let domain $\Omega$ be an interval $(0, 1)$ with $\partial \Omega = \lbrace 0, 1 \rbrace$. And the available data are assumed to be $\lbrace w(x_i) | i = 1, 2, \cdots, 20 \rbrace$.
	\item We assume that the data produced from the underlying true signal $u^{\dagger}(x) = 10 \cdot (\cos 4\pi x+1)$.
	\item The prior measure of $u$ is given by $\mu_0 = \mathcal{N}(0, \mathcal{C}_0)$, where $\mathcal{C}_0$  is given by $\mathcal{C}_0 = (\text{I} - \alpha \Delta)^{-2},$ where $\alpha = 0.05$ is a fixed constant. Here the Laplace operator is defined on $\Omega$ with zero Neumann boundary condition.
	\item In order to avoid inverse crime \cite{kaipio2006statistical}, the data is generated on a fine mesh with the number of grid points equal to $10^4$. And we use different sizes of mesh $n = 100$ in the inverse stage.
\end{itemize}

To illustrate the sampling effectiveness of cSGD-iVI and pcSGD-iVI, we compare them to pCN method, which is  wildly studied in \cite{cotter2013, dashti2013bayesian, Pillai2014SPDE}. 

\subsubsection{Numerical results}\label{subsec3.1.2}
In this subsection, we compare the computational results of cSGD-iVI, pcSGD-iVI and pCN methods.
Since the covariance operator of the stochastic gradient noise needs to be carefully handled, we follow Subsection \ref{subsec2.4} to compute the operator \( \mathcal{Q} \).
Based on the formulations of cSGD-iVI and pcSGD-iVI, we adopt different assumptions as described in Subsection C.1 of supplementary material.
Here, we briefly discuss the computational cost of the methods considered in this subsection.

\begin{table}
	\renewcommand{\arraystretch}{1.5}
	\centering
	\caption{\emph{\small The relative errors between the covariance matrix, variance function, and covariance functions derived by cSGD-iVI and pCN methods.}} \label{table:relative_cSGD}
	\begin{tabular}{c|cccc}
		\hline $\text{Function} $ & $\bm{c}$  & $\lbrace c(x_i, x_i)\rbrace^{n}_{i=1}$ & $\lbrace c(x_i, x_{i+10}) \rbrace^{n-10}_{i=1}$ & $\lbrace c(x_i, x_{i+20}) \rbrace^{n-20}_{i=1}$  \\
		\hline $\text{Relative Error}$ & $2.3057$  & $0.2586$ & $1.6433$ & $10.2072$ \\
		\hline
	\end{tabular}
\end{table}

\begin{table}
	\renewcommand{\arraystretch}{1.5}
	\centering
	\caption{\emph{\small The relative errors between the covariance matrix, variance function, and covariance functions derived by pcSGD-iVI and pCN method.}} \label{table:relative_pcSGD}
	\begin{tabular}{c|cccc}
		\hline $\text{Function} $ & $\bm{c}$  & $\lbrace c(x_i, x_i)\rbrace^{n}_{i=1}$ & $\lbrace c(x_i, x_{i+10}) \rbrace^{n-10}_{i=1}$ & $\lbrace c(x_i, x_{i+20}) \rbrace^{n-20}_{i=1}$  \\
		\hline $\text{Relative Error}$ & $0.1703$  & $0.05616$ & $0.1554$ & $0.8971$ \\
		\hline
	\end{tabular}
\end{table}

First of all, the computational cost of the sampling method is high.
In \cite{dunlop2017hierarchical, jin2010hierarchical}, the number of iterations for the sampling method is set to $2\times 10^5$ and $4 \times 10^6$, respectively.
To ensure the computational accuracy, we generate $5 \times 10^5$ samples for the parameter $u$ in this article.
Our focus is on the trace of the pCN method, as shown in sub-figure (a) of Figure $\ref{fig:trace}$.
We see that the whole sampling procedure completely explores the entire sample space.
And in sub-figure (b) of Figure \ref{fig:trace}, we plot the logarithm of the eigenvalues \( \{ c_k \}^n_{k=1} \) of the prior measure \( \mu_0 \).
The horizontal red dashed line denotes the corresponding eigenvalue \( c_M \) satisfying \( \frac{c_M}{c_1} < C_M = 10^{-3} \), where the truncated number \( M= 97 \).

\begin{figure}
	\centering
	\includegraphics[ keepaspectratio=true, width=0.6\textwidth,  clip=true, trim=50pt 2pt 80pt 30pt]{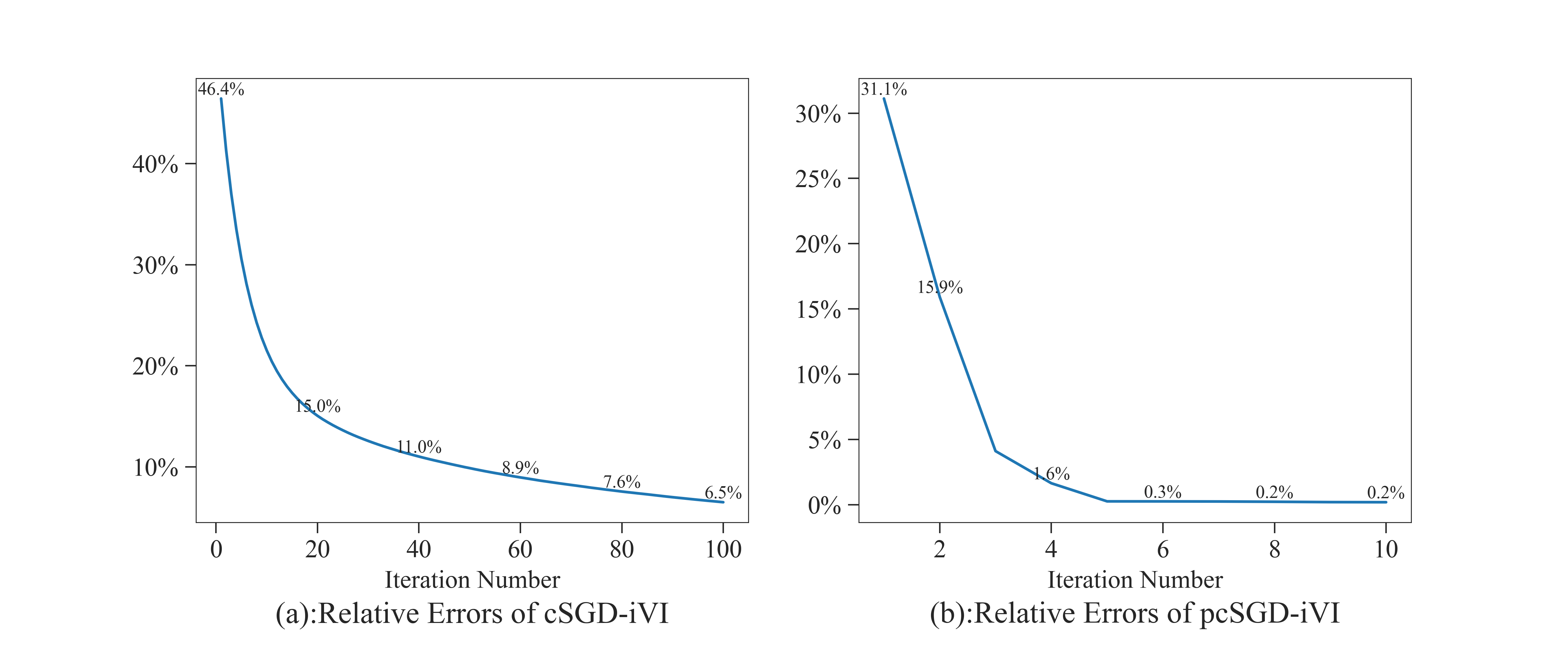}

	\caption{\emph{\small 
			(a): Relative errors of estimated posterior means of $u$ in the $L^2-norm$ of cSGD-iVI method with size of mesh $n = 100$; 
			(b): Relative errors of estimated posterior means of $u$ in the $L^2-norm$ of pcSGD-iVI method with size of mesh $n = 100$.
	}}
	
	\label{fig:Errors}
\end{figure}

Secondly, let us discuss cSGD-iVI method.
For calculating the gradient term in each iteration step, we need to solve one adjoint PDE (corresponding to computing $H^{*}$), and one forward PDE (corresponding to computing $H$).
Then each iteration we need to calculate $2$ PDEs.
Based on Algorithm $\ref{alg A}$, the total iteration number is $M\times J$, where $M$ is the number of samples, and $J$ is the number of averaging steps per sample.
In sub-figure (a) of Figure $\ref{fig:Errors}$ , we see that the relative errors between samples and background truth becomes small, and the convergence slows down at $60$ step. 
Hence we set $M = 100$ and $J=20$.
In summary, it is required to calculate $2MJ$ (each step to calculate gradient of $u$) $= 4000$ PDEs during the procedure.

At last, let us discuss pcSGD-iVI method.
At each step, we need to compute $T\mathcal{G}(u_{k})$, where $2$ PDEs are required to calculate $\mathcal{G}(u_{k})$, and $2N_{\text{ite}}$ PDEs to solve $T\mathcal{G}(u_{k})$ for accuracy.
We currently choose $N_{\text{ite}}$ as $10$.
Thus it is required to calculate $2 + 2N_{\text{ite}}=22$ PDEs in each iteration step.
In sub-figure (b) of Figure $\ref{fig:Errors}$, we see that the relative errors between samples and background truth start to converge at $3$ step, and the descending speed slows down. 
Hence we choose $M = 15$ and $J=20$.
In summary, it is required to calculate $22MJ$ (each step to calculate gradient of $u$) $= 6600$PDEs during the procedure.

As a result, we need to calculate at most $4000$ PDEs for employing cSGD-iVI method, while $6600$ PDEs for employing pcSGD-iVI method.
On the other hand, for the pCN method, $5 \times 10^5$ PDEs are needed to be calculated.
We conclude that the computational cost of these two methods are much less than the cost of pCN method.

Next, we illustrate the numerical results obtained by both methods, and compare them to the results of pCN.
As \cite{dashti2013bayesian} shows, pCN method provides a good estimate of the posterior measure of $u$.
It is necessary for us to compare the estimated posterior measures of $u$ obtined by these three methods.
Then we provide some discussions about the estimated posterior mean and covariance functions.

Let us provide some discussions of estimated posterior mean function.
As shown in sub-figure (a) of Figure $\ref{fig:Comparison_cSGD}$, the mean function of estimated posterior measure of $u$ obtained by the cSGD-iVI method and the background truth are drawn in black solid and red dashed lines, respectively.
Two green lines represent the upper and lower bounds of the $95 \%$ credibility region of the estimated posterior mean function.
We see that the $95 \%$ credibility region includes the background truth in most of the region of $u$; however on the left part, credibility region does not.
This illustrates that the estimated posterior measure derived by cSGD-iVI could not fully reflect the uncertainty of $u$ completely.
In sub-figure (a) of Figure $\ref{fig:Comparison_pcSGD}$, the mean function of estimated posterior measure of $u$ obtained by the pcSGD-iVI method and the background truth are drawn in black solid and red dashed lines, respectively.
The $95 \%$ credibility region includes the background truth, which means the estimated posterior measure obtained by pcSGD-iVI reflects the uncertainty of $u$.

\begin{figure}
	\centering
	\includegraphics[ keepaspectratio=true, width=0.8\textwidth,  clip=true, trim=120pt 2pt 100pt 20pt]{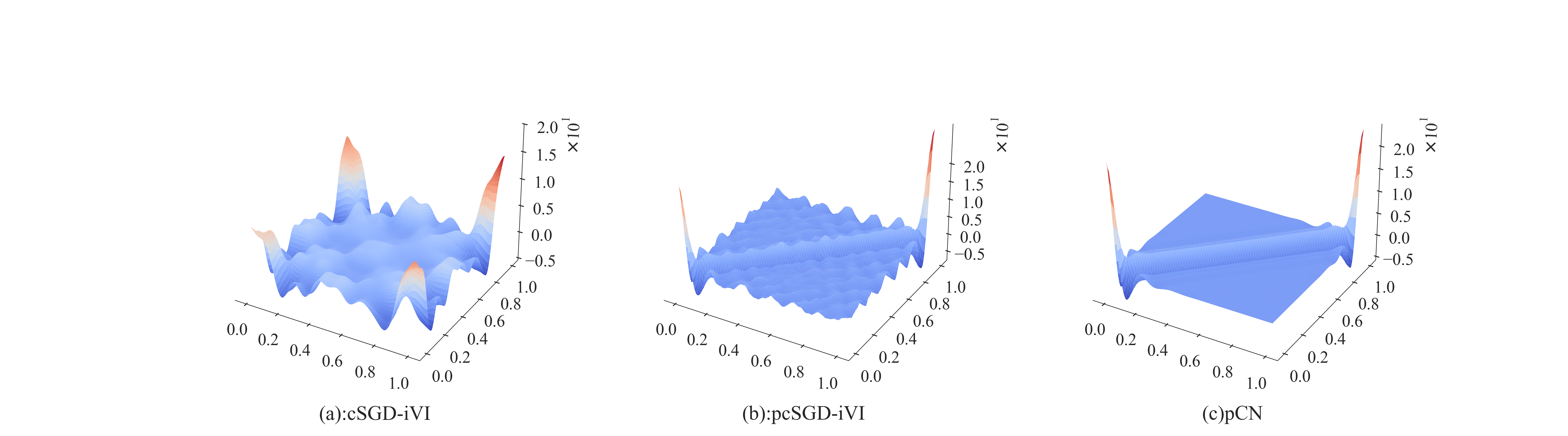}

	\caption{\emph{\small 
			The comparison of estimated posterior covariance operators of $u$ with mesh size $n=100$.
			(a): The covariance operator obtained by cSGD-iVI method; 
			(b): The covariance operator obtained by pcSGD-iVI method;
			(c): The covariance operator obtained by pCN method.
	}}
	
	\label{fig:Covariance}
\end{figure}

In Figure $\ref{fig:Errors}$, we show the relative errors of both methods, which is defined by
\begin{align}\label{equ:relative}
	\text{relative error} = \frac{\lVert u_{k, J} - u^{\dagger}\rVert^2_2}{\lVert u^{\dagger}\rVert^2_2},
\end{align}
where $J$ is the averaging converging steps defined in Algorithms $\ref{alg A}$ and $\ref{alg B}$.
As shown in sub-figure (a) of Figure $\ref{fig:Errors}$, the relative error curve illustrates that the iteration process of cSGD-iVI converges within $100$ steps, and is stable around $7 \%$ at the end of the iteration.
The convergence speed is fast since the descending trend is rapid at first $20$ steps.
As shown in sub-figure (b) of Figure $\ref{fig:Errors}$, the relative error curve illustrates that the iteration process of pcSGD-iVI converges within $10$ steps, and is stable around $0.2 \%$ at the end of the iteration.
Based on the visual (sub-figure (a) of Figures $\ref{fig:Comparison_cSGD}$ and $\ref{fig:Comparison_pcSGD}$) and quantitative (relative errors shown in Figure $\ref{fig:Errors}$) evidence,  we conclude that the pcSGD-iVI method provides estimated posterior mean functions for the parameter \(u\) that closely match the background truth; however, the cSGD-iVI mean function is less accurate near the left and right boundaries compared to pcSGD-iVI.

\begin{figure}
	\centering
	\includegraphics[ keepaspectratio=true, width=0.6\textwidth,  clip=true, trim=120pt 2pt 80pt 30pt]{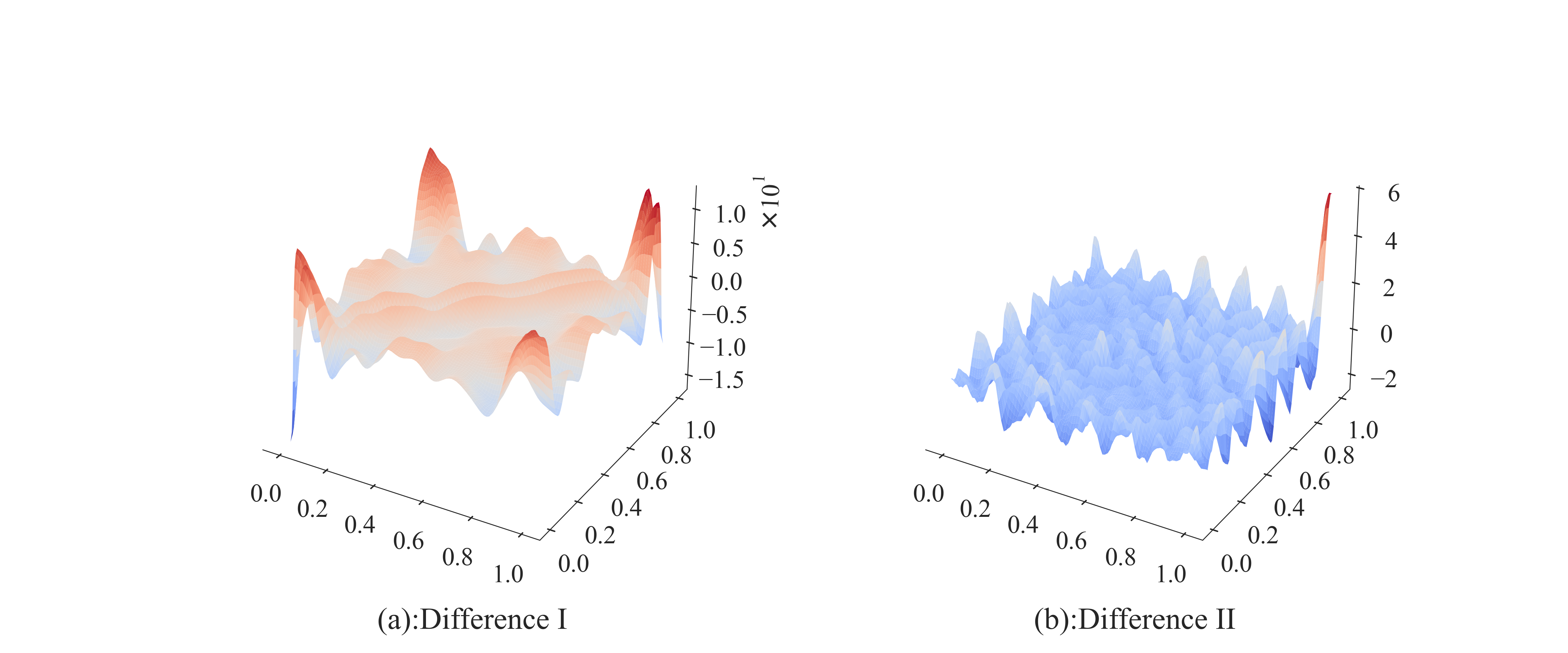}

	\caption{\emph{\small 
			(a): The difference of covariance operator, obtained by cSGD-iVI and pCN methods; 
			(b): The difference of covariance operator, obtained by pcSGD-iVI and pCN methods.
	}}
	
	\label{fig:Covariancediff}
\end{figure}

Furthermore, we provide numerical evidence to support a comparison of the estimated posterior mean functions obtained by both methods with the pCN method. 
The comparison of estimated posterior mean functions obtained by cSGD-iVI (blue solid line) and pCN (orange dashed line) is drawn in sub-figure (b) of Figure $\ref{fig:Comparison_cSGD}$. 
The differences are visually small in the middle region but become larger near the left, right, and corner areas. 
The relative error between them is given by
\begin{align}\label{equ:relate_cSGD}
	\text{relative error} := \frac{\lVert u^{*}_{\text{cSGD}} - u^{*}_{\text{pCN}} \rVert^2_{L^2}}{\lVert u^{*}_{\text{pCN}} \rVert^2_{L^2}} = 0.04755.
\end{align}
On the other hand, the comparison of pcSGD-iVI (blue solid line) and pCN (orange dashed line) is drawn in sub-figure (b) of Figure $\ref{fig:Comparison_pcSGD}$.
The differences are also small in the central, left, and right parts of the domain but slightly more pronounced in the corners.
The relative error between them is given by
\begin{align}\label{equ:relate_pcSGD}
	\text{relative error} := \frac{\lVert u^{*}_{\text{pcSGD}} - u^{*}_{\text{pCN}} \rVert^2_{L^2}}{\lVert u^{*}_{\text{pCN}} \rVert^2_{L^2}} = 0.005031.
\end{align}
Hence, the estimated posterior mean function derived by pcSGD-iVI are quantitatively similar to that by pCN method; however the estimated posterior mean obtained from cSGD-iVI is worse.
As a result, combining visual (sub-figure (b) of Figures $\ref{fig:Comparison_cSGD}$ and $\ref{fig:Comparison_pcSGD}$) and quantitative (relative errors given in ($\ref{equ:relate_cSGD}$) and ($\ref{equ:relate_pcSGD}$)) evidence, we conclude that pcSGD-iVI provides a better estimate of the posterior mean function than the cSGD-iVI method.

Next, we provide some discussions of estimated posterior covariance functions.
For the numerical convenience, we compare the estimated posterior covariance matrices, variances, and covariance functions.
The definitions and notations of covariance matrix, variance function and covariance function of the posterior covariance are provided in Subsection { A} of supplementary material.

\begin{figure}
	\centering
	\includegraphics[ keepaspectratio=true, width=0.8\textwidth,  clip=true, trim=120pt 2pt 120pt 30pt]{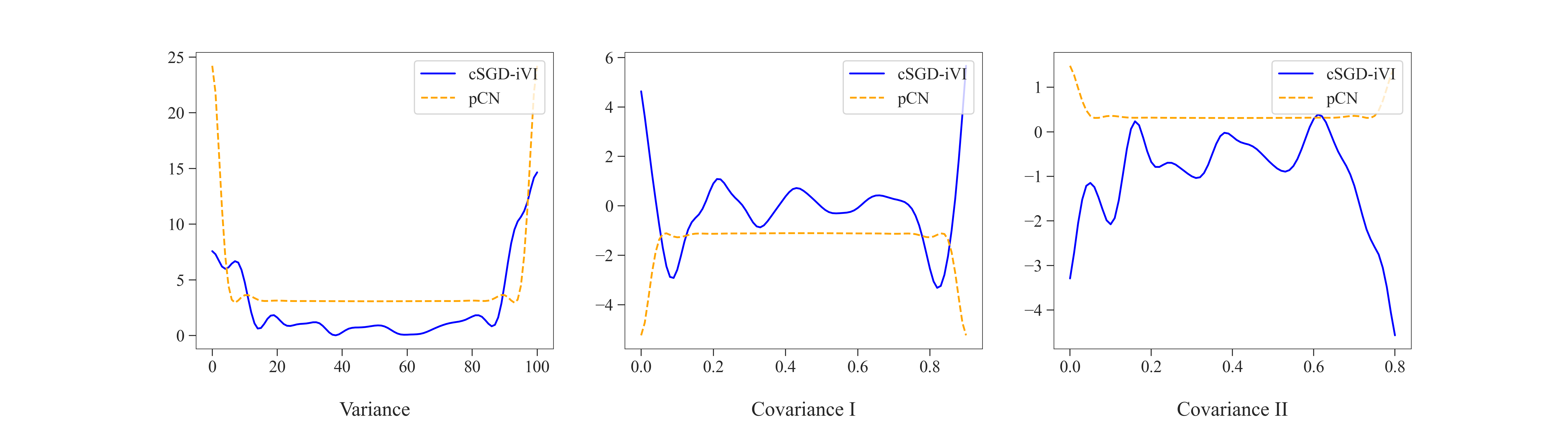}

	\caption{\emph{\small 
			Comparison of variance, covariance functions obtained by cSGD-iVI and pCN methods.
			(a): The variance function $\lbrace c(x_i, x_i)\rbrace^{100}_{i=1}$ on all the mesh point pairs $\lbrace (x_i, x_{i})\rbrace^{100}_{i=1}$;
			(b): The covariance function $\lbrace c(x_i, x_{i+10}) \rbrace^{90}_{i=1}$ on the mesh points $\lbrace (x_i, x_{i+10}) \rbrace^{90}_{i=1}$;
			(c): The covariance function $\lbrace c(x_i, x_{i+20}) \rbrace^{80}_{i=1}$ on the mesh points $\lbrace (x_i, x_{i+20}) \rbrace^{80}_{i=1}$.
	}}
	
	\label{fig:Variances_cSGD}
\end{figure}

In Tables $\ref{table:relative_cSGD}$ and $\ref{table:relative_pcSGD}$, we compare the variance, covariance functions of cSGD-iVI and pCN methods, pcSGD-iVI and pCN methods, respectively.
We see that the relative errors shown in Table $\ref{table:relative_cSGD}$ are larger than that shown in Table $\ref{table:relative_pcSGD}$.
In both tables, the relative errors for the variance functions are significantly smaller than those for the covariance functions. 
For pcSGD-iVI, the covariance matrices, variance functions, and covariance functions are quantitatively similar to those produced by the pCN method. 
Conversely, the covariance matrices, variances, and covariance functions derived by cSGD-iVI differ more notably from those obtained by pCN, particularly in the covariance functions. 
Quantitatively, this suggests that pcSGD-iVI yields a closer approximation to the posterior covariance operator than cSGD-iVI.

\begin{figure}
	\centering
	\includegraphics[ keepaspectratio=true,width=0.8\textwidth,  clip=true, trim=120pt 2pt 120pt 30pt]{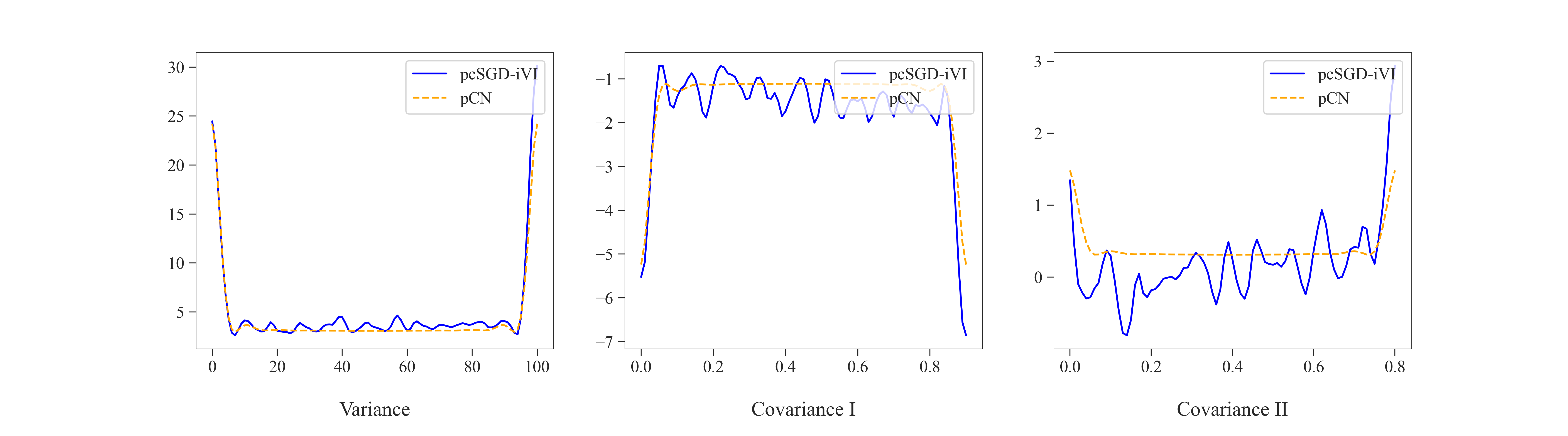}

	\caption{\emph{\small 
			Comparison of variance, covariance functions obtained by pcSGD-iVI and pCN methods.
			(a): The variance function $\lbrace c(x_i, x_i)\rbrace^{100}_{i=1}$ on all the mesh point pairs $\lbrace (x_i, x_{i})\rbrace^{100}_{i=1}$;
			(b): The covariance function $\lbrace c(x_i, x_{i+10}) \rbrace^{90}_{i=1}$ on the mesh points $\lbrace (x_i, x_{i+10}) \rbrace^{90}_{i=1}$;
			(c): The covariance function $\lbrace c(x_i, x_{i+40}) \rbrace^{80}_{i=1}$ on the mesh points $\lbrace (x_i, x_{i+20}) \rbrace^{20}_{i=1}$.
	}}
	
	\label{fig:Variances_pcSGD}
\end{figure}

In Figure $\ref{fig:Covariance}$, we draw the covariance matrix $\bm{c}_{\text{cSGD}}$, $\bm{c}_{\text{pcSGD}}$, and $\bm{c}_{\text{pCN}}$ in the sub-figures (a), (b), and (c).
In Figure $\ref{fig:Covariancediff}$, we draw the differences between $\bm{c}_{\text{cSGD}}$ and $\bm{c}_{\text{pCN}}$ in sub-figure (a); the differences between $\bm{c}_{\text{pcSGD}}$ and $\bm{c}_{\text{pCN}}$ in sub-figure (b), respectively.
Combining Figures $\ref{fig:Covariance}$ and $\ref{fig:Covariancediff}$, we see that the covariance operator obtained by pcSGD-iVI is more similar to that of the pCN method, and the differences are visually small; whereas the covariance operator derived from cSGD-iVI deviates more significantly.
Furthermore, we provide a detailed comparison of the variance and covariance functions obtained by these three method in Figures $\ref{fig:Variances_cSGD}$ and $\ref{fig:Variances_pcSGD}$, respectively.
In all the sub-figures of Figures $\ref{fig:Variances_cSGD}$ and $\ref{fig:Variances_pcSGD}$, the variance and covariance functions obtained by both iVI methods are shown in blue solid lines, while those from pCN are plotted in orange dashed lines.
In sub-figures (a) of Figures $\ref{fig:Variances_cSGD}$ and $\ref{fig:Variances_pcSGD}$, we show the variance function calculated on all the mesh point pairs $\lbrace (x_i, x_{i})\rbrace^{n}_{i=1}$ with $n = 100$.
In sub-figures (b) and (c), we show the covariance functions calculated on the pairs of points $\lbrace (x_i, x_{i+10})\rbrace^{n-10}_{i=1}$, and $\lbrace (x_i, x_{i+20})\rbrace^{n-20}_{i=1}$, respectively.

The covariance matrices, variance functions, and covariance functions obtained by cSGD-iVI show large deviations from those produced by pCN. 
In contrast, the differences between pcSGD-iVI and pCN are much smaller. 
The errors associated with cSGD-iVI are significantly larger than those of pcSGD-iVI.
This suggests that the covariance operator of cSGD-iVI diverges considerably from that of pCN, whereas the operator from pcSGD-iVI closely resembles that of pCN, consistent with the visual results. 
Furthermore, as is seen in the sub-figure (a) of Figure $\ref{fig:Comparison_pcSGD}$, the $95 \%$ credibility region of the estimated posterior mean contains the background truth, which indicates that the Bayesian setup derived by pcSGD-iVI is meaningful and in accordance with the frequentist theoretical investigations of the posterior consistency \cite{wang2019frequentist, zhang2020convergence}. 
Overall, we conclude that the pcSGD-iVI method provides a better estimate of the posterior covariance than cSGD-iVI, based on the visual (Figures $\ref{fig:Covariance}$, $\ref{fig:Covariancediff}$, $\ref{fig:Variances_cSGD}$ and $\ref{fig:Variances_pcSGD}$) and quantitative (Tables $\ref{table:relative_cSGD}$ and $\ref{table:relative_pcSGD}$) evidence.

Now we provide conclusions about the estimated posterior measures corresponding to parameter $u$ obtained by cSGD-iVI and pcSGD-iVI.
In Bayesian inference, an accurate approximation of the posterior distribution requires not only a good estimate of the mean function, but also a reliable estimate of the covariance, which quantifies the uncertainty of parameter $u$. 
As discussed previously, the posterior mean of 
\(u\) estimated by pcSGD-iVI closely matches that obtained by pCN, based on the visual (sub-figure (b) of Figure $\ref{fig:Comparison_pcSGD}$) and quantitative (relative error calculated in ($\ref{equ:relate_pcSGD}$)) evidence, moreover,  $95 \%$ credibility region of the estimated mean function includes background truth.
But the mean function of cSGD-iVI is worse based on the visual (sub-figure (b) of Figure $\ref{fig:Comparison_cSGD}$) and quantitative (relative error calculated in ($\ref{equ:relate_cSGD}$)) evidence, while $95 \%$ credibility region of the estimated mean function could not include background truth completely.
Based on the visual evidence (Figures $\ref{fig:Covariance}$ and $\ref{fig:Variances_pcSGD})$, the posterior covariance matrices, variance functions, and covariance functions derived from pcSGD-iVI and pCN appear visually similar.
And quantitative evidence (Table $\ref{table:relative_pcSGD}$) shows that the relative errors between the covariance matrices, variance, and covariance functions are quantitatively small, which indicates that the posterior covariance obtained by pcSGD-iVI method is as good as that by pCN method.
But the covariance operator obtained of cSGD-iVI is much different from that of pCN based on visual evidence (Figures $\ref{fig:Covariance}$ and $\ref{fig:Variances_cSGD})$, and quantitative evidence (Table $\ref{table:relative_cSGD}$) shows that relative errors are large.
These indicates that cSGD-iVI could not provide a reliable estimate of posterior covariance operator.
In summary, both visual and quantitative results indicate that pcSGD-iVI provides a more reliable approximation of the posterior mean function and covariance operator than cSGD-iVI.

{Furthermore, we show the performances and cost on different mesh sizes in Section B.1 of supplementary material.}


\section{Conclusion}\label{sec4}
In this paper, we develop cSGD-iVI and pcSGD-iVI methods in infinite-dimensional spaces, extending the finite-dimensional cSGD inference method proposed in \cite{Stephan2016PMLR, Stephan2017JMLR}. 
Based on Subsection $\ref{subsec2.2}$, we reformulate the cSGD method within a Bayesian inference framework. 
By solving the inference problem, we identify an optimal learning rate that minimizes the KL divergence between the estimated and true posterior measures, enabling sampling from the estimated posterior via the cSGD iteration.
We further introduce the pcSGD-iVI method, which builds upon the cSGD-iVI framework. 
We propose a learning rate that minimizes the KL divergence between the estimated and true posterior measures. 
Additionally, by treating the stochastic gradient as a random variable, we derive its corresponding probability distribution and thereby introduce the prior operator \( \mathcal{Q} \).
We also establish the regularization properties of the cSGD formulation and provide discretization error bounds between the estimated posterior mean, posterior samples, and the truth function, which depend on the learning rate and discretization level.
The proposed cSGD-iVI and pcSGD-iVI methods are applied to two inverse problems: a simple elliptic problem and a steady-state Darcy flow problem. For both cases, the pcSGD-iVI method yields estimated posterior measures for the parameter \( u \) that closely resemble the true posterior, accurately reflecting the uncertainty in \( u \), as supported by both visual and quantitative evidence. 
In contrast, the cSGD-iVI method produces less accurate posterior mean functions and covariance operators, especially in the Darcy flow problem.

Our current implementations of cSGD-iVI and pcSGD-iVI are based on linear inverse problems. 
For nonlinear problems, these methods can only address linearized formulations and produce approximate posterior measures. 
As a result, the estimated posterior may not be sufficiently accurate for highly nonlinear inverse problems.


\section*{Acknowledgments}
This work is supported by the NSFC grants 12322116, 12271428, 12326606, the Major projects of the NSFC grants 12090021, 12090020, UKRI and CSC.

\section*{Conflict of interest}
On behalf of all authors, the corresponding author states that there is no conflict of interest.

\bibliographystyle{amsplain}
\bibliography{references}
\end{document}


\abovedisplayskip=4.2pt plus 4.2pt
\abovedisplayshortskip=0pt plus 4.2pt
\belowdisplayskip=4.2pt plus 4.2pt
\belowdisplayshortskip=2.1pt plus 4.2pt
\begin{frontmatter}
	\title{Supplementary material of  ``Stochastic gradient descent based variational inference for infinite-dimensional inverse problems''
	}

	
	

\author{Jiaming Sui
\footnote{Email: sjming1997327@stu.xjtu.edu.cn, School of Mathematics and Statistics, Xi'an Jiaotong University, Xi'an, Shaanxi 710049, China},
Junxiong Jia 
\footnote{Email: jjx323@xjtu.edu.cn, School of Mathematics and Statistics, Xi'an Jiaotong University, Xi'an, Shaanxi 710049, China}, 
Jinglai Li
\footnote{Email: j.li.10@bham.ac.uk, School of Mathematics, University of Birmingham, Birmingham, B15 2TT, United Kingdom}}
	

	\begin{abstract}
		This supplementary material presents the definitions and notations introduced in Section 3. 
		Additionally, we provide the comprehensive proofs for the lemmas and theorems discussed in Section 2 of the main text.
		Next, it includes the numerical results of the simple elliptic equation on different meshes, and results of nonlinear inverse problems with the steady-state Darcy flow equation, as presented in Section 3 of the main text.
		At last, we introduce the Theorems 15 and 16 of \cite{dashti2013bayesian}.

	\end{abstract}
	
\end{frontmatter}

\section{Definitions and notations of the numerical results in Section 3}

\subsection{Definitions and notations of covariance functions in Subsection 3.1.2}\label{subsec:AppDef}

Before making the comparisons, we need to introduce the definition of covariance matrix, variance function and covariance function of the posterior covariance.
Consider a Gaussian random field $m$ on a domain $\Omega$ with mean $\bar{m}$ and the covariance function $c(\bm{x}, \bm{y})$ describing the covariance between $m(\bm{x})$ and $m(\bm{y})$:
\begin{align*}
	c(\bm{x}, \bm{y}) = \mathbb{E}((m(\bm{x})-\bar{m}(\bm{x}))(m(\bm{y})-\bar{m}(\bm{y}))), \quad \text{for} \ \bm{x}, \bm{y} \in \Omega.
\end{align*}
The corresponding covariance operator $\mathcal{C}$ is 
\begin{align*}
	(\mathcal{C}\phi^{\prime})(\bm{x}) = \int_{\Omega}c(\bm{x}, \bm{y})\phi^{\prime}(\bm{y})d\bm{y},
\end{align*}
where the functions $\phi^{\prime}$ are some sufficiently regular functions defined on $\Omega$.
To provide more numerical discussion of the posterior covariance, we introduce the matrix representation of the operator.
Let $V_h$ be the finite-dimensional subspace of $L^2(\Omega)$ originating from a finite element discretization with continuous Lagrange basis functions $\lbrace \phi_i\rbrace^n_{i=1}$, which correspond to the nodal points $\lbrace\bm{x}_j\rbrace^n_{j=1}$, such that $\phi_i(\bm{x}_j)=\delta_{ij}.$
The nodal vector of a function $m_h = \sum^n_{j=1}m_j\phi_j \in V_h$ is denoted by the boldface symbol $\bm{m} = (m_1, \cdots, m_n)^T$.
Inner products between nodal coefficient vectors are given by
\begin{align*}
	(\bm{m}_1, \bm{m}_2)_{\bm{M}} = \bm{m}^T_1\bm{M}\bm{m}_2,
\end{align*}
where the $ij$ element of mass matrix $\bm{M}$ is defined by
\begin{align*}
	M_{ij} = \int_{\Omega}\phi_i(\bm{x})\phi_j(\bm{x}) d\bm{x}.
\end{align*}
The matrix representation $\bm{c}$ of the operator $\mathcal{C}$ is given with respect to the Lagrange basis $\lbrace \phi_i\rbrace^n_{i=1}$ in $\mathbb{R}^n$ by
\begin{align*}
	\int_{\Omega} \phi_i\mathcal{C}\phi_j dx = \langle \bm{e}_i, \bm{c}\bm{e}_j\rangle_{\bm{M}},
\end{align*}
where $\bm{e}_j$ is the coordinate vector corresponding to the basis function $\phi_j$.

Let $\bm{\Phi(x)} = (\phi_1(\bm{x}), \cdots, \phi_n(\bm{x}))$, $\tilde{\bm{v}}_{k} = \bm{\mathcal{C}}^{1/2}_0\bm{v}_k$, and $\lbrace \xi_k, \bm{v}_k\rbrace^{n}_{k=1}$ be the eigenpairs of operator
$\bm{\mathcal{C}}^{1/2}_0\bm{H}^*\bm{\Gamma}^{-1}_{\text{noise}}\bm{H}\bm{\mathcal{C}}^{1/2}_0$.
In order to represent the posterior covariance operator $\mathcal{C}$, we need to calculate the posterior covariance field $c$ given by
\begin{align*}
	c(\bm{x}, \bm{y}) \approx c_{0}(\bm{x}, \bm{y}) - \sum^r_{k=1}d_k(\tilde{\bm{v}}_{kh}(\bm{x}))(\tilde{\bm{v}}_{kh}(\bm{y}))^T.
\end{align*}
Here $\tilde{\bm{v}}_{kh}(\bm{x}) = \bm{\Phi(x)}^T\tilde{\bm{v}}_{k}$, $d_k = \epsilon_i/\epsilon_i + 1$, $c_{0}(\bm{x}, \bm{y})$ is the covariance function of the prior covariance operator $\mathcal{C}_0 = (\text{I} - \alpha\Delta)^{-2}$, based on the setting in Subsection 3.1.1.
Hence, the variance function $\lbrace c(x_i, x_i)\rbrace^{n}_{i=1}$ of the covariance field $c(\bm{x}, \bm{y})$ is calculated on the pairs of points $\lbrace (x_i, x_{i})\rbrace^{n}_{i=1}$.
And the covariance function $\lbrace c(x_i, x_{i+k})\rbrace^{n-k}_{i=1}$ of the covariance field $c(\bm{x}, \bm{y})$ is calculated on the pairs of points $\lbrace (x_i, x_{i+k})\rbrace^{n-k}_{i=1}$, where $k$ is an integer.
Readers can seek more details in \cite{bui2013computational}.

Next, let us denote $\mathcal{N}(u^{*}_{\text{cSGD}}, \mathcal{C}_{\text{cSGD}})$, $\mathcal{N}(u^{*}_{\text{pcSGD}}, \mathcal{C}_{\text{pcSGD}})$, $\mathcal{N}(u^{*}_{\text{pCN}}, \mathcal{C}_{\text{pCN}})$ as the estimated posterior measures of $u$ obtained by cSGD-iVI, pcSGD-iVI and pCN methods, respectively.
We show the comparisons of the following covariance matrices, variance, and covariance functions:
\begin{itemize}
	\item $\bm{c}_{\text{cSGD}}$, $\bm{c}_{\text{pcSGD}}$ and $\bm{c}_{\text{pCN}}$, which are the matrix representations of the covariance operators $\mathcal{C}_{\text{cSGD}}$, $\mathcal{C}_{\text{pcSGD}}$ and $\mathcal{C}_{\text{pCN}}$, respectively;
	\item $\lbrace c_{\text{cCGD}}(x_i, x_i)\rbrace^{n}_{i=1}$, $\lbrace c_{\text{pcSGD}}(x_i, x_i)\rbrace^{n}_{i=1}$ and $\lbrace c_{\text{pCN}}(x_i, x_i)\rbrace^{n}_{i=1}$, which are the variance functions of the covariance operators $\mathcal{C}_{\text{cSGD}}$, $\mathcal{C}_{\text{pcSGD}}$ and $\mathcal{C}_{\text{pCN}}$, respectively;
	\item $\lbrace c_{\text{cSGD}}(x_i, x_{i+10}) \rbrace^{n-10}_{i=1}$, $\lbrace c_{\text{pcSGD}}(x_i, x_{i+10}) \rbrace^{n-10}_{i=1}$ and $\lbrace c_{\text{pCN}}(x_i, x_{i+10}) \rbrace^{n-10}_{i=1}$; 
	These are the covariance functions calculated on the pairs of points $\lbrace (x_i, x_{i+{ 10}})\rbrace^{n-10}_{i=1}$.
	\item $\lbrace c_{\text{cSGD}}(x_i, x_{i+20}) \rbrace^{n-20}_{i=1}$, $\lbrace c_{\text{pcSGD}}(x_i, x_{i+20}) \rbrace^{n-20}_{i=1}$ and $\lbrace c_{\text{pCN}}(x_i, x_{i+20}) \rbrace^{n-20}_{i=1}$;
	These are the covariance functions calculated on the pairs of points $\lbrace (x_i, x_{i+{ 20}})\rbrace^{n-20}_{i=1}$.
\end{itemize}
We calculate the relative errors according to these functions with the following definitions:
\begin{itemize}
	\item The relative errors between $\bm{c}_{\text{cSGD}}$, $\bm{c}_{\text{pcSGD}}$, and $\bm{c}_{\text{pCN}}$ are given by	
	\begin{align*}
		\text{relative error}_{\text{cSGD}} = \frac{\lVert \bm{c}_{\text{cSGD}} - \bm{c}_{\text{pCN}} \rVert^2_F}{\lVert \bm{c}_{\text{pCN}} 	\rVert^2_F}, \quad 
		\text{relative error}_{\text{pc{ SGD}}} = \frac{\lVert \bm{c}_{\text{pcSGD}} - \bm{c}_{\text{pCN}} \rVert^2_F}{\lVert \bm{c}_{\text{pCN}} \rVert^2_F},
	\end{align*}
	where $\lVert \cdot \rVert_F$ denotes the Frobenius-norm defined by
	\begin{align*}
		\lVert A \rVert_F = \bigg (\sum^{n}_{i, j=1} a_{ij}^2 \bigg )^{1/2}.
	\end{align*}
	\item The relative errors between the variance functions and the covariance functions are given by:
	\begin{align*}
		\text{relative error of variance}_{\text{cSGD}} &= \frac{\sum^{n}_{i=1}\big (c_{\text{cSGD}}(x_i, x_i) - c_{\text{pCN}}(x_i, x_i)\big )^2}{\sum^{n}_{i=1}\big (c_{\text{pCN}}(x_i, x_i)\big)^2}, \\
		\text{relative error of covariance}_{\text{cSGD}} &= \frac{\sum^{n-k}_{i=1}\big(c_{\text{cSGD}}(x_i, x_{i+k}) - c_{\text{pCN}}(x_i, x_{i+k})\big)^2}{\sum^{n-k}_{i=1}\big(c_{\text{pCN}}(x_i, x_{i+k})\big)^2}; \\
		\text{relative error of variance}_{\text{pcSGD}} &= \frac{\sum^{n}_{i=1}\big (c_{\text{pcSGD}}(x_i, x_i) - c_{\text{pCN}}(x_i, x_i)\big )^2}{\sum^{n}_{i=1}\big (c_{\text{pCN}}(x_i, x_i)\big)^2}, \\
		\text{relative error of covariance}_{\text{pcSGD}} &= \frac{\sum^{n-k}_{i=1}\big(c_{\text{pcSGD}}(x_i, x_{i+k}) - c_{\text{pCN}}(x_i, x_{i+k})\big)^2}{\sum^{n-k}_{i=1}\big(c_{\text{pCN}}(x_i, x_{i+k})\big)^2}
	\end{align*}
	and $k$ is an integer.
\end{itemize}

\section{Numerical experiments}
{ 
\subsection{The simple elliptic equation with Gaussian noise}
\subsubsection{Performance and cost on different mesh sizes}
In order to provide more evidence, we employ the mesh size \(= \{200, 300\}\) for comparsions.
And we provide the comparsions of computational cost on different meshes \(= \{100, 200, 300\}\).

\subsubsection*{Mesh size \(= 200\)}
Let us provide some discussions of estimated posterior mean function.
As shown in sub-figure (a) of Figure $\ref{fig:Comparison_SGD200}$, the mean function of estimated posterior measure of $u$ obtained by the cSGD-iVI method and the background truth are drawn in black solid and red dashed lines, respectively.
Two green lines represent the upper and lower bounds of the $95 \%$ credibility region of the estimated posterior mean function.
We see that the $95 \%$ credibility region includes the background truth in most of the region of $u$; however on the left part, credibility region does not.
This illustrates that the estimated posterior measure derived by cSGD-iVI could not fully reflect the uncertainty of $u$ completely.
In sub-figure (b) of Figure $\ref{fig:Comparison_SGD200}$, the mean function of estimated posterior measure of $u$ obtained by the pcSGD-iVI method and the background truth are drawn in black solid and red dashed lines, respectively.
The $95 \%$ credibility region includes the background truth, which means the estimated posterior measure obtained by pcSGD-iVI reflects the uncertainty of $u$.

In Figure $\ref{fig:Errors200}$, we show the relative errors of both methods.
As shown in sub-figure (a) of Figure $\ref{fig:Errors200}$, the relative error curve illustrates that the iteration process of cSGD-iVI converges within $100$ steps, and is stable around $7 \%$ at the end of the iteration.
The convergence speed is fast since the descending trend is rapid at first $20$ steps.
As shown in sub-figure (b) of Figure $\ref{fig:Errors200}$, the relative error curve illustrates that the iteration process of pcSGD-iVI converges within $10$ steps, and is stable around $0.2 \%$ at the end of the iteration.
Based on the visual (Figure $\ref{fig:Comparison_SGD200}$) and quantitative (relative errors shown in Figure $\ref{fig:Errors200}$) evidence,  we conclude that the pcSGD-iVI method provides estimated posterior mean functions for the parameter \(u\) that closely match the background truth; however, the cSGD-iVI mean function is less accurate near the left and right boundaries compared to pcSGD-iVI.
In Figure $\ref{fig:Covariance200}$, we draw the covariance matrix $\bm{c}_{\text{cSGD}}$ and $\bm{c}_{\text{pcSGD}}$ in the sub-figures (a) and (b).
Combining the credibility region shown in Figure \ref{fig:Comparison_SGD200} of cSGD-iVI and pcSGD-iVI methods.
For cSGD-iVI method (sub-figure(a)), we see that the credibility region doesn't contain background truth completely, particularly in left and right part of the region.
This indicates that cSGD-iVI can not reflect the uncertainties of parameter $u$.
On the other hand, for pcSGD-iVI method (sub-figure(b)), we see that the credibility region contains background truth completely.
This indicates that with the help of preconditioning operator $T$, the estimated posterior measure obtained by pcSGD-iVI reflects the uncertainties of $u$ accurately.

In summary, we see that pcSGD-iVI provides a good estimate of mean function and covariance operator, based on visual (Figures $\ref{fig:Comparison_SGD200}$ and $\ref{fig:Covariance200}$) and quantitative (sub-figure (b) of Figure $\ref{fig:Errors200}$) evidence; while the estimation of cSGD-iVI is much worse, based on visual (Figures $\ref{fig:Comparison_SGD200}$ and $\ref{fig:Covariance200}$) and quantitative (sub-figure (a) of Figure $\ref{fig:Errors200}$) evidence.

\begin{figure}
	\centering
	\subfloat[Comparison of \(u\) I]{
	\includegraphics[ keepaspectratio=true, width=0.30\textwidth,  clip=true, trim=80pt 25pt 450pt 30pt]{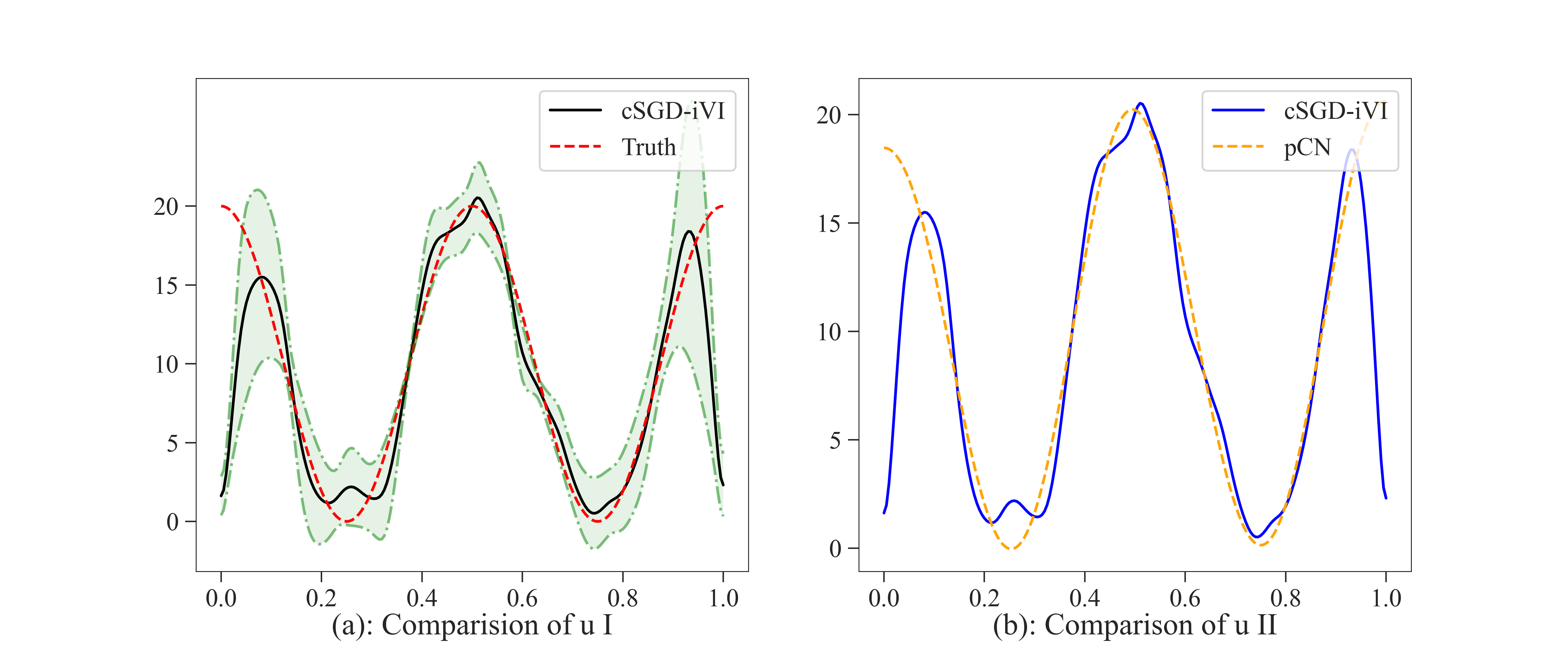}}
	\subfloat[Comparison of \(u\) II]{
	\includegraphics[ keepaspectratio=true, width=0.30\textwidth,  clip=true, trim=80pt 25pt 450pt 30pt]{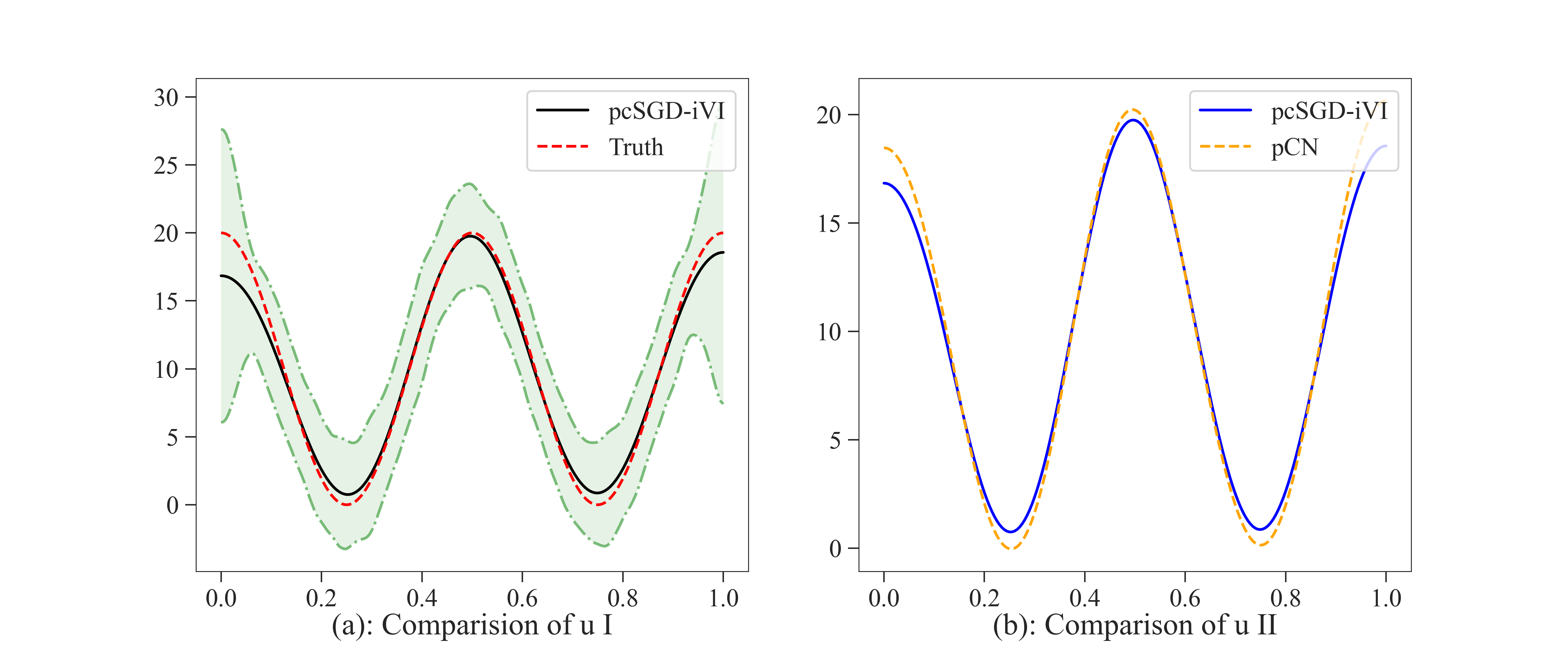}}

	\caption{\emph{\small 
			The green shade area represents the $95 \%$ credibility region of estimated posterior mean function; 
			(a): The comparison of the {estimated posterior mean function} obtained by cSGD-iVI and the background truth of $u$, respectively;
			(b): The comparison of the {estimated posterior mean function} obtained by pcSGD-iVI and the background truth of $u$, respectively.
	}}
	
	\label{fig:Comparison_SGD200}
	
\end{figure}

\begin{figure}
	\centering
	\includegraphics[ keepaspectratio=true, width=0.6\textwidth,  clip=true, trim=50pt 2pt 80pt 30pt]{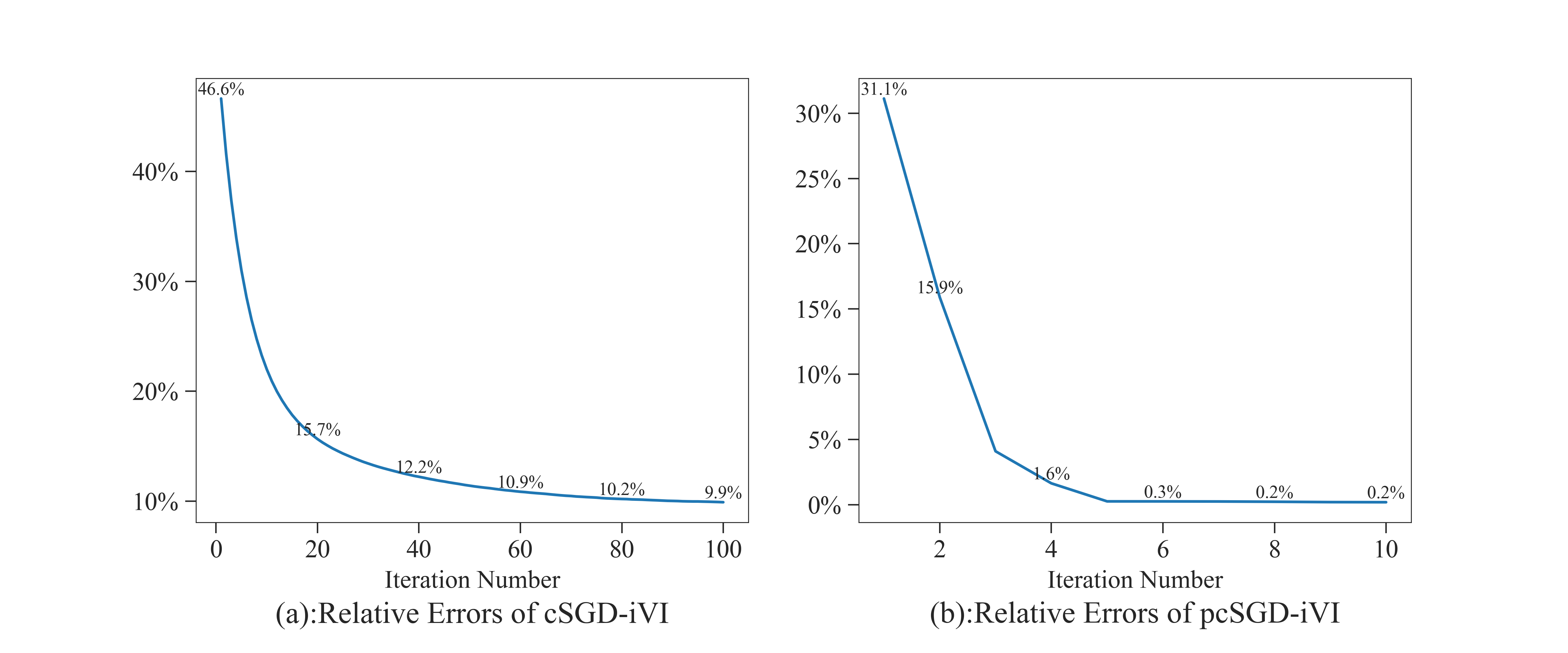}

	\caption{\emph{\small 
			(a): Relative errors of estimated posterior means of $u$ in the $L^2-norm$ of cSGD-iVI method with size of mesh $n = 200$; 
			(b): Relative errors of estimated posterior means of $u$ in the $L^2-norm$ of pcSGD-iVI method with size of mesh $n = 200$.
	}}
	
	\label{fig:Errors200}
\end{figure}

\begin{figure}
	\centering
	\includegraphics[ keepaspectratio=true, width=0.8\textwidth,  clip=true, trim=120pt 2pt 450pt 20pt]{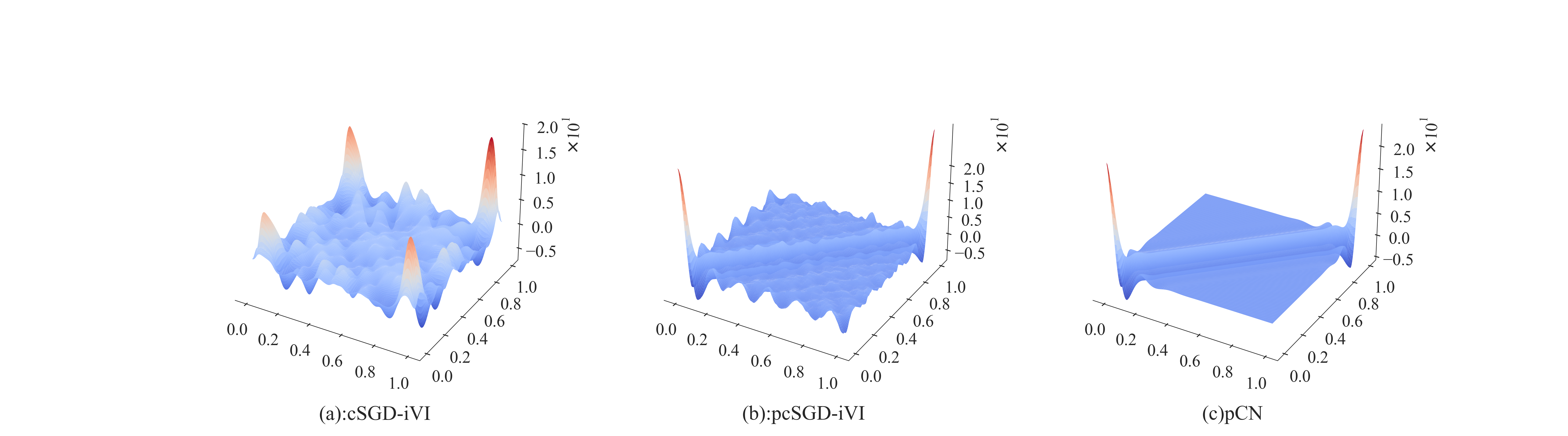}

	\caption{\emph{\small 
			The comparison of estimated posterior covariance operators of $u$ with mesh size $n=200$.
			(a): The covariance operator obtained by cSGD-iVI method; 
			(b): The covariance operator obtained by pcSGD-iVI method.
	}}
	
	\label{fig:Covariance200}
\end{figure}

\subsubsection*{Mesh size \(= 300\)}
Let us provide some discussions of estimated posterior mean function.
As shown in sub-figure (a) of Figure $\ref{fig:Comparison_SGD300}$, the mean function of estimated posterior measure of $u$ obtained by the cSGD-iVI method and the background truth are drawn in black solid and red dashed lines, respectively.
Two green lines represent the upper and lower bounds of the $95 \%$ credibility region of the estimated posterior mean function.
We see that the $95 \%$ credibility region includes the background truth in most of the region of $u$; however on the left part, credibility region does not.
This illustrates that the estimated posterior measure derived by cSGD-iVI could not fully reflect the uncertainty of $u$ completely.
In sub-figure (b) of Figure $\ref{fig:Comparison_SGD300}$, the mean function of estimated posterior measure of $u$ obtained by the pcSGD-iVI method and the background truth are drawn in black solid and red dashed lines, respectively.
The $95 \%$ credibility region includes the background truth, which means the estimated posterior measure obtained by pcSGD-iVI reflects the uncertainty of $u$.

In Figure $\ref{fig:Errors300}$, we show the relative errors of both methods.
As shown in sub-figure (a) of Figure $\ref{fig:Errors300}$, the relative error curve illustrates that the iteration process of cSGD-iVI converges within $100$ steps, and is stable around $7 \%$ at the end of the iteration.
The convergence speed is fast since the descending trend is rapid at first $20$ steps.
As shown in sub-figure (b) of Figure $\ref{fig:Errors300}$, the relative error curve illustrates that the iteration process of pcSGD-iVI converges within $10$ steps, and is stable around $0.2 \%$ at the end of the iteration.
Based on the visual (Figure $\ref{fig:Comparison_SGD300}$) and quantitative (relative errors shown in Figure $\ref{fig:Errors300}$) evidence,  we conclude that the pcSGD-iVI method provides estimated posterior mean functions for the parameter \(u\) that closely match the background truth; however, the cSGD-iVI mean function is less accurate near the left and right boundaries compared to pcSGD-iVI.
In Figure $\ref{fig:Covariance300}$, we draw the covariance matrix $\bm{c}_{\text{cSGD}}$ and $\bm{c}_{\text{pcSGD}}$ in the sub-figures (a) and (b).
Combining the credibility region shown in Figure \ref{fig:Comparison_SGD300} of cSGD-iVI and pcSGD-iVI methods.
For cSGD-iVI method (sub-figure(a)), we see that the credibility region doesn't contain background truth completely, particularly in left and right part of the region.
This indicates that cSGD-iVI can not reflect the uncertainties of parameter $u$.
On the other hand, for pcSGD-iVI method (sub-figure(b)), we see that the credibility region contains background truth completely.
This indicates that with the help of preconditioning operator $T$, the estimated posterior measure obtained by pcSGD-iVI reflects the uncertainties of $u$ accurately.

In summary, we see that pcSGD-iVI provides a good estimate of mean function and covariance operator, based on visual (Figures $\ref{fig:Comparison_SGD300}$ and $\ref{fig:Covariance300}$) and quantitative (sub-figure (b) of Figure $\ref{fig:Errors300}$) evidence; while the estimation of cSGD-iVI is much worse, based on visual (Figures $\ref{fig:Comparison_SGD300}$ and $\ref{fig:Covariance300}$) and quantitative (sub-figure (a) of Figure $\ref{fig:Errors300}$) evidence.

\begin{figure}
	\centering
	\subfloat[Comparison of \(u\) I]{
	\includegraphics[ keepaspectratio=true, width=0.30\textwidth,  clip=true, trim=80pt 25pt 450pt 30pt]{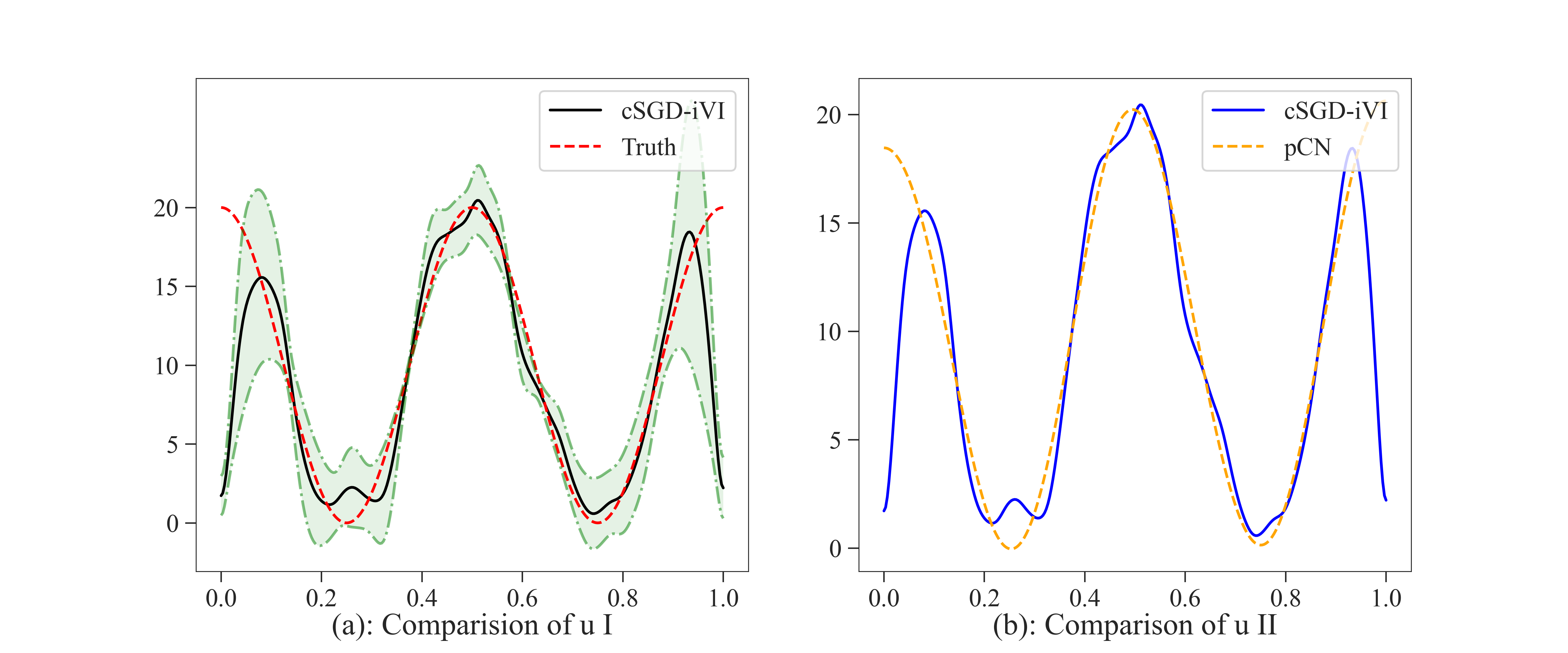}}
	\subfloat[Comparison of \(u\) II]{
	\includegraphics[ keepaspectratio=true, width=0.30\textwidth,  clip=true, trim=80pt 25pt 450pt 30pt]{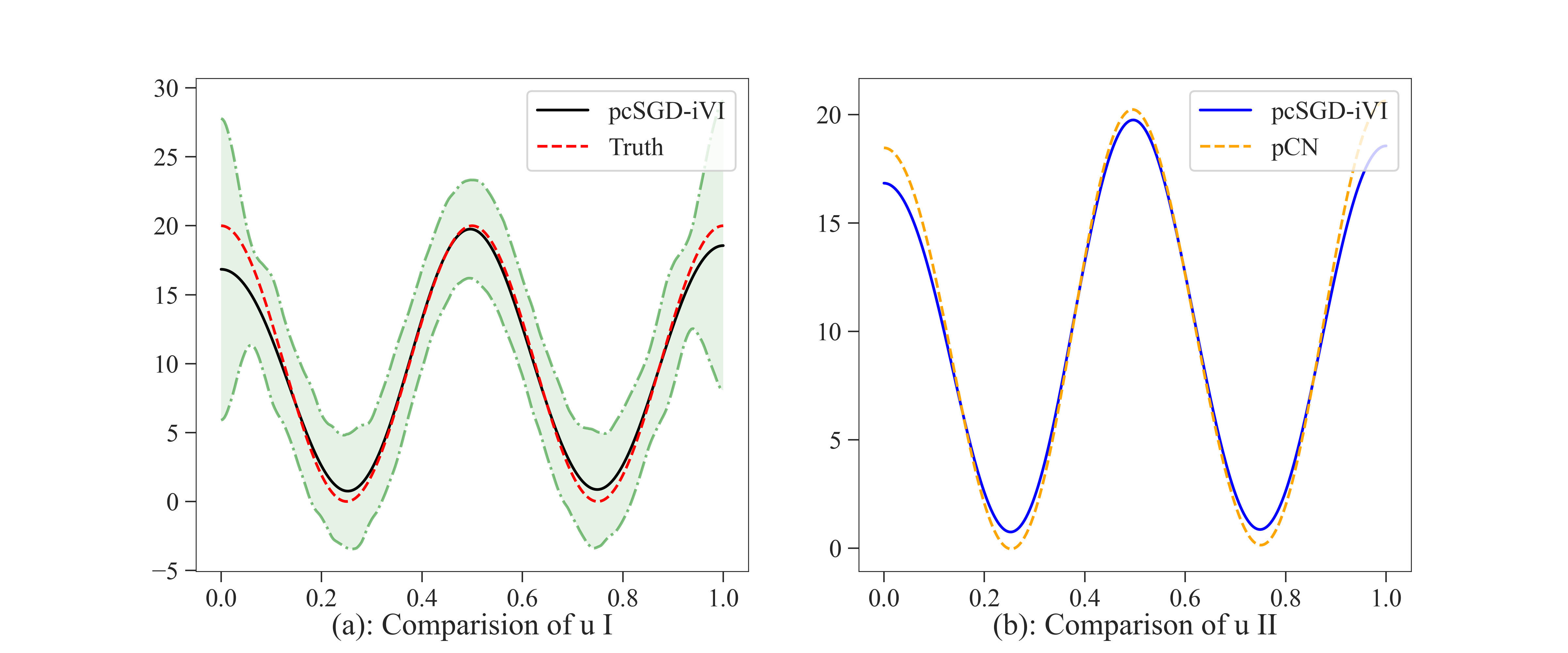}}

	\caption{\emph{\small 
			The green shade area represents the $95 \%$ credibility region of estimated posterior mean function; 
			(a): The comparison of the {estimated posterior mean function} obtained by cSGD-iVI and the background truth of $u$, respectively;
			(b): The comparison of the {estimated posterior mean function} obtained by pcSGD-iVI and the background truth of $u$, respectively.
	}}
	
	\label{fig:Comparison_SGD300}
	
\end{figure}

\begin{figure}
	\centering
	\includegraphics[ keepaspectratio=true, width=0.6\textwidth,  clip=true, trim=50pt 2pt 80pt 30pt]{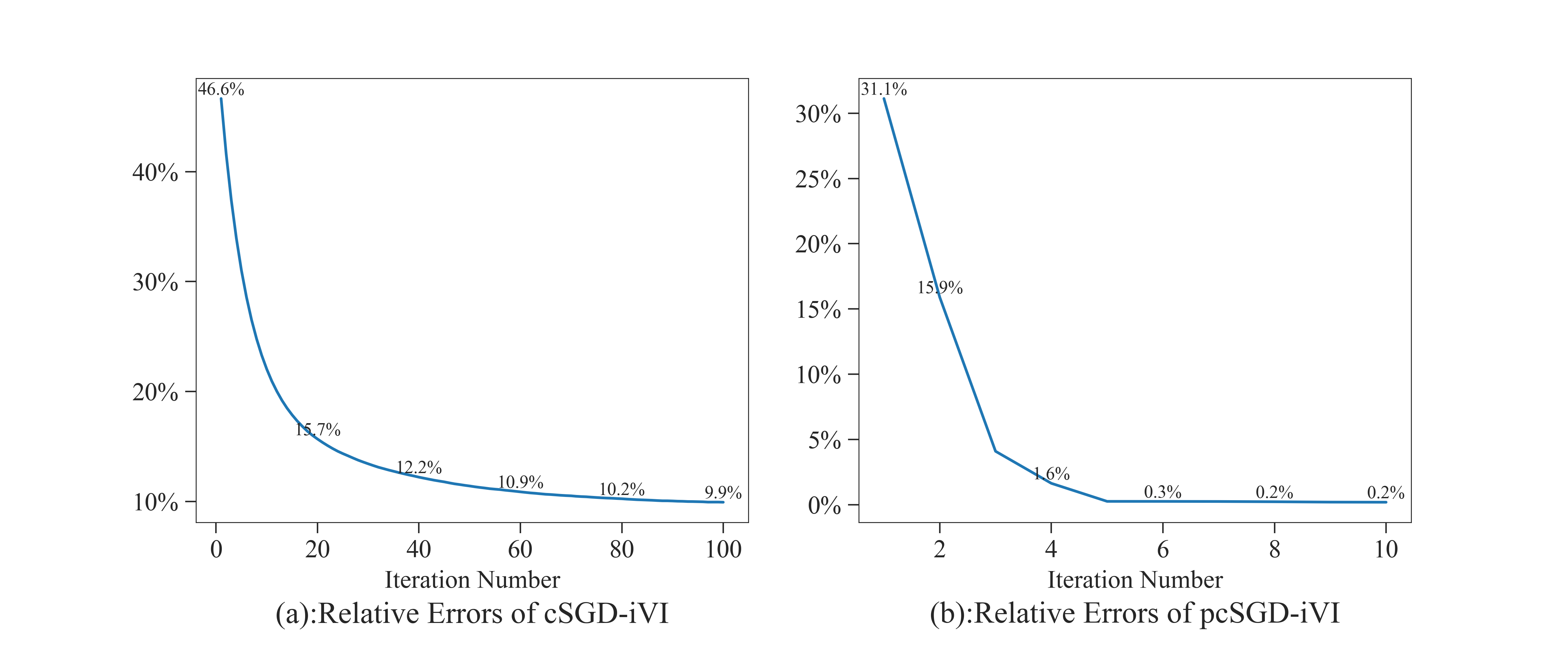}

	\caption{\emph{\small 
			(a): Relative errors of estimated posterior means of $u$ in the $L^2-norm$ of cSGD-iVI method with size of mesh $n = 300$; 
			(b): Relative errors of estimated posterior means of $u$ in the $L^2-norm$ of pcSGD-iVI method with size of mesh $n = 300$.
	}}
	
	\label{fig:Errors300}
\end{figure}

\begin{figure}
	\centering
	\includegraphics[ keepaspectratio=true, width=0.8\textwidth,  clip=true, trim=120pt 2pt 450pt 20pt]{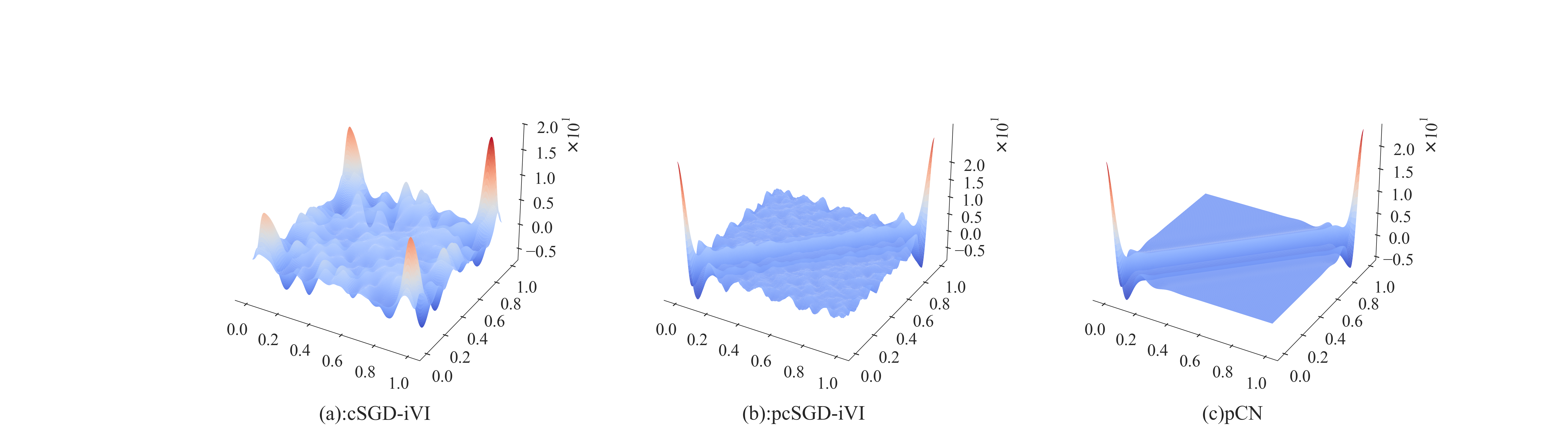}

	\caption{\emph{\small 
			The comparison of estimated posterior covariance operators of $u$ with mesh size $n=300$.
			(a): The covariance operator obtained by cSGD-iVI method; 
			(b): The covariance operator obtained by pcSGD-iVI method.
	}}
	
	\label{fig:Covariance300}
\end{figure}

\subsubsection*{Computaional cost}

\begin{table}
	\renewcommand{\arraystretch}{1.5}
	\centering
	\caption{\emph{\small Computational cost of cSGD-iVI}} \label{table:1d_cost_cSGD}
	\begin{tabular}{c|ccc}
			\hline $\text{Mesh size} $ & $n=100$  & $n=200$ & $n=300$   \\
			\hline $\text{Time (total)/s}$ & $1.9099$  & $2.5709$ & $3.0164$  \\
			\hline $\text{Time (per-step)/s}$ & $0.00955$ & $0.01285$ & $0.01508$ \\
			\hline
	\end{tabular}
\end{table}

\begin{table}
	\renewcommand{\arraystretch}{1.5}
	\centering
	\caption{\emph{\small Computational cost of pcSGD-iVI}} \label{table:1d_cost_pcSGD}
	\begin{tabular}{c|ccc}
			\hline $\text{Mesh size} $ & $n=100$  & $n=200$ & $n=300$   \\
			\hline $\text{Time (total)/s}$ & $6.5072$  & $7.6930$ & $9.2604$  \\
			\hline $\text{Time (per-step)/s}$ & $0.32536$ & $0.38465$ & $0.46302$ \\
			\hline
	\end{tabular}
\end{table}

At last, we show the computational cost for both methods on different meshes.
In Table \ref{table:Darcy_cost_cSGD}, we show the computational cost of cSGD-iVI for \(100\) steps.
And in Table \ref{table:Darcy_cost_pcSGD}, we show the computational cost of pcSGD-iVI for \(15\) steps.
For each step, the proceeding time of pcSGD-iVI is longer then cSGD-iVI, since we need to calculate more PDEs. 
Overall, the cost of pcSGD-iVI method is more expensive than cSGD-iVI method, which is consistent with our discussion of computational cost in Section 3.1.2 of the manuscript.

}

\begin{figure}
	\centering
	\includegraphics[ keepaspectratio=true, width=0.8\textwidth,  clip=true,  trim=120pt 2pt 120pt 30pt]{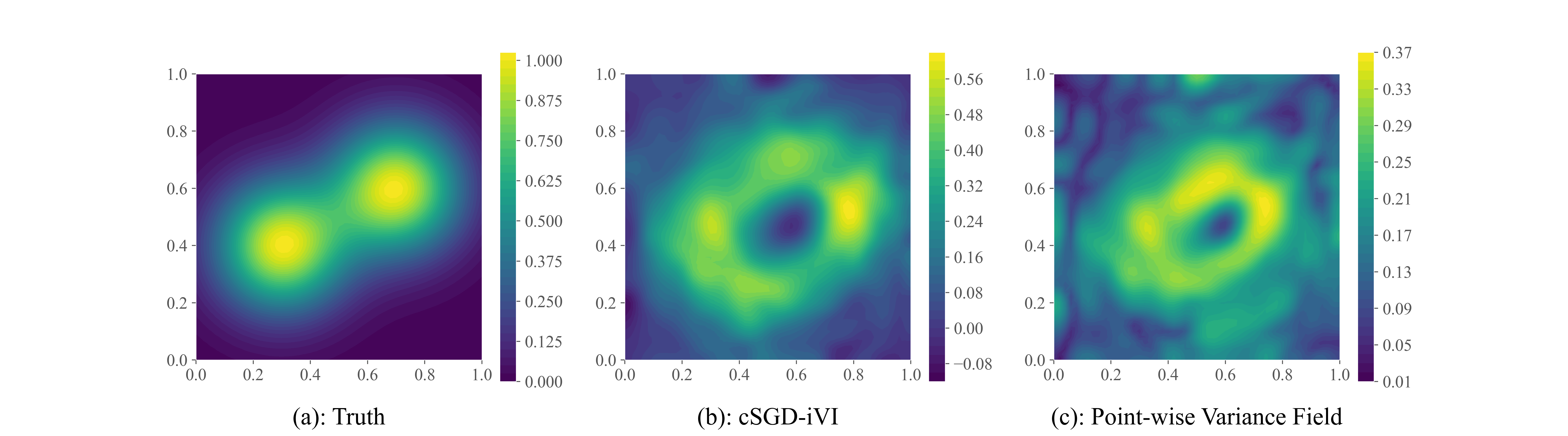}

	\caption{\emph{\small 
			(a): The background truth of $u$;
			(b): The estimated posterior mean function of $u$ obtained by cSGD-iVI;
			(c): The estimated variance function of posterior measure of $u$ obtained by cSGD-iVI.
	}}
	
	\label{fig:Comparison_cSGD_Darcy}
\end{figure}

\subsection{Non-linear inverse problem with steady-state Darcy flow equation}\label{subsec3.2}
\subsubsection{Basic settings}\label{subsec3.2.1}
In this section, we concentrate on the inverse problem of estimating the permeability distribution in a porous medium from a discrete set of pressure measurements, studied in \cite{calvetti2018iterative}. Consider the following steady-state Darcy flow equation:
\begin{align}
	\begin{split}
		-\nabla \cdot (e^u\nabla w) &= f \quad x \in \Omega, \\
		w &= 0 \quad x \in \partial \Omega,
	\end{split}
\end{align}
where $f \in H^{-1}(\Omega)$ is the source function,  $u \in \mathcal{X} := L^{\infty}(\Omega)$ called log-permeability for the computational area $\Omega = (0, 1)^2$, and denote $\bm{x} = (x, y) \in \mathbb{R}^2$. 
Then the forward operator has the following form:
\begin{align}
	Hu = (w(\bm{x}_1), w(\bm{x}_2), \cdots, w(\bm{x}_{N_d}))^T,
\end{align}
where $x_i \in \Omega$ for $i = 1, \cdots, N_d$.

Let us assume that the operator $\mathcal{F}$ is Fr\'{e}chet differentiable, and $w = \mathcal{F}(u)$. 
We linearize $\mathcal{F}(u)$ around $u^{*}$ to obtain
\begin{align}
	w \approx \mathcal{F}(u^{*}) + F^{\prime}(u - u^{*}),
\end{align}
where $F^{\prime}$ is the Fr\'{e}chet derivative of $\mathcal{F}(u)$ evaluated at $u^{*}$. Consequently, we transform the non-linear problem into the linear form, and we are able to employ our {cSGD-iVI and pcSGD-iVI} method.

\begin{figure}
	\centering
	\includegraphics[ keepaspectratio=true, width=0.8\textwidth,  clip=true,  trim=120pt 2pt 120pt 30pt]{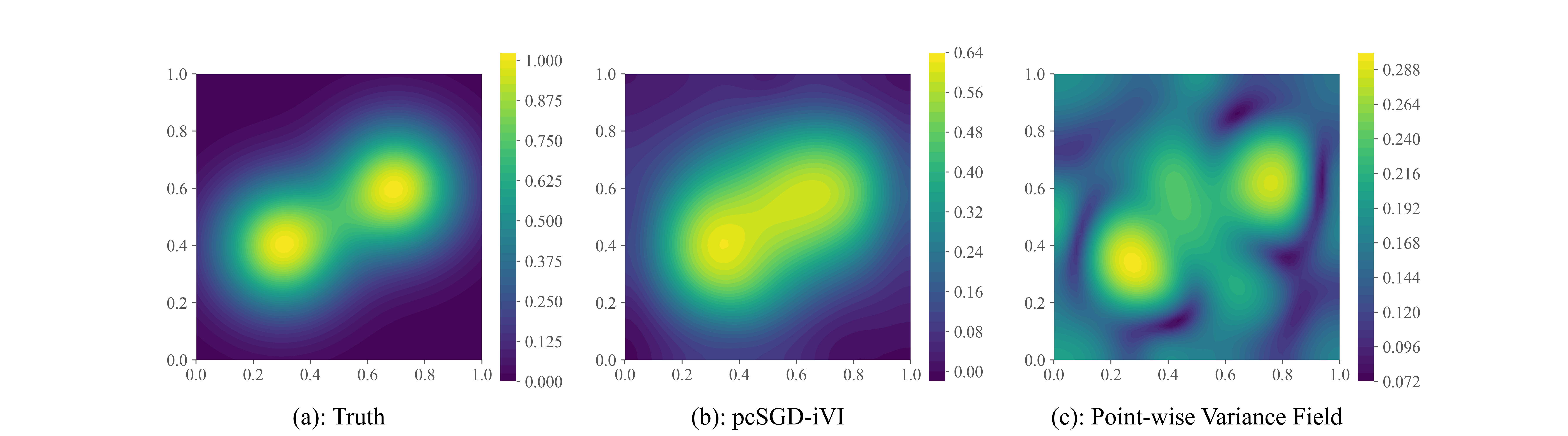}

	\caption{\emph{\small 
			(a): The background truth of $u$;
			(b): The estimated posterior mean function of $u$ obtained by pcSGD-iVI;
			(c): The estimated variance function of posterior measure of $u$ obtained by pcSGD-iVI.
	}}
	
	\label{fig:Comparison_pcSGD_Darcy}
\end{figure}

\begin{figure}
	\centering
	\includegraphics[ keepaspectratio=true, width=0.8\textwidth,  clip=true,  trim=0pt 0pt 0pt 0pt]{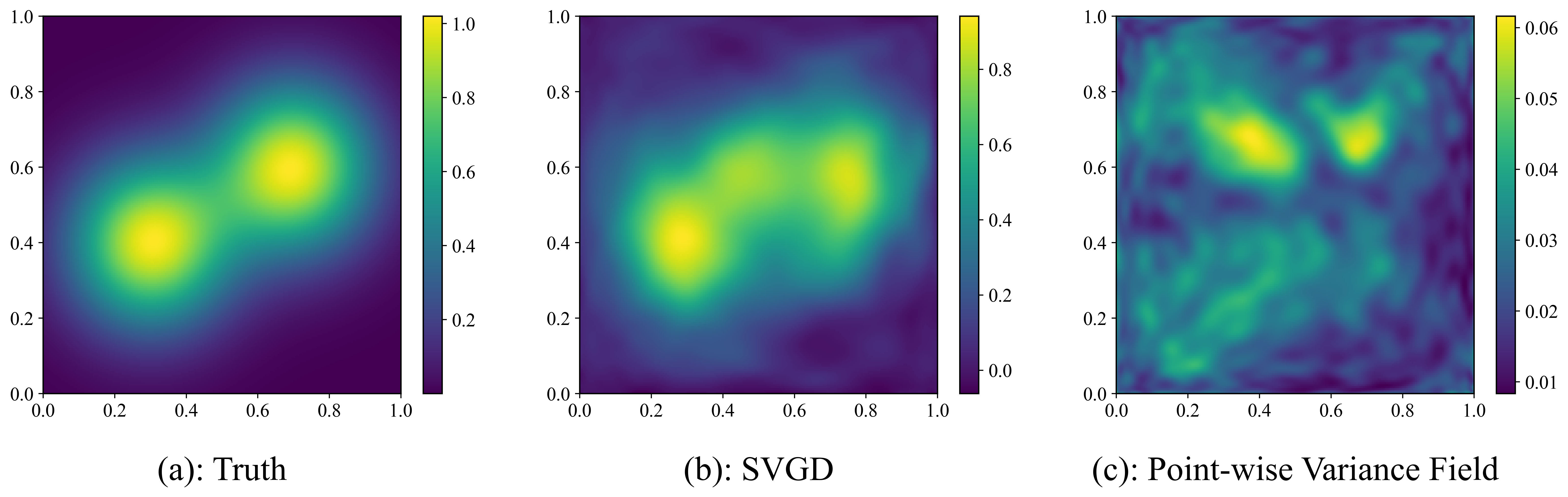}

	\caption{\emph{\small 
			(a): The background truth of $u$;
			(b): The estimated posterior mean function of $u$ obtained by SVGD;
			(c): The estimated variance function of posterior measure of $u$ obtained by SVGD.
	}}
	
	\label{fig:Comparison_SVGD_Darcy}
\end{figure}

We now denote $\delta w := F^{\prime}(\delta u)$, $w_0 = \mathcal{F}(u^*)$, and $\delta u := u - u^{*}$.
Then we obtain the linearized equation as follows:
\begin{align}
	\begin{split}
		-\nabla \cdot (e^{u^{*}}\nabla w_0 \cdot \delta u) &= \nabla \cdot (e^{u^{*}} \nabla \delta w) \quad x \in \Omega, \\
		\delta u &= 0 \quad x \in \partial \Omega,
	\end{split}
\end{align}
and the abstract form turns to 
\begin{align}
	\begin{split}
		\bm{d} &= \mathcal{F}(u^{*}) + \delta w + \bm{\epsilon},
	\end{split}
\end{align}
where $\bm{\epsilon} \sim \mathcal{N}(0, \bm{\Gamma}_{\text{noise}})$ is the random Gaussian noise.
For clarity, we list the specific choices for some parameters introduced in this subsection as follows: 
\begin{itemize}
	\item \textcolor{black}{Assume that 5$\%$ random Gaussian noise $\epsilon \sim \mathcal{N}(0, \bm{\Gamma}_{\text{noise}})$ is added, where $\bm{\Gamma}_{\text{noise}} = \tau^{-1}\textbf{I}$, and $\tau^{-1} = (0.05\max(\lvert Hu\rvert))^2$.}
	\item We assume that the data produced from the underlying log-permeability
	\begin{align*}
		u^{\dagger} = \ &\exp(-20(x_1-0.3)^2 - 20(x_2-0.4)^2) \\
		&+ \exp(-20(x_1-0.7)^2 - 20(x_2-0.6)^2).
	\end{align*}
	\item {\color{black} $\mathcal{F}(u)$ is linearized around $u^{*} = 0$}.
	\item Let domain $\Omega$ be a bounded area $(0, 1)^2$. The available data is discretized by the finite element method on a regular mesh with the number of grid points equal to $20 \times 20$.
	\item The operator $\mathcal{C}_0$ is given by $\mathcal{C}_0 = (\text{I} - \alpha\Delta)^{-2}$, where $\alpha = 0.05$ is a fixed constant. Here the Laplace operator is defined on $\Omega$ with zero Neumann boundary condition. 
	\item In order to avoid the inverse crime, a fine mesh with the number of grid points equal to $500 \times 500$ is employed for generating the data. For the inversion, a mesh with a number of grid points equal to {\color{black}$50 \times 50$} is employed.
\end{itemize}
{ To illustrate the sampling effectiveness of cSGD-iVI and pcSGD-iVI, we compare them to SVGD
method, which is proposed in \cite{jia2021stein}. }

\subsubsection{Numerical results}
In this subsection, we provide some numerical results about both two methods.
Next we show out numerical results obtained by these two methods {  and SVGD method}.

\begin{figure}
	\centering
	\includegraphics[ keepaspectratio=true, width=0.6\textwidth,  clip=true, trim=80pt 2pt 80pt 30pt]{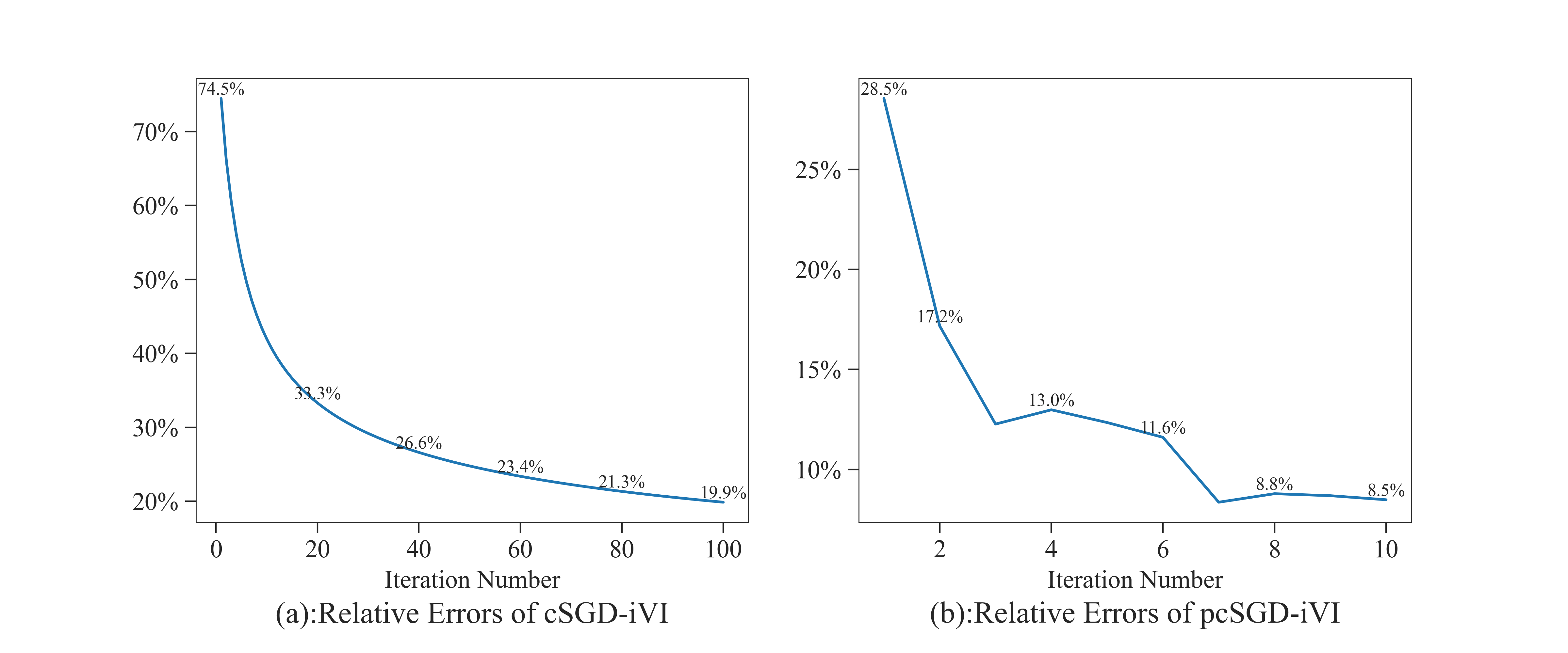}

	\caption{\emph{\small 
			(a): Relative errors of estimated posterior means of $u$ in the $L^2-norm$ of cSGD-iVI method with size of mesh $50 \times 50$; 
			(b): Relative errors of estimated posterior means of $u$ in the $L^2-norm$ of pcSGD-iVI method with size of mesh $50 \times 50$.
	}}
	
	\label{fig:Errors_Darcy}
\end{figure}

\begin{figure}[!t]
	\centering
	\subfloat[{Relative errors of SVGD I}]{
		\includegraphics[ keepaspectratio=true, width=0.32\textwidth, trim=0 0 0 0,  clip=true]{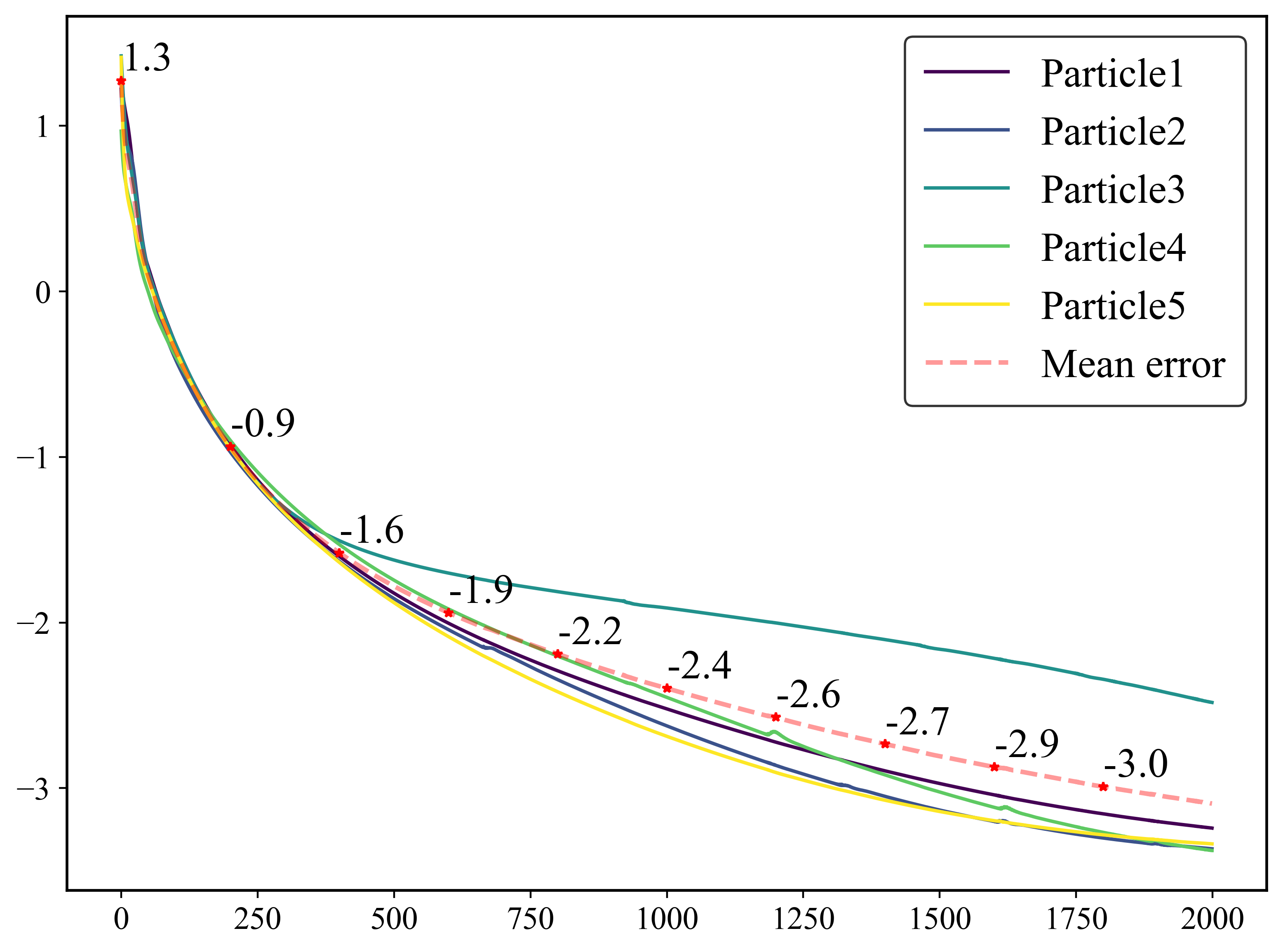}} 
	\hspace*{0.1em}
	\subfloat[{Relative errors of SVGD II}]{
		\includegraphics[ keepaspectratio=true, width=0.32\textwidth, trim=0 0 0 0,  clip=true]{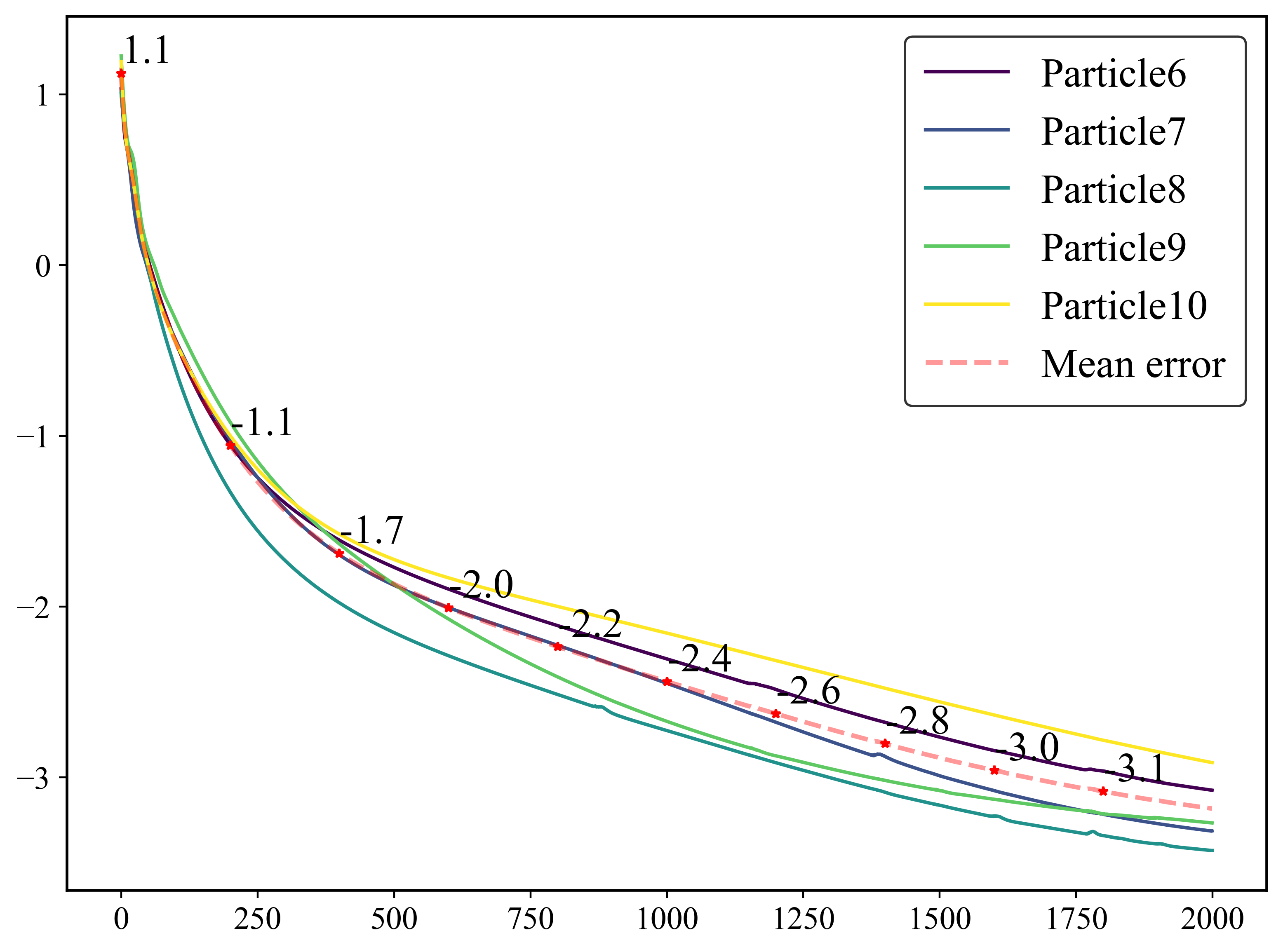}} \\
	\subfloat[{Relative errors of SVGD III}]{
		\includegraphics[ keepaspectratio=true, width=0.32\textwidth, trim=0 0 0 0,  clip=true]{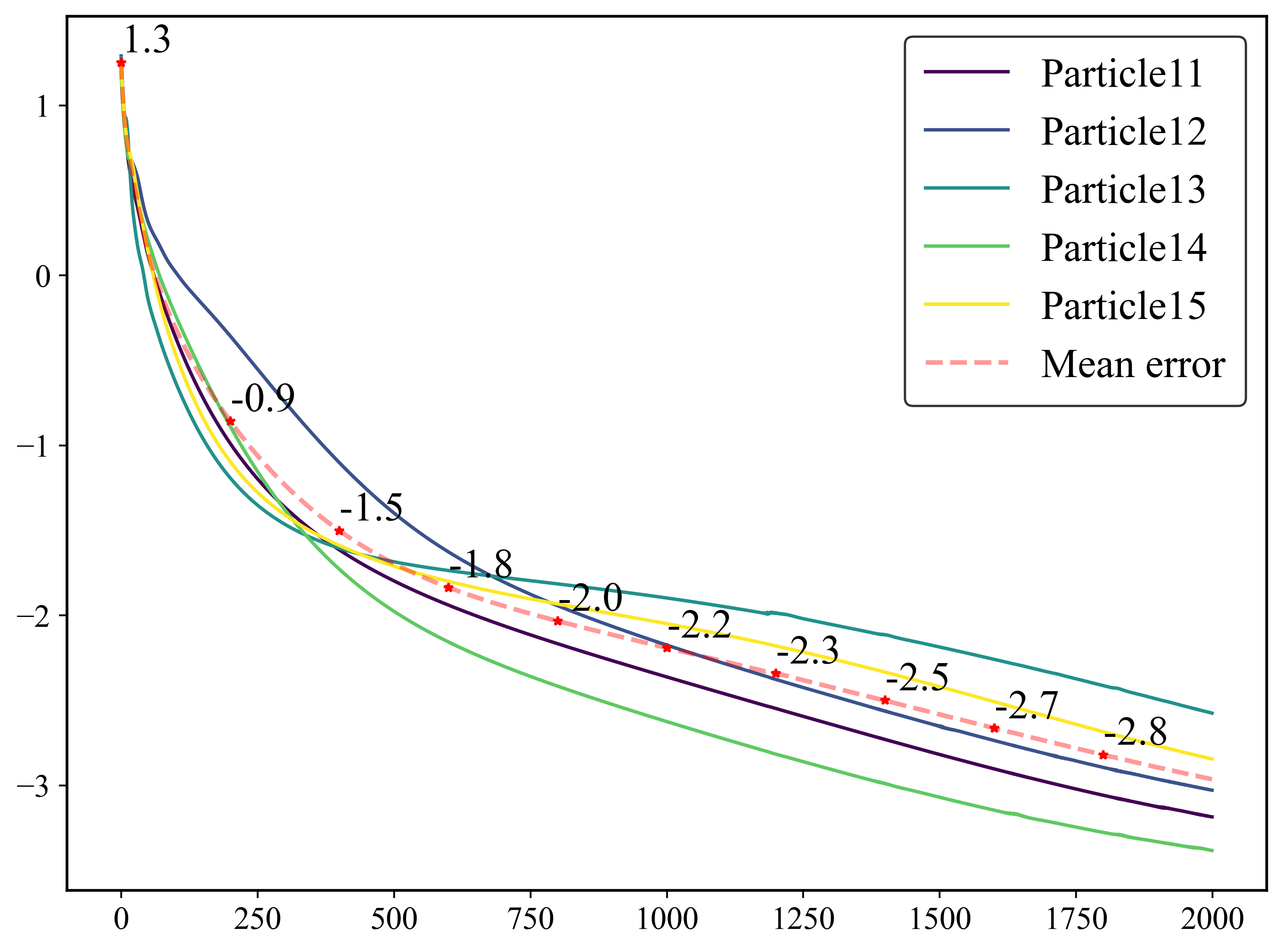}} 
	\hspace*{0.1em}
	\subfloat[{Relative errors of SVGD IV}]{
		\includegraphics[ keepaspectratio=true, width=0.32\textwidth, trim=0 0 0 0,  clip=true]{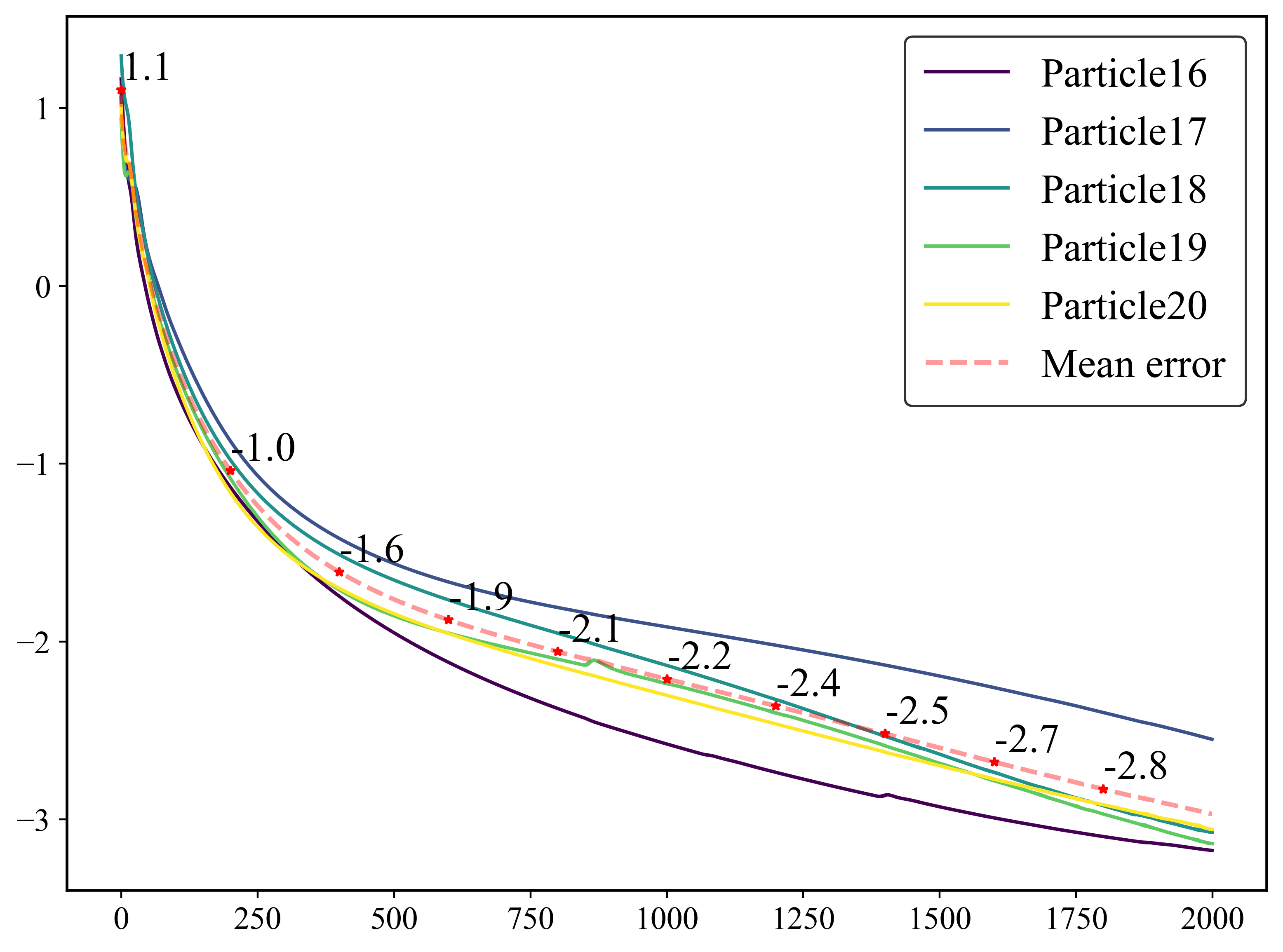}}

	\caption{\emph{\small 
			(a):the relative error curves of partilce number \(1-5\) (Logarithm);
			(b):the relative error curves of partilce number \(6-10\) (Logarithm);
			(c):the relative error curves of partilce number \(11-15\) (Logarithm);
			(d):the relative error curves of partilce number \(16-20\) (Logarithm).
	}}
	
	\label{fig:Errors_Darcy_SVGD}
\end{figure}

Firstly, in the sub-figures (a) and (b) of Figure $\ref{fig:Comparison_cSGD_Darcy}$, we show the background truth and estimated posterior mean function of $u$ obtained by cSGD-iVI method.
The estimated posterior mean function of $u$ is visually different from the true function.
In sub-figure (a) of Figure $\ref{fig:Errors_Darcy}$, we draw the relative error curve calculated in $L^2$-norm, defined in $(3.4)$.
The relative error curve illustrates that the convergence speed is fast since the descending trend is rapid at first $80$ steps, and the relative error is stable around $20 \%$ at the end of the iteration steps.
This provides quantitative evidence that the estimated posterior mean function is different from the background truth of $u$.
As a result, combining the visual (sub-figures (a), (b) of Figure $\ref{fig:Comparison_cSGD_Darcy}$) and quantitative (relative errors shown in sub-figure (a) of Figure $\ref{fig:Errors_Darcy}$) evidence, we conclude that the cSGD-iVI method provides an estimated posterior mean function of the parameter $u$ which is far different from the background truth.
On the other hand, in the sub-figures (a) and (b) of Figure $\ref{fig:Comparison_pcSGD_Darcy}$, we show the background truth and estimated posterior mean function of $u$ obtained by pcSGD-iVI mthod.
The estimated posterior mean function of $u$ is visually similar the truth.
In sub-figure (b) of Figure $\ref{fig:Errors_Darcy}$, we draw the relative error curve calculated in $L^2$-norm.
The relative error curve illustrates that the convergence speed is fast since the descending trend is rapid at first $6$ steps, and the relative error is stable around $8 \%$ at the end of the iteration steps.
This provides quantitative evidence that the estimated posterior mean function is similar to the background truth of $u$.
As a result, combining the visual (sub-figures (a), (b) of Figure $\ref{fig:Comparison_pcSGD_Darcy}$) and quantitative (relative errors shown in sub-figure (b) of Figure $\ref{fig:Errors_Darcy}$) evidence, we conclude that the pcSGD-iVI method provides an estimated posterior mean function of the parameter $u$ which is similar to the background truth.
{ Then we show the numerical results obtained by SVGD method.
We adopt the set of particles as \(\{u_i\}_{i=1}^{20}.\)
In the sub-figures (a) and (b) of Figure \ref{fig:Comparison_SVGD_Darcy}, we show the background truth and estimated posterior mean function of $u$ obtained by SVGD method.
It can be seen that the highlighted area of the left peak of the approximate posterior mean function of u is quite similar to the left peak of the true function in sub-figure (a). 
The right peak is significantly different in both highlight and shape, and there is also an additional highlighted area in the middle region. 
In all sub-figures of Figure \ref{fig:Errors_Darcy_SVGD}, we draw the relative logarithm error curves of particles calculated in $L^2$-norm.
The relative error curve illustrates that the decreasing trend of each particle is consistent with the trend of the average error.
The convergence speed is fast since the descending trend is rapid at first $250$ steps, and the relative error is stable around \(-3.0\)
($5 \%$) at the end of the iteration steps.
It takes \(2000\) steps during the convergence process, and the running time is \(11437.9929s\)(\(5.7190s\) per-step).
As a result, combining the visual (sub-figures (a), (b) of Figure \ref{fig:Comparison_SVGD_Darcy}) and quantitative (relative errors shown in Figure \ref{fig:Errors_Darcy_SVGD}) evidence, we conclude that the estimated posterior mean obtained by SVGD method differs in shape, but the relative error is relatively small.}

\begin{figure}
	\centering
	\includegraphics[ keepaspectratio=true, width=0.8\textwidth,  clip=true, trim=120pt 2pt 120pt 30pt]{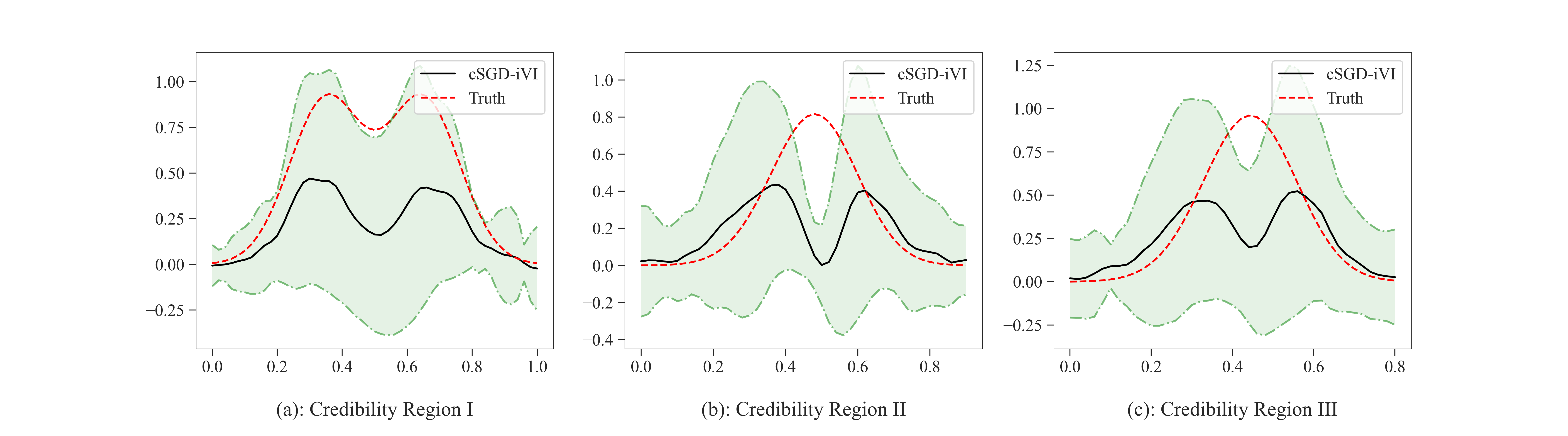}

	\caption{\emph{\small 
			The $95 \%$ credibility region of the estimated posterior mean function represented by the green shade area.
			(a): The black solid line represents the estimated mean $\lbrace u(x_i, x_i) \rbrace^{50}_{i=1}$, and red dashed line represents the true function $\lbrace u^{\dagger}(x_i, x_i) \rbrace^{50}_{i=1}$.
			(b): The black solid line represents the estimated mean $\lbrace u(x_i, x_{i+5}) \rbrace^{45}_{i=1}$, and red dashed line represents the true function $\lbrace u^{\dagger}(x_i, x_{i+5}) \rbrace^{45}_{i=1}$.
			(c): The black solid line represents the estimated mean $\lbrace u(x_i, x_{i+10}) \rbrace^{40}_{i=1}$, and red dashed line represents the true function $\lbrace u^{\dagger}(x_i, x_{i+10}) \rbrace^{40}_{i=1}$. 
	}}
	
	\label{fig:Variances_cSGD_Darcy}
\end{figure}

\begin{figure}
	\centering
	\includegraphics[ keepaspectratio=true, width=0.8\textwidth,  clip=true,trim=120pt 2pt 120pt 30pt]{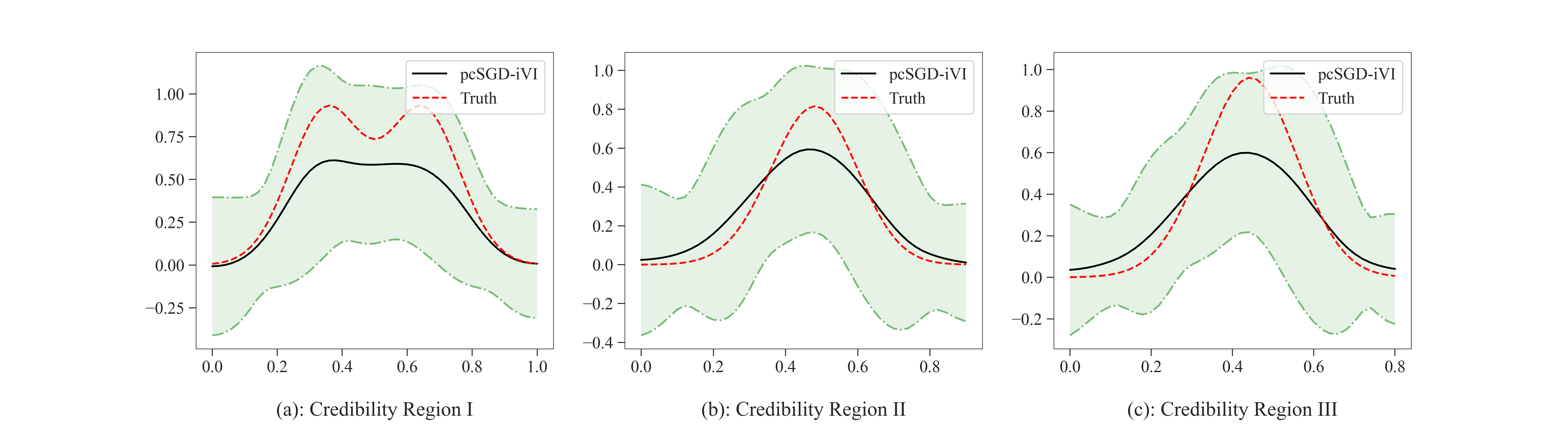}

	\caption{\emph{\small 
			The $95 \%$ credibility region of the estimated posterior mean function represented by the green shade area.
			(a): The black solid line represents the estimated mean $\lbrace u(x_i, x_i) \rbrace^{50}_{i=1}$, and red dashed line represents the true function $\lbrace u^{\dagger}(x_i, x_i) \rbrace^{50}_{i=1}$.
			(b): The black solid line represents the estimated mean $\lbrace u(x_i, x_{i+5}) \rbrace^{45}_{i=1}$, and red dashed line represents the true function $\lbrace u^{\dagger}(x_i, x_{i+5}) \rbrace^{45}_{i=1}$.
			(c): The black solid line represents the estimated mean $\lbrace u(x_i, x_{i+10}) \rbrace^{40}_{i=1}$, and red dashed line represents the true function $\lbrace u^{\dagger}(x_i, x_{i+10}) \rbrace^{40}_{i=1}$. 
	}}
	
	\label{fig:Variances_pcSGD_Darcy}
\end{figure}

\begin{figure}
	\centering
	\includegraphics[ keepaspectratio=true, width=0.8\textwidth,  clip=true,trim=120pt 30pt 120pt 0pt]{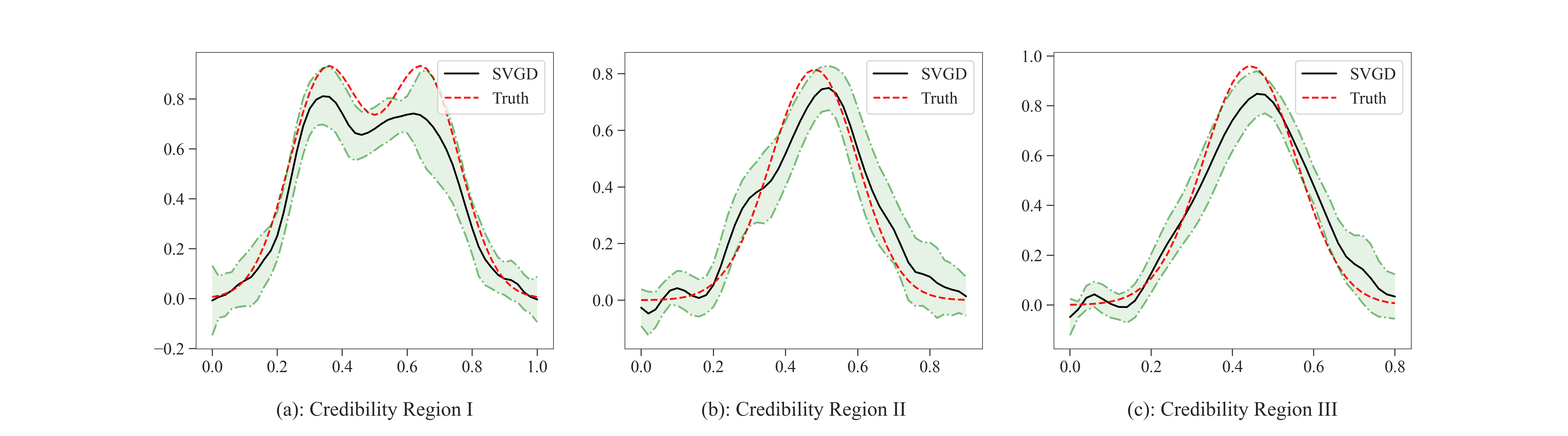}

	\caption{\emph{\small 
			The $95 \%$ credibility region of the estimated posterior mean function represented by the green shade area.
			(a): The black solid line represents the estimated mean $\lbrace u(x_i, x_i) \rbrace^{50}_{i=1}$, and red dashed line represents the true function $\lbrace u^{\dagger}(x_i, x_i) \rbrace^{50}_{i=1}$.
			(b): The black solid line represents the estimated mean $\lbrace u(x_i, x_{i+5}) \rbrace^{45}_{i=1}$, and red dashed line represents the true function $\lbrace u^{\dagger}(x_i, x_{i+5}) \rbrace^{45}_{i=1}$.
			(c): The black solid line represents the estimated mean $\lbrace u(x_i, x_{i+10}) \rbrace^{40}_{i=1}$, and red dashed line represents the true function $\lbrace u^{\dagger}(x_i, x_{i+10}) \rbrace^{40}_{i=1}$. 
	}}
	
	\label{fig:Variances_SVGD_Darcy}
\end{figure}

Secondly, we provide some discussions of the estimated posterior covariance functions about the parameter $u$.
In sub-figure (c) of Figure $\ref{fig:Comparison_cSGD_Darcy}$, we draw the point-wise variance field of the posterior measure $u$ obtained by cSGD-iVI.
To provide detailed evidence, we draw the comparisons of the estimated mean function and back-ground truth calculated by different mesh points in Figure $\ref{fig:Variances_cSGD_Darcy}$, which are given by $\lbrace (x_i, x_i) \rbrace^{n}_{i=1}$, $\lbrace (x_i, x_{i+5}) \rbrace^{n-5}_{i=1}$, and $\lbrace (x_i, x_{i+10}) \rbrace^{n-10}_{i=1}$($n = 50$ according to the mesh size $= 50 \times 50$).
And the green shade area represents the $95 \%$ credibility region according to estimated posterior mean function.
In the sub-figures (a), (b), and (c) of Figure $\ref{fig:Variances_cSGD_Darcy}$, the black solid lines represent the estimated mean $\lbrace u(x_i, x_i) \rbrace^{50}_{i=1}$, $\lbrace u(x_i, x_{i+5}) \rbrace^{45}_{i=1}$, and $\lbrace u(x_i, x_{i+10}) \rbrace^{40}_{i=1}$ obtained by cSGD-iVI method, respectively.
And red solid lines represent the true function $\lbrace u^{\dagger}(x_i, x_i) \rbrace^{50}_{i=1}$, $\lbrace u^{\dagger}(x_i, x_{i+5}) \rbrace^{45}_{i=1}$, and $\lbrace u^{\dagger}(x_i, x_{i+10}) \rbrace^{40}_{i=1}$, respectively.
We see that the credibility region doesn't contain background truth completely, particularly in sub-figures (b) and (c).
This indicates that cSGD-iVI can not reflect the uncertainties of parameter $u$.
On the other hand, in sub-figure (c) of Figure $\ref{fig:Comparison_pcSGD_Darcy}$, we draw the point-wise variance field obtained by pcSGD-iVI.
Moreover, we draw the comparisons of the estimated mean function and back-ground truth calculated by different mesh points in Figure $\ref{fig:Variances_pcSGD_Darcy}$, and the green shade area represents the $95 \%$ credibility region according to estimated posterior mean function.
We see that the credibility region contains background truth completely.
This indicates that with the help of preconditioning operator $T$, the estimated posterior measure obtained by pcSGD-iVI reflects the uncertainties of $u$ accurately.
The Bayesian setup derived by pcSGD-iVI is meaningful and in accordance with the frequentist theoretical investigations of the posterior consistency \cite{wang2019frequentist, zhang2020convergence}. 
Combining with the variance function, we conclude that the proposed pcSGD-iVI method quantifies the uncertainties of the parameter $u$.
{ 
As for SVGD method, in sub-figure (c) of Figure \ref{fig:Comparison_SVGD_Darcy}, we draw the point-wise variance field of the posterior measure $u$ obtained by SVGD.
To provide detailed evidence, we draw the comparisons of the estimated mean function and back-ground truth calculated by different mesh points in Figure \ref{fig:Variances_SVGD_Darcy}, which are given by $\lbrace (x_i, x_i) \rbrace^{n}_{i=1}$, $\lbrace (x_i, x_{i+5}) \rbrace^{n-5}_{i=1}$, and $\lbrace (x_i, x_{i+10}) \rbrace^{n-10}_{i=1}$($n = 50$ according to the mesh size $= 50 \times 50$).
And the green shade area represents the $95 \%$ credibility region according to estimated posterior mean function.
In the sub-figures (a), (b), and (c) of Figure \ref{fig:Variances_SVGD_Darcy}, the black solid lines represent the estimated mean $\lbrace u(x_i, x_i) \rbrace^{50}_{i=1}$, $\lbrace u(x_i, x_{i+5}) \rbrace^{45}_{i=1}$, and $\lbrace u(x_i, x_{i+10}) \rbrace^{40}_{i=1}$ obtained by SVGD method, respectively.
And red solid lines represent the true function $\lbrace u^{\dagger}(x_i, x_i) \rbrace^{50}_{i=1}$, $\lbrace u^{\dagger}(x_i, x_{i+5}) \rbrace^{45}_{i=1}$, and $\lbrace u^{\dagger}(x_i, x_{i+10}) \rbrace^{40}_{i=1}$, respectively.
We see that the credibility region doesn't contain background truth completely in all sub-figures.
This indicates that SVGD cannot completely reflect the uncertainties of parameter $u$.
}

{ And we derive the following results:
For posterior mean accuracy: SVGD attains the smallest relative $L^2$ error among these three methods, based on the visual (Figures \ref{fig:Comparison_cSGD_Darcy}, \ref{fig:Comparison_pcSGD_Darcy} and \ref{fig:Comparison_SVGD_Darcy}) and quantitative (Figures \ref{fig:Errors_Darcy} and \ref{fig:Errors_Darcy_SVGD}).
For uncertainty quantification: pcSGD-iVI reflects the uncertainties of $u$ accurately, however SVGD and cSGD-iVI can not reflect the uncertainties of parameter $u$, based on the visual (sub-figure (c) of Figures \ref{fig:Comparison_cSGD_Darcy}, \ref{fig:Comparison_pcSGD_Darcy} and \ref{fig:Comparison_SVGD_Darcy}) and quantitative (Figures \ref{fig:Variances_cSGD_Darcy}, \ref{fig:Variances_pcSGD_Darcy} and \ref{fig:Variances_SVGD_Darcy}).
For computational cost: under the same settings mentioned in B.2.1 of supplementary material, cSGD-iVI has the lowest cost(\(0.38312s\) per-step), followed by pcSGD-iVI(\(1.89033s\) per-step), while SVGD is the most expensive(\(5.7190s\) per-step); the computational costs of cSGD-iVI and pcSGD-iVI are shown in Tables \ref{table:Darcy_cost_cSGD} and \ref{table:Darcy_cost_pcSGD}.
}

{ 
\subsubsection{Performance and cost on different mesh sizes}
In order to provide more evidence, we employ the mesh size \(= \{40 \times 40, 60 \times 60\}\) for comparsions.
And we provide the comparsions of computational cost on different meshes \(= \{40 \times 40, 50 \times 50, 60 \times 60\}\).

\subsubsection*{Mesh size \(= 40 \times 40\)}

Firstly, in the sub-figures (a) and (b) of Figure $\ref{fig:Comparison_cSGD_Darcy40}$, we show the background truth and estimated posterior mean function of $u$ obtained by cSGD-iVI method.
The estimated posterior mean function of $u$ is visually different from the true function.
In sub-figure (a) of Figure $\ref{fig:Errors_Darcy40}$, we draw the relative error curve calculated in $L^2$-norm, defined in $(3.4)$.
The relative error curve illustrates that the convergence speed is fast since the descending trend is rapid at first $80$ steps, and the relative error is stable around $20 \%$ at the end of the iteration steps.
This provides quantitative evidence that the estimated posterior mean function is different from the background truth of $u$.
As a result, combining the visual (sub-figures (a), (b) of Figure $\ref{fig:Comparison_cSGD_Darcy40}$) and quantitative (relative errors shown in sub-figure (a) of Figure $\ref{fig:Errors_Darcy40}$) evidence, we conclude that the cSGD-iVI method provides an estimated posterior mean function of the parameter $u$ which is far different from the background truth.
On the other hand, in the sub-figures (a) and (b) of Figure $\ref{fig:Comparison_pcSGD_Darcy40}$, we show the background truth and estimated posterior mean function of $u$ obtained by pcSGD-iVI mthod.
The estimated posterior mean function of $u$ is visually similar the truth.
In sub-figure (b) of Figure $\ref{fig:Errors_Darcy40}$, we draw the relative error curve calculated in $L^2$-norm.
The relative error curve illustrates that the convergence speed is fast since the descending trend is rapid at first $6$ steps, and the relative error is stable around $8 \%$ at the end of the iteration steps.
This provides quantitative evidence that the estimated posterior mean function is similar to the background truth of $u$.
As a result, combining the visual (sub-figures (a), (b) of Figure $\ref{fig:Comparison_pcSGD_Darcy40}$) and quantitative (relative errors shown in sub-figure (b) of Figure $\ref{fig:Errors_Darcy40}$) evidence, we conclude that the pcSGD-iVI method provides an estimated posterior mean function of the parameter $u$ which is similar to the background truth.

Secondly, we provide some discussions of the estimated posterior covariance functions about the parameter $u$.
In sub-figure (c) of Figure $\ref{fig:Comparison_cSGD_Darcy40}$, we draw the point-wise variance field of the posterior measure $u$ obtained by cSGD-iVI.
To provide detailed evidence, we draw the comparisons of the estimated mean function and back-ground truth calculated by different mesh points in Figure $\ref{fig:Variances_cSGD_Darcy40}$, which are given by $\lbrace (x_i, x_i) \rbrace^{n}_{i=1}$, $\lbrace (x_i, x_{i+5}) \rbrace^{n-5}_{i=1}$, and $\lbrace (x_i, x_{i+10}) \rbrace^{n-10}_{i=1}$($n = 40$ according to the mesh size $= 40 \times 40$).
And the green shade area represents the $95 \%$ credibility region according to estimated posterior mean function.
In the sub-figures (a), (b), and (c) of Figure $\ref{fig:Variances_cSGD_Darcy40}$, the black solid lines represent the estimated mean $\lbrace u(x_i, x_i) \rbrace^{40}_{i=1}$, $\lbrace u(x_i, x_{i+5}) \rbrace^{35}_{i=1}$, and $\lbrace u(x_i, x_{i+10}) \rbrace^{30}_{i=1}$ obtained by cSGD-iVI method, respectively.
And red solid lines represent the true function $\lbrace u^{\dagger}(x_i, x_i) \rbrace^{40}_{i=1}$, $\lbrace u^{\dagger}(x_i, x_{i+5}) \rbrace^{35}_{i=1}$, and $\lbrace u^{\dagger}(x_i, x_{i+10}) \rbrace^{30}_{i=1}$, respectively.
We see that the credibility region doesn't contain background truth completely, particularly in sub-figures (b) and (c).
This indicates that cSGD-iVI can not reflect the uncertainties of parameter $u$.
On the other hand, in sub-figure (c) of Figure $\ref{fig:Comparison_pcSGD_Darcy40}$, we draw the point-wise variance field obtained by pcSGD-iVI.
Moreover, we draw the comparisons of the estimated mean function and back-ground truth calculated by different mesh points in Figure $\ref{fig:Variances_pcSGD_Darcy40}$, and the green shade area represents the $95 \%$ credibility region according to estimated posterior mean function.
We see that the credibility region contains background truth completely.
This indicates that with the help of preconditioning operator $T$, the estimated posterior measure obtained by pcSGD-iVI reflects the uncertainties of $u$ accurately.
The Bayesian setup derived by pcSGD-iVI is meaningful and in accordance with the frequentist theoretical investigations of the posterior consistency \cite{wang2019frequentist, zhang2020convergence}. 
Combining with the variance function, we conclude that the proposed pcSGD-iVI method quantifies the uncertainties of the parameter $u$.
In summary, we see that pcSGD-iVI provides a good estimate of mean function and covariance operator, based on visual (Figures $\ref{fig:Comparison_pcSGD_Darcy40}$ and $\ref{fig:Variances_pcSGD_Darcy40}$) and quantitative (sub-figure (b) of Figure $\ref{fig:Errors_Darcy40}$) evidence; while the estimation of cSGD-iVI is much worse, based on visual (Figures $\ref{fig:Comparison_cSGD_Darcy40}$ and $\ref{fig:Variances_cSGD_Darcy40}$) and quantitative (sub-figure (a) of Figure $\ref{fig:Errors_Darcy}$) evidence.

\begin{figure}
	\centering
	\includegraphics[ keepaspectratio=true, width=0.8\textwidth,  clip=true,  trim=120pt 2pt 120pt 30pt]{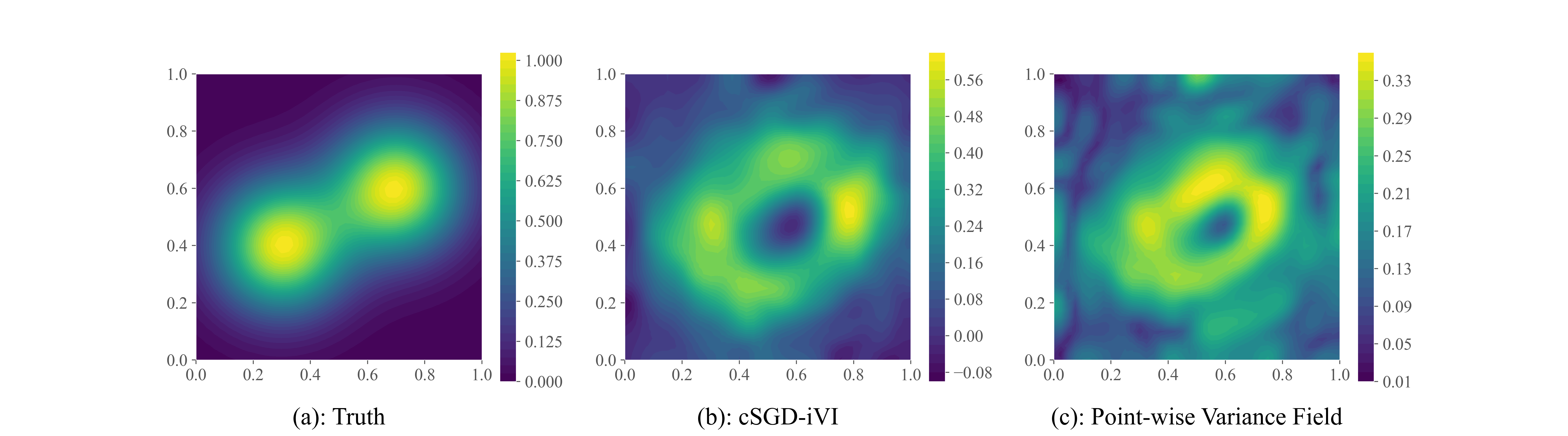}

	\caption{\emph{\small 
			(a): The background truth of $u$;
			(b): The estimated posterior mean function of $u$ obtained by cSGD-iVI;
			(c): The estimated variance function of posterior measure of $u$ obtained by cSGD-iVI.
	}}
	
	\label{fig:Comparison_cSGD_Darcy40}
\end{figure}

\begin{figure}
	\centering
	\includegraphics[ keepaspectratio=true, width=0.8\textwidth,  clip=true,  trim=120pt 2pt 120pt 30pt]{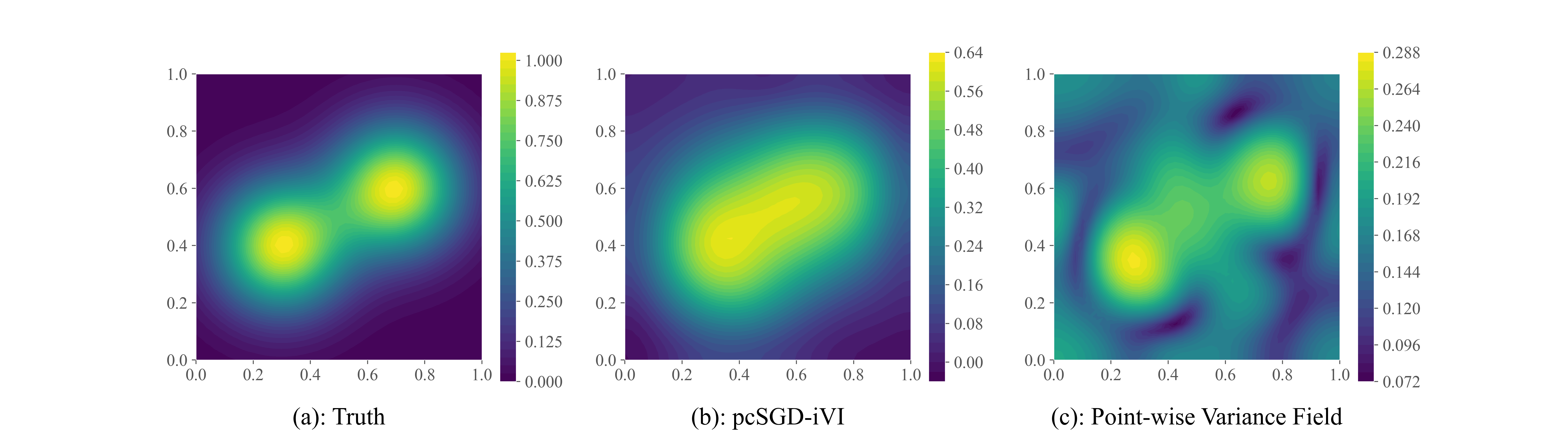}

	\caption{\emph{\small 
			(a): The background truth of $u$;
			(b): The estimated posterior mean function of $u$ obtained by pcSGD-iVI;
			(c): The estimated variance function of posterior measure of $u$ obtained by pcSGD-iVI.
	}}
	
	\label{fig:Comparison_pcSGD_Darcy40}
\end{figure}

\begin{figure}
	\centering
	\includegraphics[ keepaspectratio=true, width=0.6\textwidth,  clip=true, trim=80pt 2pt 80pt 30pt]{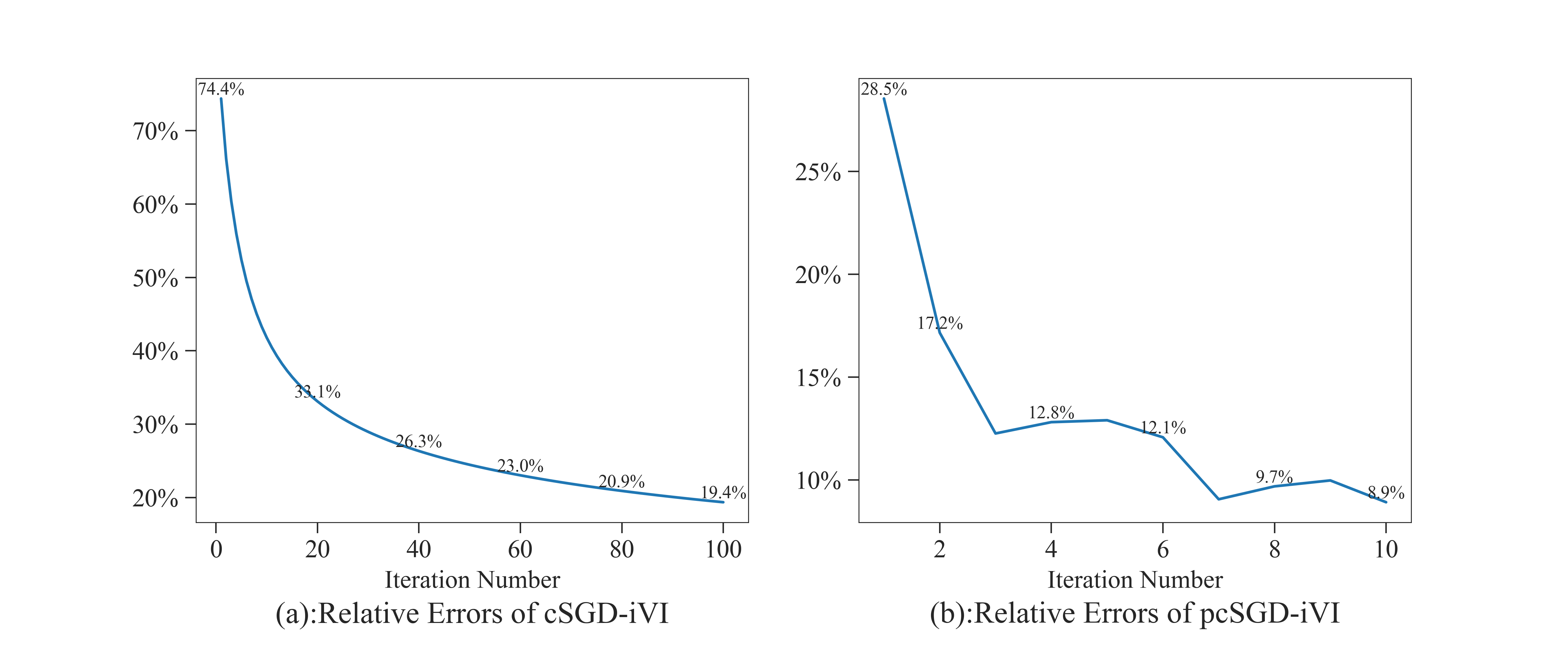}

	\caption{\emph{\small 
			(a): Relative errors of estimated posterior means of $u$ in the $L^2-norm$ of cSGD-iVI method with size of mesh $40 \times 40$; 
			(b): Relative errors of estimated posterior means of $u$ in the $L^2-norm$ of pcSGD-iVI method with size of mesh $40 \times 40$.
	}}
	
	\label{fig:Errors_Darcy40}
\end{figure}

\begin{figure}[!t]
	\centering
	\includegraphics[ keepaspectratio=true, width=0.8\textwidth,  clip=true, trim=120pt 2pt 120pt 30pt]{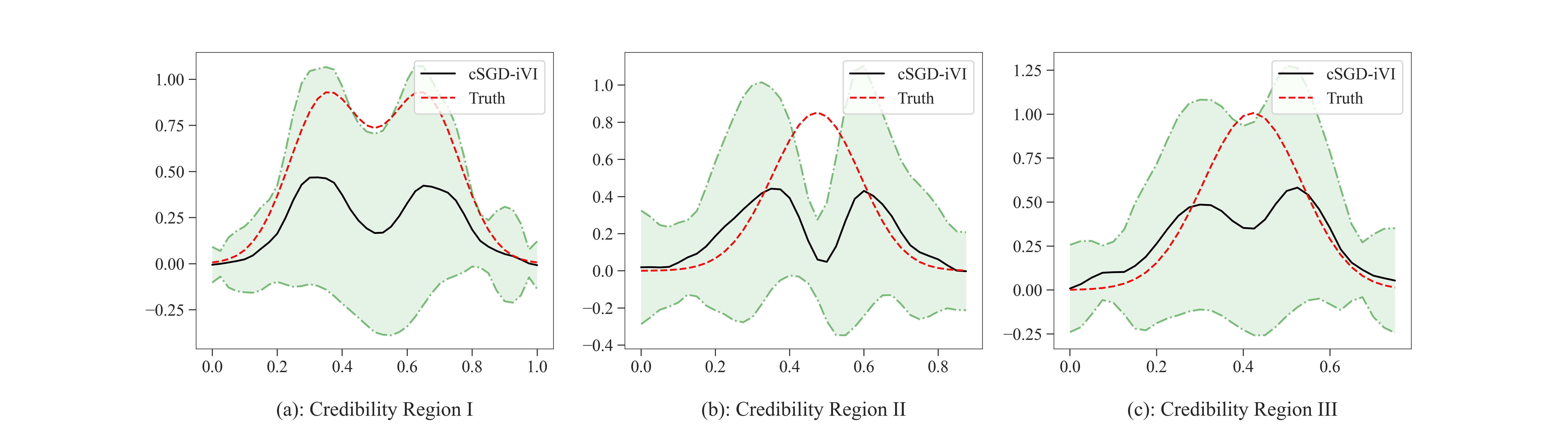}

	\caption{\emph{\small 
			The $95 \%$ credibility region of the estimated posterior mean function represented by the green shade area.
			(a): The black solid line represents the estimated mean $\lbrace u(x_i, x_i) \rbrace^{40}_{i=1}$, and red dashed line represents the true function $\lbrace u^{\dagger}(x_i, x_i) \rbrace^{40}_{i=1}$.
			(b): The black solid line represents the estimated mean $\lbrace u(x_i, x_{i+5}) \rbrace^{35}_{i=1}$, and red dashed line represents the true function $\lbrace u^{\dagger}(x_i, x_{i+5}) \rbrace^{35}_{i=1}$.
			(c): The black solid line represents the estimated mean $\lbrace u(x_i, x_{i+10}) \rbrace^{30}_{i=1}$, and red dashed line represents the true function $\lbrace u^{\dagger}(x_i, x_{i+10}) \rbrace^{30}_{i=1}$. 
	}}
	
	\label{fig:Variances_cSGD_Darcy40}
\end{figure}

\begin{figure}[!t]
	\centering
	\includegraphics[ keepaspectratio=true, width=0.8\textwidth,  clip=true,trim=120pt 2pt 120pt 30pt]{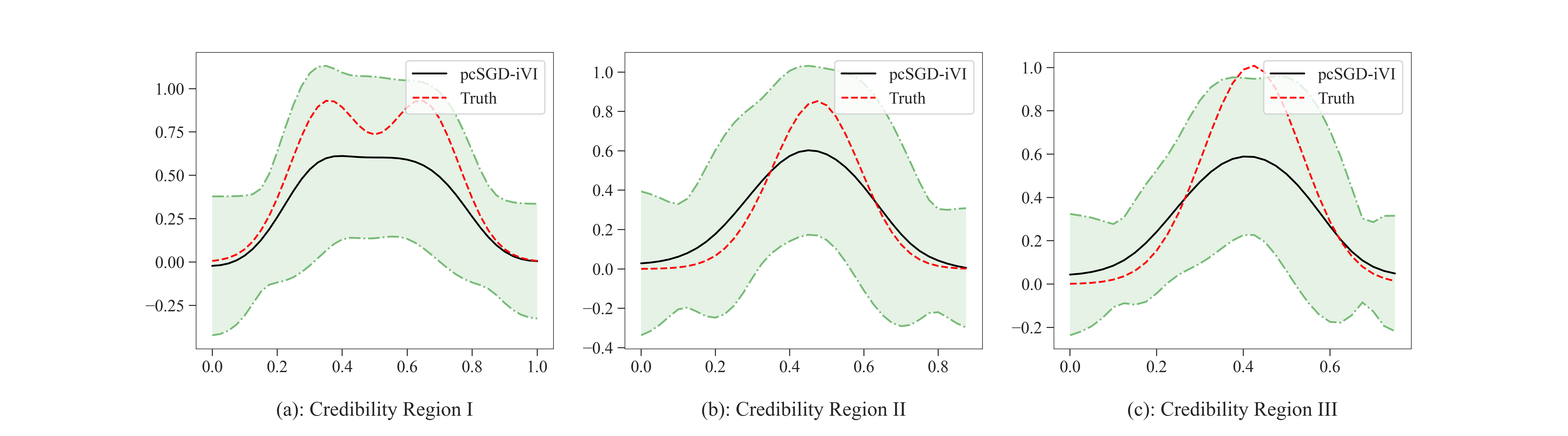}

	\caption{\emph{\small 
			The $95 \%$ credibility region of the estimated posterior mean function represented by the green shade area.
			(a): The black solid line represents the estimated mean $\lbrace u(x_i, x_i) \rbrace^{40}_{i=1}$, and red dashed line represents the true function $\lbrace u^{\dagger}(x_i, x_i) \rbrace^{40}_{i=1}$.
			(b): The black solid line represents the estimated mean $\lbrace u(x_i, x_{i+5}) \rbrace^{35}_{i=1}$, and red dashed line represents the true function $\lbrace u^{\dagger}(x_i, x_{i+5}) \rbrace^{35}_{i=1}$.
			(c): The black solid line represents the estimated mean $\lbrace u(x_i, x_{i+10}) \rbrace^{30}_{i=1}$, and red dashed line represents the true function $\lbrace u^{\dagger}(x_i, x_{i+10}) \rbrace^{30}_{i=1}$. 
	}}
	
	\label{fig:Variances_pcSGD_Darcy40}
\end{figure}

\subsubsection*{Mesh size \(= 60 \times 60\)}

Firstly, in the sub-figures (a) and (b) of Figure $\ref{fig:Comparison_cSGD_Darcy60}$, we show the background truth and estimated posterior mean function of $u$ obtained by cSGD-iVI method.
The estimated posterior mean function of $u$ is visually different from the true function.
In sub-figure (a) of Figure $\ref{fig:Errors_Darcy60}$, we draw the relative error curve calculated in $L^2$-norm, defined in $(3.4)$.
The relative error curve illustrates that the convergence speed is fast since the descending trend is rapid at first $80$ steps, and the relative error is stable around $20 \%$ at the end of the iteration steps.
This provides quantitative evidence that the estimated posterior mean function is different from the background truth of $u$.
As a result, combining the visual (sub-figures (a), (b) of Figure $\ref{fig:Comparison_cSGD_Darcy60}$) and quantitative (relative errors shown in sub-figure (a) of Figure $\ref{fig:Errors_Darcy60}$) evidence, we conclude that the cSGD-iVI method provides an estimated posterior mean function of the parameter $u$ which is far different from the background truth.
On the other hand, in the sub-figures (a) and (b) of Figure $\ref{fig:Comparison_pcSGD_Darcy60}$, we show the background truth and estimated posterior mean function of $u$ obtained by pcSGD-iVI mthod.
The estimated posterior mean function of $u$ is visually similar the truth.
In sub-figure (b) of Figure $\ref{fig:Errors_Darcy60}$, we draw the relative error curve calculated in $L^2$-norm.
The relative error curve illustrates that the convergence speed is fast since the descending trend is rapid at first $6$ steps, and the relative error is stable around $8 \%$ at the end of the iteration steps.
This provides quantitative evidence that the estimated posterior mean function is similar to the background truth of $u$.
As a result, combining the visual (sub-figures (a), (b) of Figure $\ref{fig:Comparison_pcSGD_Darcy60}$) and quantitative (relative errors shown in sub-figure (b) of Figure $\ref{fig:Errors_Darcy60}$) evidence, we conclude that the pcSGD-iVI method provides an estimated posterior mean function of the parameter $u$ which is similar to the background truth.

Secondly, we provide some discussions of the estimated posterior covariance functions about the parameter $u$.
In sub-figure (c) of Figure $\ref{fig:Comparison_cSGD_Darcy60}$, we draw the point-wise variance field of the posterior measure $u$ obtained by cSGD-iVI.
To provide detailed evidence, we draw the comparisons of the estimated mean function and back-ground truth calculated by different mesh points in Figure $\ref{fig:Variances_cSGD_Darcy60}$, which are given by $\lbrace (x_i, x_i) \rbrace^{n}_{i=1}$, $\lbrace (x_i, x_{i+5}) \rbrace^{n-5}_{i=1}$, and $\lbrace (x_i, x_{i+10}) \rbrace^{n-10}_{i=1}$($n = 60$ according to the mesh size $= 60 \times 60$).
And the green shade area represents the $95 \%$ credibility region according to estimated posterior mean function.
In the sub-figures (a), (b), and (c) of Figure $\ref{fig:Variances_cSGD_Darcy60}$, the black solid lines represent the estimated mean $\lbrace u(x_i, x_i) \rbrace^{60}_{i=1}$, $\lbrace u(x_i, x_{i+5}) \rbrace^{55}_{i=1}$, and $\lbrace u(x_i, x_{i+10}) \rbrace^{50}_{i=1}$ obtained by cSGD-iVI method, respectively.
And red solid lines represent the true function $\lbrace u^{\dagger}(x_i, x_i) \rbrace^{60}_{i=1}$, $\lbrace u^{\dagger}(x_i, x_{i+5}) \rbrace^{55}_{i=1}$, and $\lbrace u^{\dagger}(x_i, x_{i+10}) \rbrace^{50}_{i=1}$, respectively.
We see that the credibility region doesn't contain background truth completely, particularly in sub-figures (b) and (c).
This indicates that cSGD-iVI can not reflect the uncertainties of parameter $u$.
On the other hand, in sub-figure (c) of Figure $\ref{fig:Comparison_pcSGD_Darcy60}$, we draw the point-wise variance field obtained by pcSGD-iVI.
Moreover, we draw the comparisons of the estimated mean function and back-ground truth calculated by different mesh points in Figure $\ref{fig:Variances_pcSGD_Darcy60}$, and the green shade area represents the $95 \%$ credibility region according to estimated posterior mean function.
We see that the credibility region contains background truth completely.
This indicates that with the help of preconditioning operator $T$, the estimated posterior measure obtained by pcSGD-iVI reflects the uncertainties of $u$ accurately.
The Bayesian setup derived by pcSGD-iVI is meaningful and in accordance with the frequentist theoretical investigations of the posterior consistency \cite{wang2019frequentist, zhang2020convergence}. 
Combining with the variance function, we conclude that the proposed pcSGD-iVI method quantifies the uncertainties of the parameter $u$.
In summary, we see that pcSGD-iVI provides a good estimate of mean function and covariance operator, based on visual (Figures $\ref{fig:Comparison_pcSGD_Darcy60}$ and $\ref{fig:Variances_pcSGD_Darcy60}$) and quantitative (sub-figure (b) of Figure $\ref{fig:Errors_Darcy60}$) evidence; while the estimation of cSGD-iVI is much worse, based on visual (Figures $\ref{fig:Comparison_cSGD_Darcy60}$ and $\ref{fig:Variances_cSGD_Darcy60}$) and quantitative (sub-figure (a) of Figure $\ref{fig:Errors_Darcy}$) evidence.

\begin{figure}
	\centering
	\includegraphics[ keepaspectratio=true, width=0.8\textwidth,  clip=true,  trim=120pt 2pt 120pt 30pt]{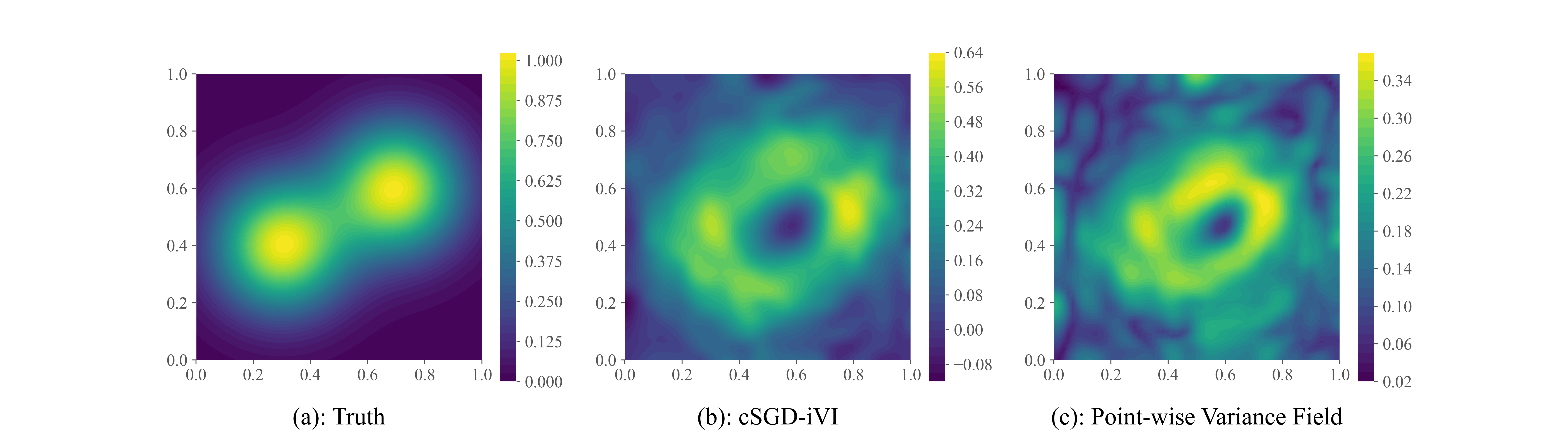}

	\caption{\emph{\small 
			(a): The background truth of $u$;
			(b): The estimated posterior mean function of $u$ obtained by cSGD-iVI;
			(c): The estimated variance function of posterior measure of $u$ obtained by cSGD-iVI.
	}}
	
	\label{fig:Comparison_cSGD_Darcy60}
\end{figure}

\begin{figure}
	\centering
	\includegraphics[ keepaspectratio=true, width=0.8\textwidth,  clip=true,  trim=120pt 2pt 120pt 30pt]{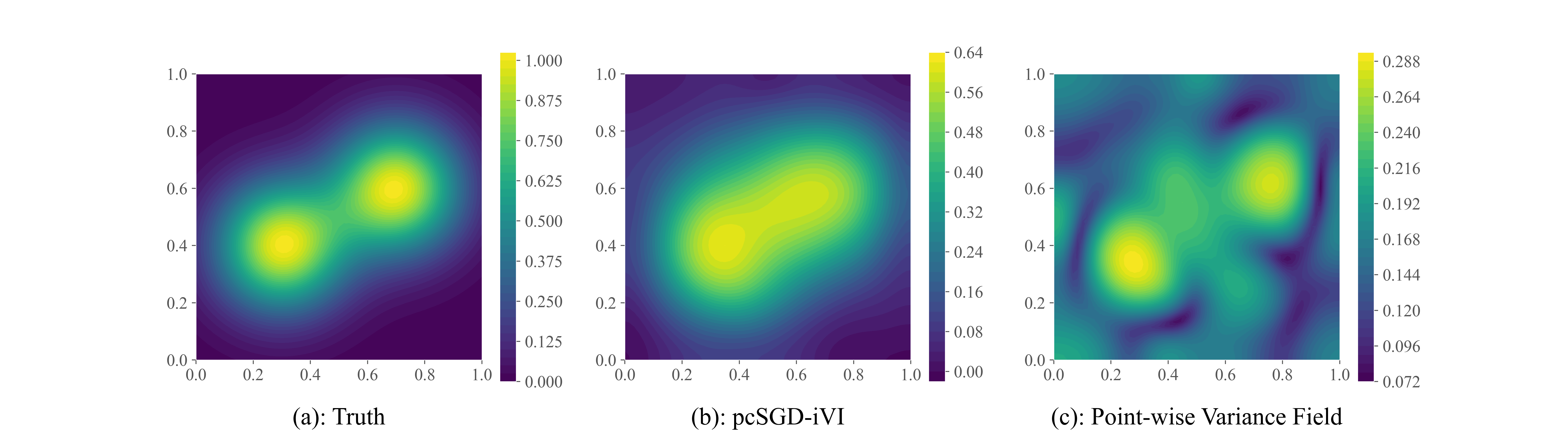}

	\caption{\emph{\small 
			(a): The background truth of $u$;
			(b): The estimated posterior mean function of $u$ obtained by pcSGD-iVI;
			(c): The estimated variance function of posterior measure of $u$ obtained by pcSGD-iVI.
	}}
	
	\label{fig:Comparison_pcSGD_Darcy60}
\end{figure}

\begin{figure}
	\centering
	\includegraphics[ keepaspectratio=true, width=0.6\textwidth,  clip=true, trim=80pt 2pt 80pt 30pt]{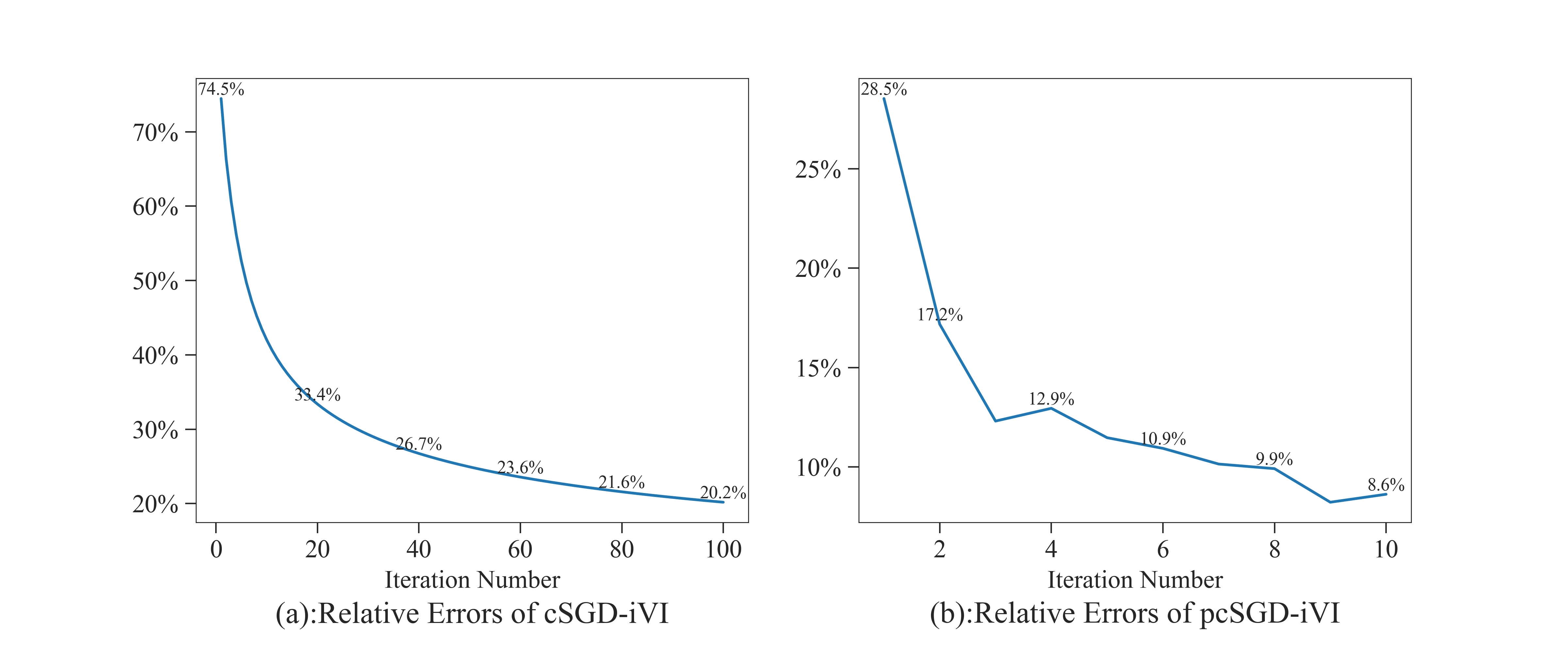}

	\caption{\emph{\small 
			(a): Relative errors of estimated posterior means of $u$ in the $L^2-norm$ of cSGD-iVI method with size of mesh $60 \times 60$; 
			(b): Relative errors of estimated posterior means of $u$ in the $L^2-norm$ of pcSGD-iVI method with size of mesh $60 \times 60$.
	}}
	
	\label{fig:Errors_Darcy60}
\end{figure}

\begin{figure}
	\centering
	\includegraphics[ keepaspectratio=true, width=0.8\textwidth,  clip=true, trim=120pt 2pt 120pt 30pt]{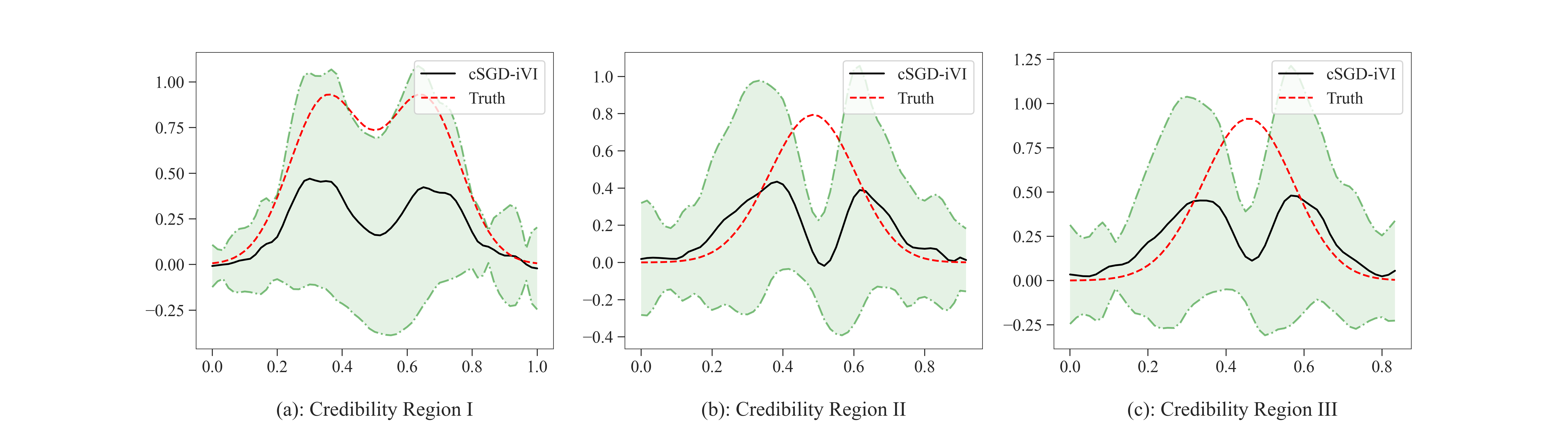}

	\caption{\emph{\small 
			The $95 \%$ credibility region of the estimated posterior mean function represented by the green shade area.
			(a): The black solid line represents the estimated mean $\lbrace u(x_i, x_i) \rbrace^{60}_{i=1}$, and red dashed line represents the true function $\lbrace u^{\dagger}(x_i, x_i) \rbrace^{60}_{i=1}$.
			(b): The black solid line represents the estimated mean $\lbrace u(x_i, x_{i+5}) \rbrace^{55}_{i=1}$, and red dashed line represents the true function $\lbrace u^{\dagger}(x_i, x_{i+5}) \rbrace^{55}_{i=1}$.
			(c): The black solid line represents the estimated mean $\lbrace u(x_i, x_{i+10}) \rbrace^{50}_{i=1}$, and red dashed line represents the true function $\lbrace u^{\dagger}(x_i, x_{i+10}) \rbrace^{50}_{i=1}$. 
	}}
	
	\label{fig:Variances_cSGD_Darcy60}
\end{figure}

\begin{figure}
	\centering
	\includegraphics[ keepaspectratio=true, width=0.8\textwidth,  clip=true,trim=120pt 2pt 120pt 30pt]{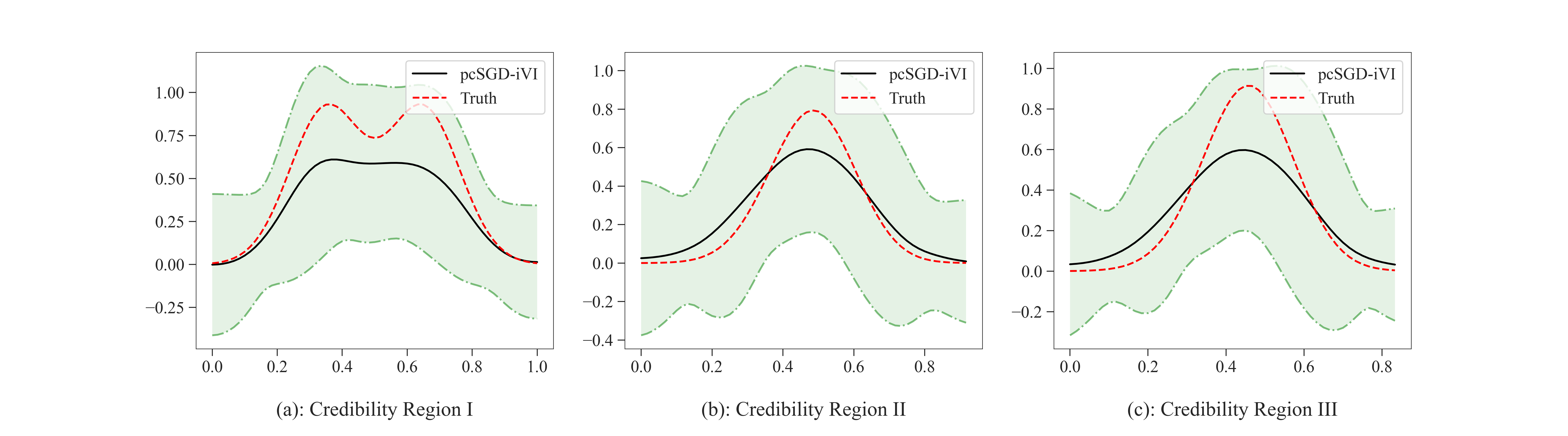}

	\caption{\emph{\small 
			The $95 \%$ credibility region of the estimated posterior mean function represented by the green shade area.
			(a): The black solid line represents the estimated mean $\lbrace u(x_i, x_i) \rbrace^{60}_{i=1}$, and red dashed line represents the true function $\lbrace u^{\dagger}(x_i, x_i) \rbrace^{60}_{i=1}$.
			(b): The black solid line represents the estimated mean $\lbrace u(x_i, x_{i+5}) \rbrace^{55}_{i=1}$, and red dashed line represents the true function $\lbrace u^{\dagger}(x_i, x_{i+5}) \rbrace^{55}_{i=1}$.
			(c): The black solid line represents the estimated mean $\lbrace u(x_i, x_{i+10}) \rbrace^{50}_{i=1}$, and red dashed line represents the true function $\lbrace u^{\dagger}(x_i, x_{i+10}) \rbrace^{50}_{i=1}$. 
	}}
	
	\label{fig:Variances_pcSGD_Darcy60}
\end{figure}

\subsubsection*{Computaional cost}

\begin{table}
	\renewcommand{\arraystretch}{1.5}
	\centering
	\caption{\emph{\small Computational cost of cSGD-iVI}} \label{table:Darcy_cost_cSGD}
	\begin{tabular}{c|ccc}
			\hline $\text{Mesh size} $ & $n=40\times 40$  & $n=50\times 50$ & $n=60\times 60$   \\
			\hline $\text{Time (total)/s}$ & $20.9677$  & $38.3116$ & $63.9130$  \\
			\hline $\text{Time (per-step)/s}$ & $0.20968$ & $0.38312$ & $0.63913$ \\
			\hline
		\end{tabular}
\end{table}

\begin{table}
	\renewcommand{\arraystretch}{1.5}
	\centering
	\caption{\emph{\small Computational cost of pcSGD-iVI}} \label{table:Darcy_cost_pcSGD}
	\begin{tabular}{c|ccc}
			\hline $\text{Mesh size} $ & $n=40\times 40$  & $n=50\times 50$ & $n=60\times 60$   \\
			\hline $\text{Time (total)/s}$ & $15.6127$  & $28.3549$ & $47.3421$  \\
			\hline $\text{Time (per-step)/s}$ & $1.56127$ & $1.89033$ & $3.15614$ \\
			\hline
		\end{tabular}
\end{table}

At last, we show the computational cost for both methods on different meshes.
In Table \ref{table:Darcy_cost_cSGD}, we show the computational cost of cSGD-iVI for \(100\) steps.
And in Table \ref{table:Darcy_cost_pcSGD}, we show the computational cost of pcSGD-iVI for \(15\) steps.
For each step, the proceeding time of pcSGD-iVI is longer then cSGD-iVI, since we need to calculate more PDEs. 
Overall, the cost of pcSGD-iVI method is more expensive than cSGD-iVI method, which is consistent with our discussion of computational cost in Section 3.1.2 of the manuscript.

}

\section{Proof details}\label{subsec:AppProof}
\subsection{Proof of Theorem 2.1}
~\\
\begin{proof}
    Following the assumptions of $\Phi(u)$, we easily know that the Assumption 1 and conditions of Theorems 15-16 in \cite{dashti2013bayesian} are satisfied, which are provided in Section D of supplementary material.
	Then the Bayes' formula $(2.2)$ is well-defined.
	Moreover, $Z_{\mu}$ is positive and finite.
\end{proof}


\subsection{Proof of Theorem 2.5}
~\\
\begin{proof}
	Based on the expression of estimated posterior covariance Eq.~\( (2.17) \), we obtain the following Rando-Nikodym derivative
	\[
	\frac{d\nu}{d\mu_0}(u) = \prod^M_{i=1}\sqrt{\frac{c_i}{s_i}} \exp \bigg(-\frac{1}{2}\| C^{-1/2}_\nu(u - \bar{u}) \|^2 +  \frac{1}{2}\| C^{-1/2}_0u \|^2 \bigg).
	\]
	By regarding prior \( \mu_0 \) as the reference measure, the KL divergence between $\nu$ and $\mu$ is rewritten as 
	\begin{align*}
		D_{\text{KL}}(\nu || \mu) &= \int_{\mathcal{H}_u} \log \bigg(\frac{d\nu}{d\mu}(u) \bigg)\nu(du)\\
		&= \int_{\mathcal{H}_u} \log \bigg(\frac{d\nu}{d\mu_0}(u) \bigg) - \log \bigg(\frac{d\mu}{d\mu_0}(u) \bigg) \nu(du)\\
		&= \int_{\mathcal{H}^M_u} \log \bigg(\frac{d\nu}{d\mu_0}(u) \bigg) - \log \bigg(\frac{d\mu}{d\mu_0}(u) \bigg) \nu^M(du) +  \int_{\mathcal{H}^{M\perp}_u}  - \log \bigg(\frac{d\mu}{d\mu_0}(u) \bigg) \nu^\perp(du),
	\end{align*} 
	where $\mathcal{H}^M_u := \text{span}\lbrace e_1, \cdots, e_M \rbrace$, and $\mathcal{H}^{M\perp}_u := \text{span}\lbrace e_{M+1}, e_{M+2}, \cdots \rbrace$.

    Noticing that the posterior measure on space \( \mathcal{H}_u^\perp \) is prior measure \( \mu_0^\perp \), the second term turns to 
    \begin{align*}
        \int_{\mathcal{H}^{\perp}_u}  - \log \bigg(\frac{d\mu}{d\mu_0}(u) \bigg) \nu^\perp(du) 
        = \frac{1}{2} \int_{\mathcal{H}^{\perp}_u} \Phi(u) \mu_0^\perp(du),
    \end{align*}
    which is a constant.
    
	The first term is then derived by
	\begin{align*}
		&\int_{\mathcal{H}^M_u} \log \bigg( \prod^M_{i=1}\sqrt{\frac{c_i}{s_i}} \exp \bigg(-\frac{1}{2}\| \mathcal{C}^{-1/2}_\nu(u - \bar{u}) \|^2 +  \frac{1}{2}\| \mathcal{C}^{-1/2}_0u \|^2 \bigg) \bigg) d\nu \\
		= &\int_{\mathcal{H}^M_u} \frac{1}{2} \bigg( \sum^{M}_{i=M}\log \frac{c_i}{s_i} - \| \mathcal{C}^{-1/2}_\nu(u - \bar{u}) \|^2 + \| \mathcal{C}^{-1/2}_0u \|^2  \bigg) + \Phi(u) d\nu \\
		= &\frac{1}{2} \bigg( \sum^M_{i = 1}\log\frac{c_i}{s_i} + \text{Tr}(\mathcal{C}_\nu \mathcal{C}^{-1}_0) + \text{Tr}(\mathcal{C}_\nu \mathcal{A}) \bigg) + \text{Const} \\
		=  &\frac{1}{2} \bigg( \sum^M_{i = 1}\log\frac{c_i}{s_i} +  \sum^M_{i = 1}\frac{s_i}{c_i} + \sum^M_{i = 1}s_ia_i \bigg) + \text{Const}
	\end{align*}
	
	Recalling the denotation \( \tilde{a}_i = a_i+\frac{1}{c_i}\) , by setting 
	\[
	\alpha_i = \frac{c_i}{s_i} = \frac{S^2(2\tilde{a}_i - \eta \tilde{a}^2_i)}{\eta q_i},
	\]
	the minimization problem becomes
	\begin{align*}
		\sum^M_{i = 1}\log\alpha_i +  \sum^M_{i = 1}\frac{1}{\alpha_i} + \sum^M_{i = 1}\frac{a_ic_i}{\alpha_i}{ +}\text{Const}.
	\end{align*}
	The derivative is written as:
	\[
	\sum^M_{i = 1}\frac{1}{\alpha_i} -  \sum^M_{i = 1}\frac{1+a_ic_i}{\alpha^2_i} = 0,
	\]
	which implies
	\[
	\sum^M_{i = 1}\alpha_i = \sum^M_{i = 1}\frac{S^2(2\tilde{a}_i - \eta \tilde{a}^2_i)}{\eta q_i} = M + \sum^M_{i = 1}a_ic_i.
	\]
	Then the optimal learning rate is given by
	\[
	\eta^\dagger = \frac{2S^2\sum^M_{i = 1}\frac{\tilde{a}_i}{q_i}}{M + \sum^M_{i = 1}a_ic_i + S^2\sum^M_{i = 1}\frac{\tilde{a}^2_i}{q_i}}
    := \frac{2S^2\sum^M_{i = 1}\frac{\tilde{a}_i}{q_i}}{M^\prime + S^2\sum^M_{i = 1}\frac{\tilde{a}^2_i}{q_i}}
	\]
    where \( \tilde{a}_i = a_i+\frac{1}{c_i} \).
\end{proof}

\subsection{Proof of Lemma 2.1}
~\\
\begin{proof}
	Eq.~$(2.22)$ can be rewritten as the following explicit form:
	\begin{align*}
		\mathbb{E}^{{\nu}_k} \bigg[ {u}_{(k+1)\eta} \bigg]= \sum^k_{j=1}\eta_j\prod^k_{i=j+1}(1-\eta_i{A}^{-1}){A}^{-1/2}\widetilde{{d}}.
	\end{align*}
	Following the definition in \cite{Peter2003IP1}, we have
	\begin{align*}
		\mathbb{E}^{{\nu}_k} \bigg[ {u}_{(k+1)\eta} \bigg] = g_{\alpha_k}({A}^{-1}){A}^{-1/2}\widetilde{{d}}
	\end{align*}
	where $g_{\alpha_k}(t) := \sum^k_{j=1}\eta_j\prod^k_{i=j+1}(1-\eta_it)$.
	
	The form of $r_{\alpha_k}$ is deduced by induction method.
	Firstly, setting $k=1$, $g_{\alpha_1}(t) = \eta_1$, and $r_{\alpha_k} = 1 - \eta_1t$.
	Assume $g_{\alpha_n}(t) = \sum^n_{j=1}\eta_j\prod^n_{i=j+1}(1-\eta_it)$ satisfying $r_{\alpha_n}(t) = 1 - tg_{\alpha_n}(t) = \prod^n_{i=1}(1 - \eta_it)$ with $k = n$.
	Then setting $k = n+1$, we have
	\begin{align*}
		g_{\alpha_{n+1}}(t) &= \sum^{n+1}_{j=1}\eta_j\prod^{n+1}_{i=j+1}(1-\eta_it) \\
		&= \sum^{n}_{j=1}\eta_j \bigg (\prod^{n}_{i=j+1}(1-\eta_it)(1 - \eta_{n+1}t) \bigg ) + \eta_{n+1} \\
		&= g_{\alpha_{n}}(t) \bigg (1 - \eta_{n+1}t \bigg ) + \eta_{n+1}.
	\end{align*}
	Following $r_{\alpha_{n+1}}(t) = 1 - tg_{\alpha_{n+1}}(t)$, we have
	\begin{align*}
		r_{\alpha_{n+1}}(t) &= 1 - \bigg (tg_{\alpha_{n}}(t)(1 - \eta_{n+1}t) + \eta_{n+1}t \bigg ) \\
		&= \bigg (1 - \eta_{n+1}t \bigg ) - tg_{\alpha_{n}}(t)(1 - \eta_{n+1}t) \\
		&= \bigg (1 - \eta_{n+1}t \bigg ) \bigg (1 - tg_{\alpha_{n}}(t) \bigg ) \\
		&= \prod^{n+1}_{i=1}(1 - \eta_it),
	\end{align*}
	which is the first conclusion.
	
	Recalling step sizes $\lbrace \eta_i \rbrace_i$ being non-increasing, we have $\eta_it < \eta_i / \eta_1 <1$.
	Then
	\begin{align*}
		| g_{\alpha_k}(t) | \leq | \sum^k_{j=1} \eta_j| =1/\alpha_k,
	\end{align*}
	which is the second conclusion.
	Similarly, we obtain the Third conclusion.
	
	Since 
	\begin{align*}
		\prod^k_{i=j+1}(1 - \eta_it)\rho(t) \leq \exp(-t\sum^k_{i=j+1}\eta_i)\rho(t),
	\end{align*}
	we have
	\begin{align*}
		e^{-t\sum^k_{i=j+1}\eta_i}\rho(t) \leq C_{\rho}\rho \bigg(\frac{1}{\sum^k_{i=j+1}\eta_i} \bigg).
	\end{align*}
	Applying $j=0$, this yields
	\begin{align*}
		\prod^k_{i=1}(1 - \eta_it)\rho(t) \leq C_{\rho}\rho \bigg(\frac{1}{\sum^k_{i=1}\eta_i} \bigg).
	\end{align*}
	Taking the absolute value of both sides of the equation, we obtain the last conclusion.
\end{proof}

\subsection{Proof of Theorem 2.7}
~\\
\begin{proof}
	 
    Taking the constant learning rate into Lemma 2.1, we have $\alpha_k = 1 / (k\eta^{\dagger})$.
	Then sequence $\mathbb{E}^{{\nu}_k} [u_{k+1}^M ]$ obeys Landweber iteration for operator ${A} = \mathcal{A} + \mathcal{C}_0^{-1}$.
	Based on Assumption \( { 2.6} \), function $\varphi^2$ is concave.
	With the help of Theorem 2 \cite{Peter2003IP1}, we obtain conclusion.
\end{proof}

\subsection{Proof of Theorem 2.8}
~\\
\begin{proof}
	Following the proof of Theorem 2.5, the minimization problem becomes
	\[
	\sum^M_{i = 1}\log\alpha_i +  \sum^M_{i = 1}\frac{1 + a_ic_i}{\alpha_i} + C,
	\]
	where
	\[
	\alpha_i = \frac{c_i}{s_i} = \frac{S^2(2\tilde{a}_i - \eta t_i\tilde{a}^2_i)}{\eta t_iq_i} 
	\]
	The derivative is written as:
	\[
	\sum^M_{i = 1}\frac{1}{\alpha_i} -  \sum^M_{i = 1}\frac{1 + a_ic_i}{\alpha^2_i} = 0,
	\]
	implying
	\[
	\sum^M_{i = 1}\alpha_i = \sum^M_{i = 1}\frac{S^2(2\tilde{a}_i - \eta t_i\tilde{a}^2_i)}{\eta t_iq_i} = M + \sum^M_{i = 1}a_ic_i.
	\]
	Then the optimal learning rate is given by
	\[
	\eta^\dagger = \frac{2S^2\sum^M_{i = 1}\frac{\tilde{a}_i}{t_iq_i}}{M + \sum^M_{i = 1}a_ic_i + S^2\sum^M_{i = 1}\frac{\tilde{a}^2_i}{q_i}}
    := \frac{2S^2\sum^M_{i = 1}\frac{\tilde{a}_i}{t_iq_i}}{M^\prime + S^2\sum^M_{i = 1}\frac{\tilde{a}^2_i}{q_i}}
	\]
    where \( \tilde{a}_i = a_i+\frac{1}{c_i} \).
\end{proof}

\subsection{Proof of Theorem 2.9}
~\\ 
\begin{proof}
	For simplicity, we denote
	\[
	g(u) = \nabla \tilde{\mathcal{L}}(x) = A^{1/2*}P^*PA^{1/2}(u - \bar{u}).
	\]
	
	Assume that the entries of \(P\) are i.i.d. with
	\[
	P_{ij}\sim \mathcal{N}\left(0,1/p\right).
	\]
	In order to calculate \( \mathbb{E}[g(u)] \), we aim to calculate \( \mathbb{E}[P^* P] \) first.
	The \((i,j)\)-th entry of the matrix \(P^T P \in \mathbb{R}^{N_d \times N_d}\) is:
	\[
	(P^T P)_{ij} = \sum_{k=1}^p P_{ki} P_{kj}.
	\]
	For \( i = j \):
	\[
	\mathbb{E}[(P^T P)_{ii}] = \mathbb{E} \left[ \sum_{k=1}^p P_{ki}^2 \right] 
	= \sum_{k=1}^p \mathbb{E}[P_{ki}^2] 
	= \sum_{k=1}^p \frac{1}{p} = 1,
	\]
	and others are zero.
	Then we obtain
	\[
	\mathbb{E}[P^T P] = I_{M},
	\]
	where \( I_{M}\) is the identity matrix.

%
	Then the expected gradient is
	\[
	\begin{aligned}
		\mathbb{E}[g(u)] &= \mathbb{E} \bigg[ A^{1/2*}P^*PA^{1/2}(u - \bar{u}) \bigg] \\
		&= A^{1/2*}\mathbb{E}[P^*P]A^{1/2}(u - \bar{u})  \\
		&=  A^{1/2*}\Bigl( I_{M}\Bigr)A^{1/2}(u - \bar{u}) \\
		&= A(u - \bar{u}).
	\end{aligned}
	\]
	
	We now compute the covariance of
	\[
	g(u) = A^{1/2*} P^* Pz,
	\]
	where $z = A^{1/2}( u - \bar{u})$.
	In the following derivation it is assumed that $P_{ij}\sim N(0,\sigma^2)$.
	Define the deviation of \(P^* P\) from its mean by
	\[
	\Delta = P^* P - \mathbb{E}[P^* P] = P^* P - p\,\sigma^2\, I_{M}.
	\]
	Then, we can write
	\[
	g(u) =  A^{1/2*}\Bigl(p\,\sigma^2\, I_{M} + \Delta\Bigr)z 
	= p\,\sigma^2\,  A^{1/2*}z +  A^{1/2*}(\Delta z).
	\]
	
	Since
	\[
	\mathbb{E}[g(u)] = p\,\sigma^2\,  A^{1/2*}z,
	\]
	the fluctuation is
	\[
	g(u) - \mathbb{E}[g(u)] = A^{1/2*}(\Delta z).
	\]
	Thus, the covariance is
	\[
	\operatorname{Cov}(g(u)) = \mathbb{E}\Bigl[\bigl(  A^{1/2*}(\Delta z)\bigr)\bigl( A^{1/2*}(\Delta z)\bigr)^*\Bigr]
	=  A^{1/2*}\,\mathbb{E}\Bigl[(\Delta z)(\Delta z )^*\Bigr]\, A^{1/2}.
	\]
	
	Next our goal is to calculate \(\mathbb{E}[(\Delta z)(\Delta z)^*]\).
	Let \(Q=P^* P = P^T P\). Then \(\Delta=Q-p\,\sigma^2\,I_{M}\). For a Wishart-distributed matrix \(Q \sim W_m(p,\,\sigma^2I_{M})\), we want to calculate \( \operatorname{Cov}(\Delta\Delta^*) \).
	Firstly, we denote that term \(p\,\sigma^2\,I_{M}\) is not random, then
	\[
	\operatorname{Cov}(\Delta\Delta^*) = \operatorname{Cov}(QQ^*).
	\]
	Recalling matrix \( Q = P^T P \) symmetric, we derive
	\[
	Q_{ij} = \sum_{a=1}^p P_{ai} P_{aj}, \quad Q_{rk} = \sum_{b=1}^p P_{br} P_{bk},
	\]
	where \( P_{ai} \) is the element of matrix \( P \) in the \( a \)-th row and \( i \)-th column. Each element \( P_{ai} \) follows a normal distribution with mean 0 and variance \( \sigma^2 \), and the elements of \( P \) are independent for different rows.
	Expand the product:
	\[
	Q_{ij} Q_{rk} = \left(\sum_{a=1}^p P_{ai} P_{aj}\right) \left(\sum_{b=1}^p P_{br} P_{bk}\right)
	\]
	Taking the expectation:
	\[
	\mathbb{E}[Q_{ij} Q_{rk}] = \mathbb{E}\left[\sum_{a=1}^p \sum_{b=1}^p P_{ai} P_{aj} P_{br} P_{bk}\right].
	\]
	Hence
	\begin{align*}
		\operatorname{Cov}(Q_{ij},Q_{rk})&=
		\mathbb{E}[Q_{ij} Q_{rk}]  - \mathbb{E}[Q_{ij}] \mathbb{E}[Q_{rk}]  \\
		&= \sum_{a=1}^p \sum_{b=1}^p \bigg(\mathbb{E}[P_{ai}P_{aj}P_{br}P_{bk} ] - \mathbb{E} [ P_{ai}P_{aj}]\mathbb{E}[P_{br}P_{bk} ]\bigg),
	\end{align*}
	The expectation of \(\mathbb{E}[P_{ai}P_{aj}P_{br}P_{bk} ] \) automatically accounts for both the covariance between the components of the vectors, and the relationship between the components across the two vectors \cite{bishop2006pattern}, then the covariance of the product of the random variables can be directly calculated from the fourth moments \(\mathbb{E}[P_{ai}P_{aj}P_{br}P_{bk} ]  \).
	Then the calculation can be simplified by
	\[
	\operatorname{Cov}(Q_{ij},Q_{rk}) =  \sum_{a=1}^p \sum_{b=1}^p \mathbb{E}[P_{ai}P_{aj}P_{br}P_{bk} ].
	\]
	
	Since the elements of \( P \) are independent, the expectation will be non-zero only when the indices match according to the covariance structure of the normal distribution:
	\begin{align*}
		\mathbb{E}[P_{ai} P_{br}] = \sigma^2 \delta_{ab} \delta_{ir}, \quad \mathbb{E}[P_{ai} P_{bk}] = \sigma^2 \delta_{ab} \delta_{ik}, \\
		\mathbb{E}[P_{aj} P_{br}] = \sigma^2 \delta_{ab} \delta_{jr}, \quad \mathbb{E}[P_{aj} P_{bk}] = \sigma^2 \delta_{ab} \delta_{jk}
	\end{align*}
	
	Hence:
	\[
	\mathbb{E}[Q_{ij} Q_{rk}] = \sigma^4\sum_{a=b=1}^p(\delta_{ab}\delta_{ir} \delta_{jk} + \delta_{ab}\delta_{ik} \delta_{jr})
	= p \sigma^4 (\delta_{ir} \delta_{jk} + \delta_{ik} \delta_{jr})
	\]
	Then we derive
	\[
	\operatorname{Cov}(\Delta_{ij},\Delta_{rk}) = p \sigma^4 (\delta_{ir} \delta_{jk} + \delta_{ik} \delta_{jr}).
	\]

	The \((i,r)\)-th entry of \( M=\mathbb{E}[(\Delta z)(\Delta z)^*] \) is given by
	\[
	\begin{aligned}
		M_{ir} &= \mathbb{E}\Biggl[\Bigl(\sum_{j=1}^{N_d} \Delta_{ij}\,z_j\Bigr)\Bigl(\sum_{k=1}^{N_d} \Delta_{rk}\,z_k\Bigr)\Biggr] \\
		&=\sum_{j,k=1}^{N_d} z_j z_k\,\operatorname{Cov}(\Delta_{ij},\Delta_{rk})\\
		&= p\,\sigma^4 \sum_{j,k=1}^{N_d} z_jz_k\Bigl(\delta_{ir}\delta_{jk}+\delta_{ik}\delta_{jr}\Bigr)\\[1mm]
		&= p\,\sigma^4\Bigl[\delta_{ir}\sum_{j=1}^{N_d} z_j^2+z_i z_r\Bigr]\\[1mm]
		&= p\,\sigma^4\Bigl(\|z\|^2\,\delta_{ir}+z_i z_r\Bigr).
	\end{aligned}
	\]
	In matrix notation,
	\[
	\mathbb{E}\Bigl[(\Delta z)(\Delta z)^*\Bigr] = p\,\sigma^4\Bigl(\|z\|^2 I_{M}+zz^*\Bigr).
	\]
	
	Substitute back into the covariance of \(g(x)\):
	\begin{align*}
		\operatorname{Cov}(g(u)) = A^{1/2*}\left[p\,\sigma^4\Bigl(\|z\|^2 I_{M}+zz^T\Bigr)\right]A^{1/2}.
	\end{align*}
	Taking \(\sigma^2=\frac{1}{p}\), we obtain
	\begin{align*}
		\operatorname{Cov}(g(u)) 
		&= p\left(\frac{1}{p^2}\right) A^{1/2*}\left[\Bigl(\|z\|^2 I_{M}+zz^T\Bigr)\right]A^{1/2}\\
		&= \frac{1}{p} A^{1/2*}\Bigl(\|A^{1/2}( u - \bar{u})\|^2 I_{M}+(A^{1/2}( u - \bar{u}))(A^{1/2}( u - \bar{u}))^*\Bigr)A^{1/2}.
	\end{align*}
	
	With \(P_{ij}\sim \mathcal{N}\left(0,\frac{1}{p}\right)\), the mean function and covariance operator of the stochastic gradient are
	\begin{itemize}
		\item Mean function:
		\[
		\mathbb{E}[\nabla \tilde{\mathcal{L}}(u)] = A(u - \bar{u})
		\]
		\item Covariance operator:
		\[
		\operatorname{Cov}(\nabla \tilde{\mathcal{L}}(u)) = \frac{1}{p} A^{1/2*}\Bigl(\|A^{1/2}( u - \bar{u})\|^2 I_{M}+(A^{1/2}( u - \bar{u}))(A^{1/2}( u - \bar{u}))^*\Bigr)A^{1/2}.
		\]
	\end{itemize}
\end{proof}

\section{Assumption 1 and Theorem 15-16 in \cite{dashti2013bayesian}}

Here we provide the Assumption 1 and Theorem 15-16 in \cite{dashti2013bayesian} needed in proving Theorem 2.1, which is shown in Assumption $\ref{assump1app}$, Theorem $\ref{theorem1app}$ and $\ref{theorem2app}$, respectively.
Here we need to state that the symbols in \cite{dashti2013bayesian} have been changed to the notations in our paper for the assumption and theorems below.

\begin{assumption}\label{assump1app}
	(Assumption 1 in \cite{dashti2013bayesian})
	Let us denote $X = \mathcal{H}_u, Y = \mathbb{R}^{N_d}$, and assume $\Phi \in C(X \times Y; \mathbb{R})$.
	Assume further that there are functions $M_i:\mathbb{R}^{+}\times \mathbb{R}^{+} \rightarrow \mathbb{R}^{+}, i=1, 2$, monotonic non-decreasing separately in each argument, and with $M_2$ strictly positive, such that for all $u \in \mathcal{H}_u$, $\bm{d}, \bm{d}_1, \bm{d}_2 \in B_{\mathbb{R}^{N_d}}(0, r)$,
	\begin{align*}
		\Phi(u; \bm{d}) &\geq -M_1(r, \lVert u \rVert_{\mathcal{H}_u}), \\
		\lvert \Phi(u; \bm{d}_1) - \Phi(u; \bm{d}_2)\rvert &\leq M_2(r, \lVert u \rVert_{\mathcal{H}_u})\lVert \bm{d}_1 - \bm{d}_2 \rVert_{\mathbb{R}^{N_d}}.
	\end{align*}
\end{assumption}

\begin{theorem}\label{theorem1app}
	(Theorem 15 in \cite{dashti2013bayesian})
	Let Assumption $\ref{assump1app}$ hold.
	Assume that $\mu_0(X) = 1$, and that $\mu_0(X \cap B) > 0$ for some bounded set $B$ in $X$.
	Assume additionally that, for every fixed $r>0$,
	\begin{align*}
		\exp(M_1(r, \lVert u \rVert_{\mathcal{H}_u})) \in L^1_{\mu_0}(X; \mathbb{R}),
	\end{align*}
	where $L_{\mu_0}^1(X;\mathbb{R})$ represents $X$-valued integrable functions under the measure $\mu_0$.
	Then for every $\bm{d} \in \mathbb{R}^{N_d}$, $Z_{\mu} = \int_X \exp(-\Phi(u ;\bm{d}))\mu_0(du)$ is positive and finite, and the probability measure $\mu$ is well defined, which is given by
	\begin{align*}
		\frac{d\mu}{d\mu_0}(u) = \frac{1}{Z_{\mu}}\exp(-\Phi(u; \bm{d})).
	\end{align*}
\end{theorem}

\begin{theorem}\label{theorem2app}
	(Theorem 16 in \cite{dashti2013bayesian})
	Let Assumption $\ref{assump1app}$ hold.
	Assume that $\mu_0(X) = 1$ and that $\mu_0(X \cap B) > 0$ for some bounded set $B$ in $X$.
	Assume additionally that, for every fixed $r > 0$,
	\begin{align*}
		\exp(M_1(r, \lVert u \rVert_{\mathcal{H}_u}))(1 + M_2(r, \lVert u \rVert_{\mathcal{H}_u})^2) \in L^1_{\mu_0}(X; \mathbb{R}).
	\end{align*}
	Then there is $C = C(r) > 0$ such that, for all $\bm{d}, \bm{d}^{\prime} \in B_Y(0, r)$
	\begin{align*}
		\sqrt{\frac{1}{2}\int \bigg (1 -  \sqrt{\frac{d\mu^{\prime}}{d\mu}} \bigg )^2d\mu} \leq C\lVert \bm{d} - \bm{d}^{\prime} \rVert_{\mathbb{R}^{N_d}},\\
	\end{align*}
	where measures $\mu^{\prime}$, $\mu$ represent the posterior measure according to $\bm{d}^{\prime}$, $\bm{d}$, respectively.
\end{theorem}

\bibliographystyle{amsplain}
\bibliography{references}